\documentclass[twoside]{article}

\usepackage[sc]{mathpazo} 
\usepackage[T1]{fontenc} 
\usepackage{microtype} 

\usepackage[english]{babel} 
\usepackage{subfig}
\usepackage{amsmath}
\usepackage{mathtools}
\usepackage{xcolor}
\usepackage[hmarginratio=1:1,top=32mm,columnsep=20pt]{geometry} 
\usepackage{booktabs}
\usepackage{algorithm}
\usepackage{algpseudocode}
\usepackage{enumitem} %
\usepackage[output-decimal-marker={,}]{siunitx}
\setlist[itemize]{noitemsep} 

\usepackage{abstract} 

\usepackage{titlesec} 

\titleformat{\section}[block]{\large\scshape\centering}{\thesection.}{1em}{} 
\titleformat{\subsection}[block]{\large \scshape}{\thesubsection.}{1em}{} 
\usepackage{graphicx}
\graphicspath{{images/}}
\newcommand{\bu}{\boldsymbol{u}}
\newcommand{\bw}{\boldsymbol{w}}
\newcommand{\bp}{\boldsymbol{p}}

\usepackage[usenames,dvipsnames]{xcolor}
\definecolor{mediumpurple}{rgb}{0.58, 0.0, 0.83}
\newcommand{\reviewerA}[1]{{\color{black}#1}}
\newcommand{\reviewerB}[1]{{\color{black}#1}}
\newcommand{\reviewerALL}[1]{{\color{black}#1}}

\usepackage{fancyhdr} 
\usepackage{multirow}
\usepackage{titling} 

\usepackage{hyperref} 
\usepackage{comment}

\usepackage{makecell}
\usepackage[bigfiles]{pdfbase}
\ExplSyntaxOn
\NewDocumentCommand\embedvideo{smm}{
  \group_begin:
  \leavevmode
  \tl_if_exist:cTF{file_\file_mdfive_hash:n{#3}}{
    \tl_set_eq:Nc\video{file_\file_mdfive_hash:n{#3}}
  }{
    \IfFileExists{#3}{}{\GenericError{}{File~`#3'~not~found}{}{}}
    \pbs_pdfobj:nnn{}{fstream}{{}{#3}}
    \pbs_pdfobj:nnn{}{dict}{
      /Type/Filespec/F~(#3)/UF~(#3)
      /EF~<</F~\pbs_pdflastobj:>>
    }
    \tl_set:Nx\video{\pbs_pdflastobj:}
    \tl_gset_eq:cN{file_\file_mdfive_hash:n{#3}}\video
  }
  \pbs_pdfobj:nnn{}{dict}{
    /Type/RichMediaInstance/Subtype/Video
    /Asset~\video
    /Params~<</FlashVars (
      source=#3&
      skin=SkinOverAllNoFullNoCaption.swf&
      skinAutoHide=true&
      skinBackgroundColor=0x5F5F5F&
      skinBackgroundAlpha=0.75
    )>>
  }
  \pbs_pdfobj:nnn{}{dict}{
    /Type/RichMediaConfiguration/Subtype/Video
    /Instances~[\pbs_pdflastobj:]
  }
  \pbs_pdfobj:nnn{}{dict}{
    /Type/RichMediaContent
    /Assets~<<
      /Names~[(#3)~\video]
    >>
    /Configurations~[\pbs_pdflastobj:]
  }
  \tl_set:Nx\rmcontent{\pbs_pdflastobj:}
  \pbs_pdfobj:nnn{}{dict}{
    /Activation~<<
      /Condition/\IfBooleanTF{#1}{PV}{XA}
      /Presentation~<</Style/Embedded>>
    >>
    /Deactivation~<</Condition/PI>>
  }
  \hbox_set:Nn\l_tmpa_box{#2}
  \tl_set:Nx\l_box_wd_tl{\dim_use:N\box_wd:N\l_tmpa_box}
  \tl_set:Nx\l_box_ht_tl{\dim_use:N\box_ht:N\l_tmpa_box}
  \tl_set:Nx\l_box_dp_tl{\dim_use:N\box_dp:N\l_tmpa_box}
  \pbs_pdfxform:nnnnn{1}{1}{}{}{\l_tmpa_box}
  \pbs_pdfannot:nnnn{\l_box_wd_tl}{\l_box_ht_tl}{\l_box_dp_tl}{
    /Subtype/RichMedia
    /BS~<</W~0/S/S>>
    /Contents~(embedded~video~file:#3)
    /NM~(rma:#3)
    /AP~<</N~\pbs_pdflastxform:>>
    /RichMediaSettings~\pbs_pdflastobj:
    /RichMediaContent~\rmcontent
  }
  \phantom{#2}
  \group_end:
}
\ExplSyntaxOff

\newcommand{\ddEF}{DD-EF}

\newcommand{\energyEFR}{\ensuremath{\mathcal{E}\text{-DD-EFR}}}

\newcommand{\enstrophyEFR}{\ensuremath{\mathcal{E}\mathcal{Z}\text{-DD-EFR}}}

\newcommand{\chienergy}{\ensuremath{\chi_{\mathcal{E}\text{-opt}}}}

\newcommand{\chienstrophy}{\ensuremath{\chi_{\mathcal{E} \mathcal{Z}\text{-opt}}}}

\newcommand{\chienstrophyonly}{\ensuremath{\chi_{\mathcal{Z}\text{-opt}}}}

\newcommand{\Ttrain}{\ensuremath{T_{\text{train}}}}

\newcommand{\Itrain}{\ensuremath{I_{\text{train}}}}

\newcommand{\Itest}{\ensuremath{I_{\text{test}}}}
\title{\Large\bfseries A new data-driven energy-stable Evolve-Filter-Relax model for turbulent flow simulation} 

\author{%
\normalsize \textsc{Anna Ivagnes}$^{\star}$\\[1ex] 
\small International School for Advanced Studies (SISSA), Via Bonomea 265, Trieste, 34136, Italy \\ 
\small \href{mailto:aivagnes@sissa.it}{aivagnes@sissa.it} 
\and
\normalsize \textsc{Toby van Gastelen}$^{\star}$\\[1ex] 
\small Scientific Computing, CWI, Science Park 123, Amsterdam, 1098XG, The Netherlands \\ 
\small \href{mailto:tvg@cwi.nl}{t.v.gastelen@cwi.nl} 
\and
\normalsize \textsc{Syver D{\o}ving Agdestein}\\[1ex] 
\small Scientific Computing, CWI, Science Park 123, Amsterdam, 1098XG, The Netherlands \\ 
\small \href{mailto:sda@cwi.nl}{s.d.agdestein@cwi.nl}
\and
\normalsize \textsc{Benjamin Sanderse}\\[1ex] 
\small Scientific Computing, CWI, Science Park 123, Amsterdam, 1098XG, The Netherlands \\ 
\small \href{mailto:b.sanderse@cwi.nl}{b.sanderse@cwi.nl} 
\and
\normalsize \textsc{Giovanni Stabile}\\[1ex] 
\small Biorobotics Institute, Sant'Anna School of Advanced Studies,\\
\small Viale Rinaldo Piaggio 34, 56025 Pontedera, Pisa, Italy \\ 
\small \href{mailto:giovanni.stabile@santannapisa.it}{giovanni.stabile@santannapisa.it} 
\and
\normalsize \textsc{Gianluigi Rozza}\\[1ex] 
\small International School for Advanced Studies (SISSA), Via Bonomea 265, Trieste, 34136, Italy \\ 
\small \href{mailto:gianluigi.rozza@sissa.it}{gianluigi.rozza@sissa.it} 
\and \small $^{\star}$ These authors contributed equally to this work.
}

\date{} 

\begin{document}

\maketitle
\begin{abstract}
We present a novel approach to define the filter and relax steps in the evolve-filter-relax (EFR) framework for simulating turbulent flows.
The EFR main advantages are its ease of implementation and computational efficiency.
However, as it only contains two parameters (one for the filter step and one for the relax step) its flexibility is rather limited.
In this work, we propose a data-driven approach in which the optimal filter is found based on DNS data in the frequency domain. 
The optimization step is computationally efficient and only involves one-dimensional least-squares problems for each wavenumber.
Across both decaying turbulence and Kolmogorov flow, our learned filter decisively outperforms the standard differential filter and the Smagorinsky model, yielding significantly improved accuracy in energy spectra and in the temporal evolution of both energy and enstrophy.
In addition, the relax parameter is determined by requiring energy and/or enstrophy conservation, which enforces stability of the method and reduces the appearance of numerical wiggles, especially when the filter is built in scarce data regimes.
Applying the learned filter is also more computationally efficient compared to traditional differential filters, as it circumvents solving a linear system.

\begin{figure}
    \centering
    \includegraphics[width=1\linewidth]{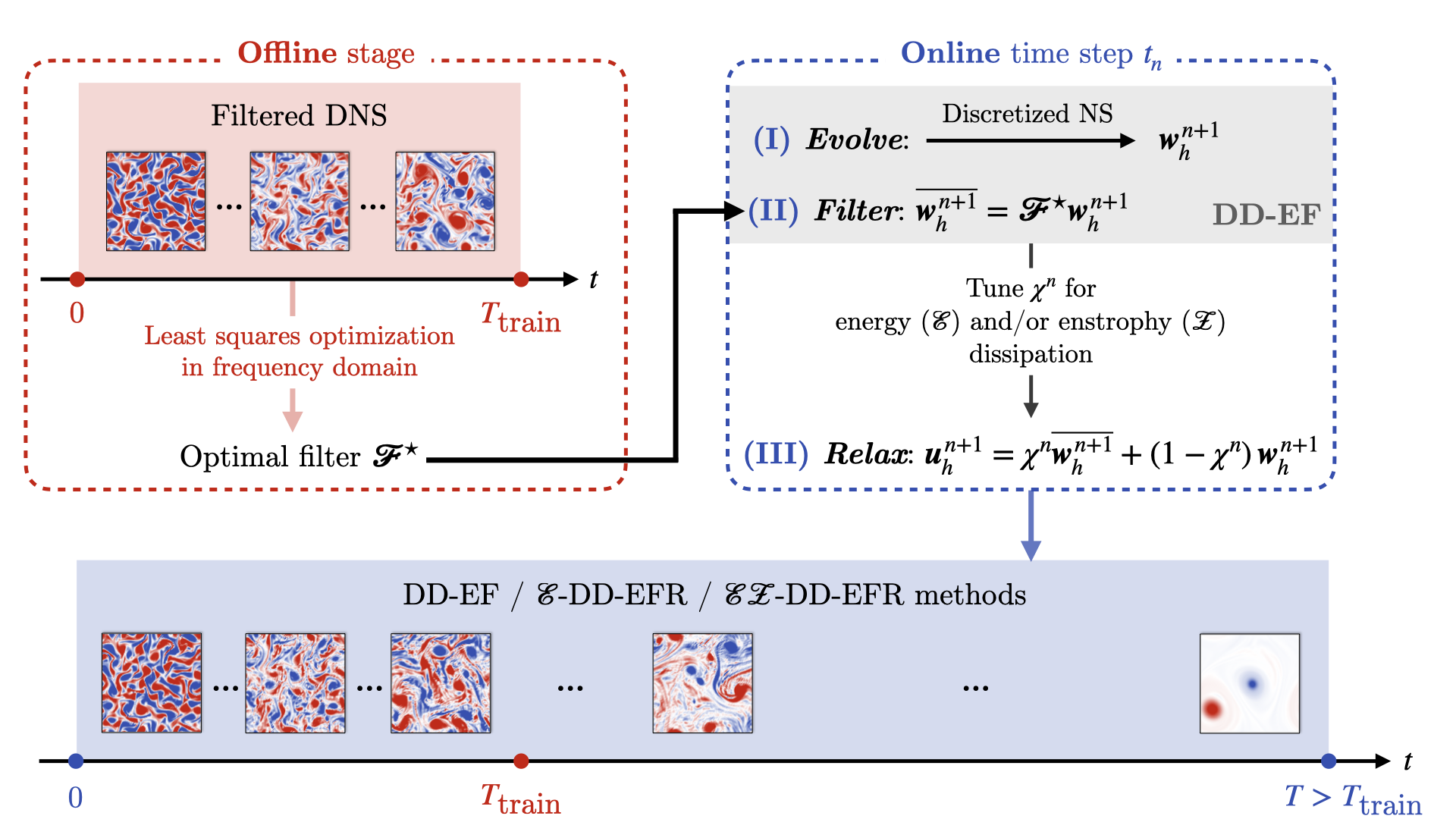}
    \caption{Workflow of the DD-EFR framework presented in the article. An optimal filter $\boldsymbol{\mathcal{F}}^{\star}$ is learned \emph{offline} from filtered DNS data. The relax parameter is tuned \emph{online} at each time step to control energy and/or enstrophy dissipation.}
    \label{fig:flowchart}
\end{figure}
\end{abstract}

\section{Introduction}
\label{sec:introduction}
It is well known that modeling turbulent flows poses significant challenges for traditional numerical discretization methods due to the emergence of turbulent eddies~\cite{roache1998verification}. In incompressible flows, this complexity is primarily governed by the Reynolds number, which dictates the scale of the smallest eddies present in the flow~\cite{pope2000turbulence}. Accurately capturing these features requires a computational grid fine enough to resolve all relevant scales. As the Reynolds number increases, so does the need for finer resolution, leading to a substantial rise in computational cost. This often renders \textbf{direct numerical simulations (DNS)} impractical for high Reynolds number flows~\cite{moin1998DNS}. \reviewerB{
To address this, various modeling strategies have been developed. These include \textbf{Reynolds-Averaged Navier–Stokes (RANS)} methods \cite{wilcox2006turbulence}, which model all turbulent fluctuations through additional transport equations for the averaged quantities, and \textbf{Large Eddy Simulation (LES)} \cite{sagaut2002LES}, which resolves the large, energy-containing eddies while modeling only the smaller, subgrid-scale motions. The two approaches thus differ not only in the degree of resolution but also in their underlying treatment of turbulence and computational requirements. Beyond turbulence modeling, \textbf{numerical stabilization techniques} \cite{brooks1982upwind, collis2024analysissupgmethodsolution}, such as the streamline upwind/Petrov–Galerkin (SUPG) method, are often employed to mitigate numerical oscillations in convection-dominated or under-resolved flows. This improves stability without explicitly modeling turbulence.}

One of the most popular numerical stabilization techniques is the \textbf{Evolve-Filter-Relax (EFR)} framework \cite{quaini2024bridginglargeeddysimulation,ervin2012EFR}. It is a widely used approach for stabilizing numerical simulations by suppressing spurious oscillations. It has proven its worth in a multitude of applications such as nuclear engineering design, simulation of ocean circulation, and cardiovascular modeling \cite{girfoglio2019finite, girfoglio2023novel, besabe2025linear, xu2018coupled, quaini2024bridginglargeeddysimulation, xu2020backflow,  tsai2023parametric, layton2012approximate}. It consists of three steps: (1) evolving the solution using a coarse spatial discretization, (2) applying a differential filter to the evolved state to dampen high-frequency numerical noise, and (3) computing a convex combination of the evolved and filtered states to mitigate excessive diffusion introduced by the filter. This process is repeated at each time step to maintain stability and accuracy throughout the simulation. Key advantages of EFR are its computational efficiency and ease of implementation. It only requires access to a discrete diffusion operator to construct the filter, making it compatible with a variety of numerical methods, including finite volumes, finite elements, and \textbf{Reduced Order Models (ROMs)} \cite{benner2015ROM, canuto1988some, fischer2001filter, germano1986differential, girfoglio2019finite,  girfoglio2021fluid, xie2018evolve, tsai2025time}. In ROM applications, the filter is typically projected onto the reduced basis. In this work, however, we focus on the application of the EFR framework in LES context \cite{quaini2024bridginglargeeddysimulation,girfoglio2019finite}, where the differential filter is constructed using a coarsened or modified diffusion operator. This allows the method to effectively emulate subgrid-scale dissipation while remaining flexible and easy to integrate into existing simulation codes.

Applying the differential filter involves solving a linear system, where the filtered field is obtained by applying the inverse of a differential operator to the original field. In \cite{quaini2024bridginglargeeddysimulation} it is shown that solving an elliptic filter equation is equivalent to applying an eddy-viscosity model in the LES context using a backward Euler scheme. 
In EFR context different eddy-viscosity models are referred to as indicator functions. These functions act as spatially-varying weights that determine where and to what extent the filter is applied \cite{berselli2006mathematics}. Furthermore, they are designed to identify regions within the domain with strong gradients or under-resolved features, such as near sharp interfaces or turbulent structures. This allows the differential filter to act locally rather than globally. This adaptive approach improves the accuracy and efficiency of the EFR method, making it suitable for complex flow simulations.
Moreover, although it is sparse, the filter step still requires solving a linear system to apply the differential filter, which is computationally expensive. 

The EFR framework only contains two parameters, namely the filter radius $\delta$ and the relaxation parameter $\chi$ \cite{quaini2024bridginglargeeddysimulation,ervin2012EFR}. Choosing the values for these parameters requires careful consideration. If the filter radius is too large the solution fields become too diffused, while a too small value might fail in adequately suppressing the noise. Different ways of choosing these parameters have been suggested: e.g. setting the filter radius either equal to the mesh size \cite{girfoglio2019finite,strazzullo2022consistency}, or equal to the Kolmogorov lengthscale \cite{bertagna2016deconvolution}, and setting the relaxation parameter proportional to the size of the time step \cite{strazzullo2022consistency,bertagna2016deconvolution}.
\reviewerB{The sensitivity of regularized and Leray-type models to the choice of the filter radius has also been analyzed in detail \cite{bertagna2018sensitivity}, while the impact of filtering parameters in patient-specific LES simulations has been investigated through global sensitivity analysis \cite{xu2021global}.}

However, these approaches for choosing the parameters do not fully alleviate the limited flexibility of the model, as highlighted by recent results in \cite{ivagnes2025datadrivenoptimizationevolvefilterrelaxregularization}. 
In that work, it was shown that the optimal values of the parameters are highly variable during the course of a simulation. By introducing 
a data-driven approach to determine the optimal parameters at each time step, the results more closely resembled the DNS results. However, extrapolation capabilities of such an approach are limited, as DNS data for each time step is required to determine the optimal parameters settings.
Making use of different indicator functions has been shown to improve the performance of the EFR framework. Different indicator functions are available, which effectively make the differential filter non-linear \cite{bertagna2016deconvolution}. An overview of such functions can be found here \cite{quaini2024bridginglargeeddysimulation} specialized to different use cases.

While traditional EFR approaches rely on \reviewerB{physically or mathematically motivated} indicator functions to control dissipation and filter strength, recent years have seen growing interest in data-driven alternatives. These data-driven approaches are especially researched in the LES community. The LES and EFR frameworks are related in the fact that both aim to suppress spurious oscillations during simulation time and do this in a similar manner \cite{quaini2024bridginglargeeddysimulation}. In particular, machine learning approaches have come forward as a way to improve the accuracy of LES simulations \cite{beck2022MLLES_1,kurz2023deep,shankar2023differentiable,maulik2018SGSmodeling2D,kochkov2021machine,park2021MLLES,Boral2023MLLES,guan2022MLLES,Agdestein2025MLLES}. An example of this is the work in \cite{List2022MLclosure}, where a convolutional neural network was added to the right hand side of the system to improve LES simulation. It was found that the accuracy of the resulting models relies heavily on the way the neural network parameters are optimized. Using a differentiable solver, gradient-based optimization techniques can be applied to optimize these parameters with respect to how well the coarsened DNS solution is reproduced. This is often required, as training neural networks to solely reproduce the right hand side, although requiring significantly less computational effort, often displays inaccuracies and instabilities \cite{List2022MLclosure,frezat2021MLLES,macart2021MLLES,syver2022MLLES,melchers2022MLLES}. The clear downside of the machine learning approach is the required computational effort to generate the training data and train the neural networks, and the required access to a differentiable solver. The latter makes the approach incompatible with many of the existing codebases and requires specialized solvers such as \cite{incompressibleNSElibrary}. Furthermore, the neural networks are less interpretable than the differentiable filter used in the EFR framework.

In this work we aim to combine the increased flexibility of data-driven approaches with the modularity and computational efficiency of the EFR framework. We therefore propose the following idea: \textit{learn the optimal data-driven linear filter} which is trained on DNS data. This allows us to construct a filter which is tailored to the problem at hand.
We combine the learned filter with an energy-conserving finite volume discretization such that we can analyze the stability of the method.
Next to being more flexible than the differentiable filter this approach has some other key advantages: Making use of a \textbf{fast Fourier transform (FFT)} the filter can be applied in a computationally efficient manner and no longer involves solving a linear system. Furthermore, computing the optimal filter coefficients involves only the solution of one-dimensional least-squares problems for each coefficient. This makes finding the optimal coefficients much more efficient than training a neural network \cite{kingma2017adammethodstochasticoptimization}. 
Furthermore, we introduce a clever way to use the relax parameter $\chi$ to preserve both energy and enstrophy. This is done with the goal of ensuring both accuracy and stability when little DNS data is available.

The paper is structured as follows: In Section \ref{sec:preliminaries} we start off by introducing the incompressible Navier-Stokes equations, discuss the used spatial discretization, and introduce the EFR framework. These are regarded as \emph{preliminaries}. Furthermore, in Section \ref{sec:methods} we introduce our data-driven methodology and the way we impose the energy and enstrophy constraints. After this we present our results in Section \ref{sec:results}, where we compare to existing approaches. We conclude our work in Section \ref{sec:conclusions}.

\section{Preliminaries}
\label{sec:preliminaries}

In this Section we introduce the Navier-Stokes equations and describe the discretization in space and time. Moreover, we introduce the \textbf{Evolve-Filter-Relax (EFR)} framework, to account for coarse grids. 

\subsection{The Navier-Stokes equations}
\label{subsec:NS}

We model the dynamics of an incompressible fluid using the \textbf{incompressible Navier--Stokes Equations (NSE)}. Let $\Omega \subset \mathbb{R}^d$ be a fixed spatial domain, where $d = 2$ or $3$. The fluid velocity is denoted by $\boldsymbol{u} \doteq \boldsymbol{u}(\boldsymbol{x}, t) \in \mathbb{U}$, and the pressure by $p \doteq p(\boldsymbol{x}, t) \in \mathbb{Q}$. The governing equations for the fluid motion are:
\begin{equation}
\begin{cases}\label{eq:NSE}
 \displaystyle \frac{\partial \boldsymbol{u}}{\partial t} + (\boldsymbol{u} \cdot \nabla) \boldsymbol{u} - \nu \Delta \boldsymbol{u} + \nabla p = f & \text{in }  \Omega \times (t_0, T), \\
\nabla \cdot \boldsymbol{u} = 0  & \text{in }  \Omega \times (t_0, T), 
\end{cases}
\end{equation}
with the initial condition $\boldsymbol{u} = \boldsymbol{u}_0$ in $\Omega \times { t_0 }$. Here, $f$ denotes external forcing, $\nu$ is the kinematic viscosity, and $\mathbb{U}$ and $\mathbb{Q}$ are appropriate Hilbert spaces for velocity and pressure, respectively. The simulation runs over the time interval $(t_0, T)$.

In this study, we adopt periodic boundary conditions and therefore omit an explicit discussion of boundary terms. The incompressibility condition, $\nabla \cdot \boldsymbol{u} = 0$, enforces mass conservation.
The nature of the flow is characterized by the dimensionless Reynolds number:
\begin{equation}
\label{eq:Re}
Re \doteq \frac{UL}{\nu},
\end{equation}
where $U$ and $L$ are the characteristic velocity and the length scale of the system, respectively. A high Reynolds number indicates that inertial effects outweigh viscous effects, leading to a convection-dominated regime.

\subsection{Finite volume discretization}
\label{subsec:FV_discretization}

In the numerical tests, we employ a \textbf{Finite Volume (FV)} discretization on a staggered grid presented in \cite{benjamin_thesis,harlow1965numerical},
and an explicit fourth-order Runge-Kutta time integration scheme. 
Using this discretization the velocity field is approximated on the cell-faces and is discretely represented by the vector $\bu_h \in \mathbb{R}^{N_u}$, whereas the pressure field is approximated in the cell-centers as $\bp_h \in \mathbb{R}^{N_p}$. $N_u$ and $N_p$ are the number of discrete points for the velocity field and the pressure, respectively.
Using the FV discretization \eqref{eq:NSE} is written in matrix representation as
\begin{equation}\label{eq:momentum-matrix}
    \frac{d \bu_h}{dt} = - \boldsymbol{C}_h(\bu_h)\bu_h+ \nu \boldsymbol{D}_h \bu_h + \boldsymbol{f}_h - \boldsymbol{G}_h \bp_h,
\end{equation}
where $\boldsymbol{C}_h(\bu_h) \in \mathbb{R}^{N_u \times N_u}$ represents the nonlinear convection operator, $\boldsymbol{D}_h \in \mathbb{R}^{N_u \times N_u}$ the linear diffusion operator, and $\boldsymbol{f}_h  \in \mathbb{R}^{N_u}$ the external forces. Moreover, $\boldsymbol{G}_h \in \mathbb{R}^{N_u \times N_p}$ is the pressure gradient operator.
In matrix representation the divergence freeness condition is written as
\begin{equation}
\boldsymbol{M}_h \boldsymbol{u}_h = \boldsymbol{0}_h,
    \label{eq:div-free-matrix}
\end{equation}
where $\boldsymbol{M}_h \in \mathbb{R}^{N_p \times N_u}$ represents the divergence operator and $\boldsymbol{0}_h \in \mathbb{R}^{N_p}$ represents a vector of zeros. 
The staggered grid discretization is chosen as it preserves the kinetic energy in the inviscid limit. To see this, we start by defining the (kinetic) energy as
\begin{equation}\label{eq:continuous_E}
    \mathcal{E} := \frac{1}{2 |\Omega|}\int_\Omega \bu \cdot \bu \, \text{d}\Omega.
\end{equation}
Discretely, we represent this as
\begin{equation}\label{eq:discrete_E}
    \mathcal{E}_h := \frac{h^2}{2|\Omega|} \lVert \bu_h \lVert_2^2,
\end{equation}
\reviewerB{where
\begin{equation}
    \lVert \bu_h \lVert_2^2 := \bu_h^T \bu_h,
\end{equation}
assuming a uniform 2D grid with grid spacing $h$ in each direction, as used in this work. The appearance of $h^2$ is such that the unit of the discrete energy is the same as its continuous counterpart, i.e. \eqref{eq:discrete_E} discretely approximates the integral in \eqref{eq:continuous_E}.}
The discretization in \eqref{eq:momentum-matrix} produces the following energy behavior:
\begin{equation}\label{eq:energy_behavior}
    \frac{\text{d}\mathcal{E}_h}{\text{d}t} = \frac{h^2}{2|\Omega|} \bu_h^T \frac{\text{d}\bu_h}{\text{d}t} = \frac{h^2}{2|\Omega|} (- \nu \lVert \boldsymbol{Q}_h \bu_h \lVert_2^2 + \bu_h^T\boldsymbol{f}_h).
\end{equation}
The energy behavior in \eqref{eq:energy_behavior} is derived by using the product rule, starting from \eqref{eq:discrete_E}, and using the fact that $\boldsymbol{D}_h$ can be Cholesky decomposed as $\boldsymbol{D}_h := - \boldsymbol{Q}_h  \boldsymbol{Q}_h$, while the pressure gradient and convective term contributions cancel, see \cite{benjamin_thesis,harlow1965numerical}. From \eqref{eq:energy_behavior} it is clear that this discretization provides stability, as in the absence of forcing the energy can never increase. Note that this energy behavior is true for periodic boundary conditions. In the case of boundary conditions of different type, boundary contributions will appear in the energy equation.

\subsection{Evolve-Filter-Relax}\label{subsec:efr}

In under-resolved or marginally-resolved simulations central spatial discretization schemes, such as the one adopted in this work, tend to produce spurious oscillations in convection-dominated regimes.
To mitigate these artifacts, we employ the \textbf{Evolve-Filter-Relax (EFR)} algorithm. Let the time step be denoted by $\Delta t$, with $t_n = t_0 + n\Delta t$ for $n = 0, \dots, N$, and total simulation time $T = t_0 + N\Delta t$. We write $y^n$ to denote the numerical approximation of a generic quantity $y$ at time $t_n$.

The EFR procedure consists of the following three steps:
\begin{eqnarray}
&\text{\bf (I)} & \text{\emph{Evolve:}} \quad 
\begin{cases}
\displaystyle \frac{\bw_h^{n + 1} - \bu_h^n}{\Delta t} = - \boldsymbol{C}_h(\bu_h^n)\bu_h^n + \nu \boldsymbol{D}_h \bu_h^n + \boldsymbol{f}_h - \boldsymbol{G}_h \boldsymbol{p}_h^n, \\
\boldsymbol{M}_h \bw_h^{n+1} = \boldsymbol{0}_h,
\end{cases}
\label{eqn:ef-rom-1} \nonumber \\[0.3cm]
&\text{\bf (II)} & \text{\emph{Filter:}} \quad 
(\boldsymbol{I} - 2\delta^2 \boldsymbol{D}_h) \overline{\bw}_h^{n+1} = \bw_h^{n+1} 
\label{eqn:ef-rom-2} \nonumber \\[0.3cm]
&\text{\bf (III)} & \text{\emph{Relax:}} \quad 
\bu_h^{n+1} = (1 - \chi)\bw_h^{n+1} + \chi \, \overline{\bw}_h^{n+1}
\label{eqn:ef-rom-3} \nonumber
\end{eqnarray}
Here, $\bw_h^{n+1}$ is the velocity field obtained from the evolve step, $\overline{\bw}_h^{n+1}$ is its filtered counterpart, and $\bu_h^{n+1}$ is the final relaxed velocity field. The parameters $\delta$ and $\chi$ denote the \textit{filter radius} and \textit{relaxation parameter}, respectively, with $\chi \in [0,1]$.

In step \textbf{(I)}, the velocity field is advanced in time using a standard temporal discretization of the NSE. For simplicity, a forward Euler method is shown here, though in practice higher-order schemes are commonly used to enhance accuracy and stability~\cite{Butcher2007RK4,benjamin_thesis}.

Step \textbf{(II)} applies a \textit{differential filter} to smooth the velocity field. The filter is defined via an elliptic operator involving the discrete Laplacian $\boldsymbol{D}_h$ and an explicit spatial scale $\delta$. This step removes small-scale (high-frequency) components and has been shown to be mathematically robust~\cite{berselli2006mathematics}. The filtered field $\overline{\bw}_h^{n+1}$ is obtained by solving a linear system.

Step \textbf{(III)} combines the evolved and filtered velocities using a convex combination controlled by the relaxation parameter $\chi$. The extremes of this procedure are:
\begin{itemize}
    \item[(i)] $\chi = 1$: the fully filtered velocity is used, i.e., $\bu_h^{n+1} = \overline{\bw}_h^{n+1}$. We will refer to this approach as EF;
    \item[(ii)] $\chi = 0$: no filtering is applied, and $\bu_h^{n+1} = \bw_h^{n+1}$, and we refer to this as noEFR.
\end{itemize}

Thus, the parameter $\chi$ modulates the influence of filtering, allowing one to tune the balance between stability and accuracy~\cite{ervin2012numerical, fischer2001filter, mullen1999filtering}. This is particularly useful when filtering leads to excessive numerical diffusion, as discussed in~\cite{bertagna2016deconvolution}.

The filter radius $\delta$ is typically chosen in relation to either the minimum grid size $h_{\min}$ or the \textit{Kolmogorov length scale} $\eta = L\, Re^{-3/4}$, and usually $\eta \simeq h_{\min}$ for well-resolved DNS simulations \cite{strazzullo2022consistency, bertagna2016deconvolution}.

\subsubsection{Stability of the differential filter}\label{subsec:stability_diff_filter}

In Section \ref{subsec:FV_discretization} we discussed the stability of the spatial discretization. In this Section we do the same for the EFR algorithm. In order to derive stability of the differential filter we consider the relaxed velocity field at time step $n+1$ as
\begin{equation}\label{eq:relax-new}
    \bu^{n+1}_h = (1-\chi)\bw_h^{n+1} + \chi \boldsymbol{\mathcal{F}} \bw_h^{n+1},
\end{equation}
where $\boldsymbol{\mathcal{F}} = (\boldsymbol{I}- 2\delta^2 \boldsymbol{D}_h)^{-1}$ represents the filter such that $\overline{\bw_{h}}^{n+1} = \boldsymbol{\mathcal{F}}  \bw_h^{n+1}$. In order to prove stability we aim to show the following:
\begin{equation}\label{eq:relax_inequality}
\lVert \bu_h^{n+1}\lVert _2^2 \leq \lVert \bw_h^{n+1}\lVert _2^2,
\end{equation}
i.e. the relax step dissipates energy from the system. 
Since the relax parameter satisfies $0 \leq \chi \leq 1$, expression \eqref{eq:relax-new} is a convex combination of vectors, and the following holds \cite{convex2009linearalgebra}:
\begin{equation}
    \min(\Vert \bw_h^{n+1} \Vert_2^2,\Vert \boldsymbol{\mathcal{F}} \bw_h^{n+1} \lVert_2^2) \leq \Vert \bu^{n+1}_h \lVert_2^2 \leq \max(\Vert \bw_h^{n+1} \lVert_2^2,\Vert \boldsymbol{\mathcal{F}} \bw_h^{n+1} \lVert_2^2).
\end{equation}
From this we find that 
\begin{equation}\label{eq:stability_inequality}
    \lVert {\boldsymbol{\mathcal{F}}}\bw_h\lVert _2^2 \leq \lVert \bw_h \lVert _2^2
\end{equation}
is the required condition for the filter to be dissipative, where we omitted the superscripts $n$ and $n+1$ for ease of notation. 

To show this for the differential filter we express the diffusion operator $\boldsymbol{D}_h$ in terms of its eigenvalues $\boldsymbol{\Lambda}_{ii} = \lambda_i$ and orthonormal eigenvectors $P$ and use this to rewrite the filter:
\begin{equation}
\begin{split}
\boldsymbol{\mathcal{F}}  &= (\boldsymbol{I}-2\delta^2 \boldsymbol{D}_h)^{-1} = (\boldsymbol{P} \boldsymbol{I} \boldsymbol{P}^\dag -2\delta^2 \boldsymbol{P} \boldsymbol{\Lambda} \boldsymbol{P}^\dag )^{-1} = (\boldsymbol{P} (\boldsymbol{I}-2\delta^2 \boldsymbol{\Lambda}) \boldsymbol{P}^\dag )^{-1} \\ &= (\boldsymbol{P}^\dag )^{-1} (\boldsymbol{I}-2\delta^2 \boldsymbol{\Lambda})^{-1} \boldsymbol{P}^{-1} = \boldsymbol{P} (\boldsymbol{I}-2\delta^2 \boldsymbol{\Lambda})^{-1} \boldsymbol{P}^\dag ,
\end{split}
\label{eq:diff-filter-ref}
\end{equation}
where the eigenvectors form an orthogonal basis due to the symmetry of $\boldsymbol{D}_h$, i.e. $\boldsymbol{P}^\dag \boldsymbol{P} = \boldsymbol{P}\boldsymbol{P}^\dag  =\boldsymbol{I}$, with the superscript $\dag$ indicating the Hermitian conjugate. Note that for periodic boundary conditions and uniform grids the eigenvectors in $\boldsymbol{P}$ simply correspond to the Fourier modes \cite{Olson2014circulant}. Next, we express the filtering procedure in eigenvector space: 
\begin{equation}\label{eq:relax_eigenvector}
    \begin{split}
        \boldsymbol{\mathcal{F}}  \bw_h =  \boldsymbol{P} (\boldsymbol{I}-2\delta^2 \boldsymbol{\Lambda})^{-1} \boldsymbol{P}^\dag  \bw_h. 
    \end{split}
\end{equation}
The diagonal entries of this matrix, which correspond to the eigenvalues of $\boldsymbol{\mathcal{F}}$, are given by 
\begin{equation}\label{eq:definition_f}
    \hat{\mathcal{F}}_{ii} = \hat{f}_i :=   \frac{1}{1-2\delta^2\lambda_{i}} \in \mathbb{C}, \quad i=1, \dots, N_u.
\end{equation}
Next, we use \eqref{eq:relax_eigenvector} to rewrite the stability inequality \eqref{eq:stability_inequality}:
\begin{equation}\label{eq:deriv_stability_req}
\begin{split}
     \lVert \boldsymbol{P} \hat{\boldsymbol{\mathcal{F}}} \boldsymbol{P}^\dag  \bw_h\lVert ^2_2 = (\bw_h)^\dag \boldsymbol{P}\hat{\boldsymbol{\mathcal{F}}}^\dag \boldsymbol{P}^\dag \boldsymbol{P} \hat{\boldsymbol{\mathcal{F}}} \boldsymbol{P}^\dag  \bw_h = (\bw_h)^\dag \boldsymbol{P}\hat{\boldsymbol{\mathcal{F}}}^\dag \hat{\boldsymbol{\mathcal{F}}} \boldsymbol{P}^\dag  \bw_h  = \sum_{i=1} |\hat{f}_i|^2 |\hat{w}_{h,i}|^2 \leq \sum_{i=1}  |\hat{w}_{h,i}|^2,
\end{split}
\end{equation}
where 
\begin{equation}
\boldsymbol{P}^\dag  \bw_h = \hat{\bw}_h \in \mathbb{C}^{N_u}
\end{equation}
and 
\begin{equation}
    |\hat{f}_i|^2 := \hat{f}_i^\dag \hat{f}_i.
\end{equation}
In order for this inequality to be satisfied we require $|\hat{f}_i|^2 \leq 1 \, \forall i=1,\dots, N_u$. Looking at \eqref{eq:definition_f}, this is satisfied as the eigenvalues for the diffusion operator are both real and negative \cite{sanderse2020non}. This condition is therefore sufficient to show diffusivity/stability of the EFR formulation, and independent of the value of $\delta$. \reviewerB{An alternative viewpoint is that the filter $\boldsymbol{\mathcal F} = (\boldsymbol{I} - 2\delta^2 \boldsymbol{D}_h)^{-1}$ 
is the discrete analogue of the resolvent of the Laplace operator
$J_\alpha = (I - \alpha \Delta)^{-1}$ 
\cite{Brezis2010}, which satisfies 
$\|J_\alpha\|_{L^2 \to L^2} \le 1$ for $\alpha \geq 0$.
Since $\boldsymbol{D}_h$ is a symmetric negative semidefinite discretization of $\Delta$, the same nonexpansive property holds in the discrete setting, 
i.e.\ $\|\boldsymbol{\mathcal F}\|_2 \le 1$.}

\section{Methods}
\label{sec:methods}
In this section we aim to increase the flexibility of the EFR framework with a data-driven approach and suggest a structure-preserving extension using the relaxation parameter $\chi$. A graphical version of the methodological pipeline can be found in Figure \ref{fig:flowchart}.

\subsection{Data-driven linear filter}\label{dd-filter}

In the previous Section we expressed the filter in the EFR framework as the matrix vector product $\boldsymbol{P}^\dag  \hat{\boldsymbol{\mathcal{F}}} \boldsymbol{P}$, where $\hat{\boldsymbol{\mathcal{F}}}$ is a diagonal matrix.
In this Section we present the following idea, which is one of the key novelties of our work: We can derive the diagonal elements $\hat{f}_i$ using a data-driven approach.
We will indicate with $\hat{\boldsymbol{\mathcal{F}}}^{\star}$ the \emph{optimal} filter, fitted to some reference data, and with $\hat{f}_i^{\star}, i=1, \dots, N_u$ its diagonal elements. It is important to emphasize that the optimal filter needs to be computed only \emph{once}, serving as the \emph{offline} stage in the typical surrogate modeling fashion, as depicted in the red dashed box of Figure \ref{fig:flowchart}.

For periodic boundary conditions and uniform grids the eigenvectors in $\boldsymbol{P}$ simply correspond to the Fourier modes. Using the notation $\hat{.}$ for the Fourier coefficients, we can use a \textbf{fast Fourier transform (FFT)} to efficiently obtain $\hat{\bw}_h$ \cite{DUHAMEL1990FFT}, i.e.
\begin{align}
    \hat{\bw}_h &= \boldsymbol{P}^\dag  \bw_h = \text{FFT}  (\bw_h), \\
    \bw_h &= \boldsymbol{P} \hat{\bw}_h = \text{FFT}^{-1}   (\hat{\bw}_h),
\end{align}
For non-periodic boundary conditions this is more complicated as we require the eigenvectors to satisfy the boundary conditions. This can possibly be done through the use of lifting functions like in \cite{rosenberger2023no}, but in this work we restrict ourselves for simplicity to the periodic case.

To derive the values for the complex filter coefficients $\hat{f}_i^{\star}$ we use reference data $\boldsymbol{U}^{\text{true}}$ and $\boldsymbol{W}^{\text{true}}$, collected within the time window $[0, \Ttrain]$.

In particular, we consider a parametrized flow case with parameter $\boldsymbol{\mu}$, given in our numerical tests by different initial conditions. For the derivation of the optimal filter, we consider:
\begin{itemize}
    \item $\boldsymbol{\mu}^{(j)}, j=1, \dots, \Itrain$ different initial conditions;
    \item $n=1, \dots, N_{\text{train}}$ time instances for each configuration, where $N_{\text{train}}$ corresponds to the training time $\Ttrain$. 
\end{itemize}

Matrix $\boldsymbol{U}^{\text{true}}$ contains in its columns the snapshots of the \emph{true} velocity field $\boldsymbol{u}^{n}_{\text{true}}(\boldsymbol{\mu}^{(j)})$, $n=2, \dots, N_{\text{train}}$; $j=1, \dots, \Itrain$. Then, $\boldsymbol{U}^{\text{true}} \in \mathbb{R}^{N_u \times M}$, with $M=N_{\text{train}}\,(\Itrain-1)$.
 The true velocity fields are obtained from fine-grid DNS data that is coarse-grained using a face-averaging filter, as described in \cite{Agdestein2025MLLES}.
Matrix $\boldsymbol{U}^{\text{true}}$ may be written as:

\begin{equation}
\boldsymbol{U}^{\text{true}} = \begin{bmatrix}
    |&  &|&&|&&|\\
    \boldsymbol{u}^2_{\text{true}}(\boldsymbol{\mu}^{(1)})& \dots &\boldsymbol{u}^{N_{\text{train}}}_{\text{true}}(\boldsymbol{\mu}^{(1)})&\dots&\boldsymbol{u}^{2}_{\text{true}}(\boldsymbol{\mu}^{(\Itrain)})&\dots&\boldsymbol{u}^{N_{\text{train}}}_{\text{true}}(\boldsymbol{\mu}^{(\Itrain)})\\
    |&  &|&&|&&|\\
\end{bmatrix}.
    \label{eq:ref-matrix}
\end{equation}

On the other hand, matrix $\boldsymbol{W}^{\text{true}}$ is obtained by evolving each true velocity field $\boldsymbol{u}^n_{\text{true}}$, $n=1, \dots, N_{\text{train}}-1$ for a single time step with the coarse grid solver, namely:

\begin{equation}
\boldsymbol{W}^{\text{true}} = \begin{bmatrix}
    |&  &|\\
    \boldsymbol{w}^2_{\text{true}}(\boldsymbol{\mu}^{(1)})& \dots &\boldsymbol{w}^{N_{\text{train}}}_{\text{true}}(\boldsymbol{\mu}^{(\Itrain)})\\
    |&  &|\\
\end{bmatrix}=\text{evolve} \left( \begin{bmatrix}
    |& &|\\
    \boldsymbol{u}^1_{\text{true}}(\boldsymbol{\mu}^{(1)}) & \dots & \boldsymbol{u}^{N_{\text{train}-1}}_{\text{true}}(\boldsymbol{\mu}^{(\Itrain)})\\
    |& &|\\
\end{bmatrix} \right).
    \label{eq:evolve-matrix}
\end{equation}
\reviewerB{Note that for the initial state $\boldsymbol{u}^1_{\text{true}} = \boldsymbol{w}^1_{\text{true}}$, which is why it is omitted from the snapshot matrices.}
In the numerical results, we will consider two different \emph{data regimes} to find the optimal filter: the \emph{scarce} data regime ($\Itrain=1$), and the \emph{full} data regime ($\Itrain=10$).

In this subsection, we focus on the filter step (equivalent to choosing $\chi=1$ in the relax step). The goal is to find filter coefficients such that the filtered coarse-grid solution is close to the reference solution $\boldsymbol{U}^\text{true}$:
\begin{equation}
    \boldsymbol{U}^\text{true} \approx \boldsymbol{\mathcal{F}}\boldsymbol{W}^{\text{true}}.
\end{equation}
A suitable error to satisfy this approximation is
\begin{equation}
    \mathcal{L}_{\hat{f}}(\boldsymbol{W}^\text{true},\boldsymbol{U}^\text{true}):=\lVert  \boldsymbol{\mathcal{F}} \boldsymbol{W}^{\text{true}} - \boldsymbol{U}^\text{true} \lVert ^2_F ,
\end{equation}
where the subscript $F$ indicates the Frobenius norm. As $\boldsymbol{P}^\dag $ is orthogonal we can place it inside the norm
\begin{equation}
    \mathcal{L}_{\hat{f}}(\boldsymbol{W}^\text{true},\boldsymbol{U}^\text{true}) = \lVert   \boldsymbol{P}^\dag  \boldsymbol{\mathcal{F}} \boldsymbol{W}^\text{true} - \boldsymbol{P}^\dag  \boldsymbol{U}^\text{true} \lVert _F^2 = \lVert   \hat{\boldsymbol{\mathcal{F}}} \hat{\boldsymbol{W}}^\text{true} - \hat{\boldsymbol{U}}^\text{true} \lVert _F^2,
    \label{eq:loss-simple}
\end{equation}
where $\hat{.}$ indicates the Fourier coefficients and we use the fact that $\boldsymbol{\mathcal{F}}=\boldsymbol{P}\hat{\boldsymbol{\mathcal{F}}} \boldsymbol{P}^\dag$. As $\hat{\boldsymbol{\mathcal{F}}}$ is diagonal Equation \eqref{eq:loss-simple} can be rewritten as
\begin{equation}
\begin{split}
    \mathcal{L}_{\hat{f}}(\boldsymbol{W}^\text{true},\boldsymbol{U}^\text{true}) &= \lVert   \hat{\boldsymbol{\mathcal{F}}} \hat{\boldsymbol{W}}^\text{true} - \hat{\boldsymbol{U}}^\text{true} \lVert _F^2 = \sum_i \lVert   \hat{f}_i \hat{\boldsymbol{W}}^\text{true}_{i} - \hat{\boldsymbol{U}}^\text{true}_{i} \lVert _2^2  \\ &= \sum_i ( |\hat{f}_i|^2 (\hat{\boldsymbol{W}}^\text{true}_{i})^\dag   \hat{\boldsymbol{W}}^\text{true}_{i} -  \hat{f}^\dag_i (\hat{\boldsymbol{W}}^\text{true}_{i})^\dag  \hat{\boldsymbol{U}}^\text{true}_i - (\hat{\boldsymbol{U}}^\text{true}_i)^\dag  \hat{f}_i \hat{\boldsymbol{W}}^\text{true}_{i}  + (\hat{\boldsymbol{U}}^\text{true}_i)^\dag \hat{\boldsymbol{U}}^\text{true}_i,
\end{split}
\end{equation}
where the subscript $i$ represents the $i$-th row of the matrices. The optimal values for $\hat{f}_i$ can be found by computing the gradient of $\mathcal{L}_{\hat{f}}$ with respect to $\hat{f}_i^\dag$ and setting it to zero:

\begin{equation}
\label{eq:loss-derivative}
    \frac{\partial \mathcal{L}_{\hat{f}}(\boldsymbol{W}^\text{true},\boldsymbol{U}^\text{true})}{\partial \hat{f}_i^\dag} = 2 \hat{f}_i (\hat{\boldsymbol{W}}^\text{true}_{i})^\dag \hat{\boldsymbol{W}}^\text{true}_{i}  - 2(\hat{\boldsymbol{W}}^\text{true}_{i})^\dag \hat{\boldsymbol{U}}^\text{true}_i = 0.
\end{equation}

Note we take the gradient with respect to $\hat{f}_i^\dag$ and not $\hat{f}_i$, to obtain an equation which solely depends on $\hat{f}_i$ and not $\hat{f}_i^\dag$. 
Solving \eqref{eq:loss-derivative} yields
\begin{equation}
\boxed{
    \hat{f}_i^* = \frac{(\hat{\boldsymbol{W}}^\text{true}_{i})^\dag  \hat{\boldsymbol{U}}^\text{true}_i}{ (\hat{\boldsymbol{W}}^\text{true}_{i})^\dag  \hat{\boldsymbol{W}}^\text{true}_{i}  }}
    \label{eq:f_star}
\end{equation}
for the optimal values $\hat{f}_i^*$, $i=1, \dots, N_u$. As $\boldsymbol{U}^\text{true}$ is real-valued the resulting $\boldsymbol{\mathcal{F}}^{\star}$ is also real. 
Simulations are carried out using this filter and setting $\chi = 1$, such that we effectively skip the relax step. This will be referred to as \ddEF{}. The filtered velocity field will be denoted as
\begin{equation}
    \overline{\boldsymbol{w}}_h := \boldsymbol{\mathcal{F}}^{\star} \boldsymbol{w}_h = \boldsymbol{P}^\dag  \hat{\boldsymbol{\mathcal{F}}}^{\star} \boldsymbol{P} \boldsymbol{w}_h ,
\end{equation}
which is computed efficiently in the frequency domain using a FFT, as $\hat{\boldsymbol{\mathcal{F}}}^{\star}$ is diagonal.

We finally remark that this data-driven filter can be built offline from a limited amount of data, and then employed in online simulations in extrapolation regimes, both in time and with respect to varying parameters such as initial conditions or forcing terms. This separation between offline training and online deployment enables fast and flexible evaluations across multiple scenarios, without the need for recomputing the filter.

\subsection{Conservation of energy using $\chi$}
\label{sec:tuning-chi-energy}

As discussed in \ref{subsec:stability_diff_filter} one the of the desirable properties of the existing differential filter is that it guarantees stability of the system by dissipating kinetic energy. This is guaranteed by the fact that \( |\hat{f}_i| \leq 1,\ \forall i \) for this filter.

However, this condition is in general not met by the filter in \eqref{eq:f_star}. Requiring the same bound on the \ddEF{} approach seems natural, but can lead to overdiffusive behavior, as it enforces energy dissipation for each frequency.
In fact, $|\hat{f}_i^{\star}|^2 \leq 1$ is a \emph{sufficient} but not a \emph{necessary} condition to achieve stability. Instead, we therefore impose a \emph{global} constraint on the energy, which ensures stability while still allowing for some frequencies to have an increasing energy.

This constraint can be enforced by a smart choice of the relaxation parameter $\chi$ at each time step \( t^n \), as follows:
\begin{equation}
\mathcal{E}^n_h \propto \lVert \bu_h^n \rVert_2^2 = \lVert (1 - \chi) \bw_h^n + \chi \overline{\bw}_h^n \rVert_2^2 \leq \lVert \bw_h^n \rVert_2^2, 
\label{eq:global-constraint-space}
\end{equation}
where $\mathcal{E}^n_h$ is the discrete global energy at $t_n$. In these expressions, \( \bw_h^n \), \(  \overline{\bw}_h^n\), and \( \bu_h^n \) correspond to the \emph{online} solution approximations at time step \( t^n \), and not to the snapshot matrices \( \boldsymbol{W}_{\text{true}} \) and \( \boldsymbol{U}_{\text{true}} \) used in the offline computation of \( \hat{\boldsymbol{\mathcal{F}}}^{\star} \) in Section~\ref{dd-filter}. 

One effective approach to enforce the global dissipativity constraint is to compute the filter \emph{offline}, as in equation~\eqref{eq:f_star}, and adaptively tune \( \chi^n = \chi(t^n) \) at each \emph{online} time step, as represented in the blue dashed box in Figure \ref{fig:flowchart}.
From the constraint~\eqref{eq:global-constraint-space}, we derive the following inequality:
\begin{equation}
    (\chi^n)^2 \tilde{a}^n  + 2 \chi^n \tilde{b}^n \leq 0,
\label{eq:chi-online-1}
\end{equation}
where 
\[
\tilde{a}^n = \lVert  \bw_h^n - \overline{\bw}_h^n \rVert_2^2 \quad \text{and} \quad \tilde{b}^n =   (\bw_h^n)^T  \overline{\bw}_h^n - \lVert \bw_h^n \lVert_2^2.
\]
Given that \( 0 \leq \chi^n \leq 1 \), to remain consistent with the relax step in the EFR approach (see Section \ref{subsec:efr}), inequality~\eqref{eq:chi-online-1} yields the bound:

\begin{equation}
\boxed{
\begin{cases}
   & \chi^n \leq -\frac{2\tilde{b}^n}{\tilde{a}^n} \quad \text{if } \tilde{b}^n \leq 0, \\
   & \chi^n = 0 \quad \text{otherwise}.
\end{cases}
}
\label{eq:condition-energy}
\end{equation}

Accordingly, during the online phase, we select the optimal relaxation parameter (we refer to it as \chienergy{}) at each time step \( t^n \) as follows:
\begin{itemize}
    \item If constraint~\eqref{eq:global-constraint-space} is satisfied for $\chi = 1$:
    \[
    \chienergy = 1,
    \]
    i.e., we only perform a filtering step. Relaxation is not needed since the filtered solution satisfies the global energy inequality. 
    \item If constraint~\eqref{eq:global-constraint-space} is \emph{not} satisfied:
    \begin{itemize}
        \item \( \chienergy = -\frac{2\tilde{b}^n}{\tilde{a}^n} \), if \( \tilde{b}^n \leq 0 \) (the least dissipative value that ensures the constraint);
        \item $ \chienergy = 0 $, if $\tilde{b}^n > 0$. This means that we simply use the coarse grid solution without applying any filtering or relaxation. While this yields an energy-stable result, it may exhibit spurious oscillations due to the coarse resolution, making it a less desirable option in practice.
    \end{itemize}
\end{itemize}
This approach will be referred to as \energyEFR. Note that these computations are all carried out in physical space, as the FFT is only required to efficiently apply the data-driven filter. Moreover, computing the inner products in this procedure in the frequency domain does not offer any computational speedup.

\subsection{Conservation of enstrophy using $\chi$}
\label{sec:tuning-chi-enstrophy}

The derivation in Section \ref{sec:tuning-chi-energy} does not inherently ensure a monotonic decrease in enstrophy, a property that should hold for the viscous 2D NSE \cite{Palha2017enstrophy}. Preserving this property is important, as an increase in enstrophy often indicates the emergence of unphysical oscillations the numerical solution. This occurs because enstrophy is sensitive to high-frequency modes, as it involves the squared vorticity norm:
\begin{equation}
\mathcal{Z} := \frac{1}{2|\Omega|} \int_\Omega \lVert \nabla \times \mathbf{u} \rVert_2^2 \, \mathrm{d}\Omega.
\end{equation}
 In Section~\ref{sec:results}, and specifically in Figure~\ref{fig:scarce-test-general}, we demonstrate that enforcing this physical constraint improves simulation accuracy, particularly in data-scarce training regimes.

To satisfy this physical constraint, we impose the following inequality at each online time step $t_n$:
\begin{equation}
    \mathcal{Z}^n_h := \frac{h^2}{2|\Omega|} \| \boldsymbol{B}_h \bu_h^n \|^2_2 \leq \frac{h^2}{2|\Omega|}  \| \boldsymbol{B}_h \bw_h^n \|^2_2,
    \label{eq:enstrophy-inequality}
\end{equation}
where $\boldsymbol{B}_h \in \mathbb{R}^{N_p \times N_u}$ is a discrete curl operator, such that $\mathcal{Z}^n_h$ approximates the global enstrophy.
The above expression may be rewritten in a form similar to the one found for the energy-preserving inequality, namely:
\begin{equation}
    (\chi^n)^2 \, \tilde{a}^n_{\mathcal{Z}}  + 2\chi^n \tilde{b}^n_{\mathcal{Z}} \leq 0,
    \label{eq:enstrophy-ineq-compact}
\end{equation}
where 
\[
\tilde{a}^n_{\mathcal{Z}} = \lVert \boldsymbol{B}_h\bw_h^n - \boldsymbol{B}_h\overline{\bw}_h^n \rVert_2^2 \quad \text{and} \quad \tilde{b}^n_{\mathcal{Z}} = (\boldsymbol{B}_h \bw_h^n)^T \boldsymbol{B}_h \overline{\bw}_h^n - \lVert \boldsymbol{B}_h \bw_h^n \lVert_2^2.
\]
Exactly as for the energy inequality, the optimal $\chi$ satisfying inequality \eqref{eq:enstrophy-ineq-compact} is $\chienstrophyonly = -2\frac{\tilde{b}^n_{\mathcal{Z}}}{\tilde{a}^n_{\mathcal{Z}}}$.
Also in this case, we have:
\begin{equation}
\boxed{
\begin{cases}
   & \chi^n \leq -\frac{2\tilde{b}_{\mathcal{Z}}^n}{\tilde{a}_{\mathcal{Z}}^n} \quad \text{if } \tilde{b}_{\mathcal{Z}}^n \leq 0, \\
   & \chi^n = 0 \quad \text{otherwise}.
\end{cases}
}
\label{eq:condition-enstrophy}
\end{equation}

As for $\tilde{b}^n$, $\tilde{b}^n_{\mathcal{Z}}$ may also be negative. In that case, one would simply retain $\chi_{\mathcal{Z}-\text{opt}}=0$. This approach can easily be combined with the energy-conserving approach by simply setting $\chi$ to $\chienstrophy = \min \{\chienergy,\chienstrophyonly\}$ to minimize dissipation.
While conceptually straightforward, its implementation is slightly more involved due to the need to compute discrete gradients within the staggered grid framework.
We refer to this approach as \enstrophyEFR{}.


\section{Numerical results}
\label{sec:results}

The methodology presented in Section \ref{sec:methods} is tested here on two test cases: 2D decaying homogeneous turbulence (\ref{test-case-1}), and the 2D Kolmogorov flow (\ref{test-case-2}). Both cases are performed in a square domain $\Omega=[0, 1]^2$ at Reynolds number $Re=\num{4e4}$, with periodic boundary conditions. Time integration is performed using a fourth-order explicit Runge–Kutta method~\cite{sanderse2012accuracy}.

The DNS simulation is carried out on a refined grid, whereas both the proposed and state-of-the-art methods are implemented and compared on a coarse grid. The goal is to find a filtering strategy that allows to obtain accurate results despite the reduced resolution.
More details \reviewerA{on the different methods compared in the article and} on the \reviewerA{corresponding} mesh sizes are available in Table \ref{tab:grids}.

We will refer to the DNS simulation on the coarse grid as \emph{noEFR}.

\begin{table}[htpb!]
    \centering
    \setlength{\tabcolsep}{8pt} 
    \renewcommand{\arraystretch}{1.2} 
    \begin{tabular}{>{\centering\arraybackslash}m{1.7cm}  
                    >{\centering\small\arraybackslash}m{8.5cm}    
                    >{\centering\arraybackslash}m{1.5cm}  
                    >{\centering\arraybackslash}m{1.8cm}} 

    \toprule
         \textbf{Simulation type} & \reviewerA{\textbf{Description}} & \textbf{Resolution} & \textbf{Mesh size}  \\
         \midrule
         DNS &\reviewerA{Fully-resolved simulation of all turbulence scales}& $512^2$ & $\num{1.95e-3}$ \\
         \midrule
         \reviewerA{Filtered DNS} &\reviewerA{DNS with small scales removed via explicit filter}& \multirow{14}{*}{$128^2$} & \multirow{14}{*}{$\num{7.81e-3}$} \\
         \cmidrule{1-2}
          \reviewerA{noEFR} & \reviewerA{DNS on coarse grid (i.e., no filter, no regularization)} && \\
          \cmidrule{1-2}
          \reviewerA{Standard EF} & \reviewerA{Standard approach for \emph{Evolve--Filter} using the differential filter in Equation \eqref{eq:diff-filter-ref} (with optimized $\delta$)} && \\
          \cmidrule{1-2}
          \reviewerA{Standard EFR} & \reviewerA{Standard approach for \emph{Evolve--Filter--Relax} using the differential filter in Equation \eqref{eq:diff-filter-ref} (with prescribed $\delta=h$ and optimized $\chi$)} && \\
          \cmidrule{1-2}
          \reviewerA{Smagorinsky} & \reviewerA{Smagorinsky model (with optimized parameter $\theta$)} && \\
          \cmidrule{1-2}
        \reviewerA{DD-EF} & \reviewerA{\emph{Evolve--Filter} using
        the data-driven filter introduced in Section \ref{dd-filter}} && \\
          \cmidrule{1-2}
          \reviewerA{\energyEFR{}} & \reviewerA{Energy-conserving \emph{Evolve--Filter--Relax} with:
		  \begin{itemize}[topsep=0pt, partopsep=0pt, itemsep=0pt, parsep=0pt]
              \item the data-driven filter introduced in Section \ref{dd-filter};
              \item the relax parameter adjusted to ensure energy dissipation, as in Section \ref{sec:tuning-chi-energy}.
          \end{itemize}}&& \\
          \cmidrule{1-2}
          \reviewerA{\enstrophyEFR{}} & \reviewerA{Energy and enstrophy-conserving \emph{Evolve--Filter--Relax} with:
		  \begin{itemize}[topsep=0pt, partopsep=0pt, itemsep=0pt, parsep=0pt]
              \item the data-driven filter introduced in Section \ref{dd-filter};
              \item the relax parameter adjusted to ensure energy dissipation as in Section \ref{sec:tuning-chi-energy} and enstrophy dissipation as in Section \ref{sec:tuning-chi-enstrophy}.
          \end{itemize}} && \\
         \bottomrule
    \end{tabular}
    \caption{\reviewerA{Brief description of the different types of simulation, including the }details on the grids employed.}
    \label{tab:grids}
\end{table}

The performance of the proposed methodology is compared with state-of-the-art approaches with respect to the filtered DNS on the coarse grid, where the filter is a simple face-average~\cite{Agdestein2025MLLES}.

The accuracy is measured as the discrepancy with respect to the filtered DNS in the global kinetic energy $\mathcal{E}^n_h$, in the global enstrophy $\mathcal{Z}^n_h$, for $n=1, \dots, N$, and in the energy spectrum.

More in detail, we employ a range of initial conditions to evaluate the robustness of the methods, and we analyze the results in terms of averages and 95 $\%$ confidence intervals across multiple initializations.

Table \ref{tab:train-test} shows train and test settings adopted for the different methodologies, specifying the training time window \Ttrain{}, the simulation time $T$, the number of training and test initializations (\Itrain{} and \Itest{}, respectively), and the corresponding results' Sections.
We also specify that \Ttrain{} has different meanings for DD-EF(R) methods and for state-of-the-art approaches \reviewerA{(Standard EFR, Standard EF, and Smagorinsky)}. On the one hand, in the DD-EF(R) methods, $[0, \Ttrain]$ is the simulation time window used to collect the snapshots' matrices $\boldsymbol{U}^{\text{true}}$ and $\boldsymbol{W}^{\text{true}}$ employed in the optimization problem (Section \ref{dd-filter}). On the other hand, for state-of-the-art approaches it is the simulation time window considered to perform hyperparameter tuning. For more details about the state-of-the-art approaches, we refer the reader to \ref{literature-tuning}.

\begin{table}[htpb!]
    \centering
    \begin{tabular}{cccccc}
    \toprule
       \multirow{2}{*}{\textbf{Method}}& \multicolumn{2}{c}{\textbf{Train} }&\multicolumn{2}{c}{\textbf{Test} }&\multirow{2}{*}{\textbf{Section(s)}} \\
        \cmidrule{2-5}
         &$\Ttrain$ & $\Itrain$ &$T$ & $\Itest$& \\
         \midrule
         State-of-the-art methods&$1$ s& 10& \multirow{4}{*}{$10$ s} & \multirow{4}{*}{5}&\ref{literature-tuning}\\
         \cmidrule{1-3}  \cmidrule{6-6}
         \ddEF{}, \energyEFR{}  & \multirow{2}{*}{$1$ or $3$ s} & 10 (full data)&&&\ref{results-dd-efr-case1}, \ref{results-dd-efr-case2}\\
         \cmidrule{3-3}
       \cmidrule{1-1} \cmidrule{6-6}
         \ddEF{}, \energyEFR{}, \enstrophyEFR{} && 1 (scarce data)&&&\ref{results-scarce-dd-efr-case1}, \ref{results-scarce-dd-efr-case2}\\
         \bottomrule
    \end{tabular}
    \caption{Train and test setup for DD-EF(R) and for state-of-the-art approaches, and the corresponding results' Sections where the methods are applied.}
    \label{tab:train-test}
\end{table}

\subsection{Test case 1: two-dimensional decaying homogeneous turbulence}
\label{test-case-1}

The first test case is decaying homogeneous isotropic turbulence, \reviewerB{already treated }
in \cite{Agdestein2025MLLES, kochkov2021machine, kurz2023deep}.
The initial DNS velocity is sampled from a random velocity field with prescribed energy spectrum, \reviewerB{as in } \cite{orlandi2000fluid, san2012high}.
More in detail, as in \cite{Agdestein2025MLLES}, the velocity field is initialized with the following procedure:

\begin{enumerate}
    \item The velocity is sampled in spectral space to match a prescribed energy $\hat{\mathcal{E}}_{\kappa}$ at each wave number $\kappa$. Each component of the spectral velocity is indeed initialized such that $\|\hat{\boldsymbol{u}}_{\kappa} \|= |a_{\kappa}|$, where $a_{\kappa}=\sqrt{2\hat{\mathcal{E}}_{\kappa}} \text{e}^{2\pi i \tau_{\kappa}}$ and $\tau_{\kappa}$ is a random phase shift.
    \item The velocity is then projected to make it divergence-free: $\hat{\boldsymbol{u}}_{\kappa}=\frac{a_{\kappa} \hat{\boldsymbol{P}}_{\kappa} \boldsymbol{e}_{\kappa}}{\hat{\boldsymbol{P}}_{\kappa} \boldsymbol{e}_{\kappa}}$, where $\hat{\boldsymbol{P}}_{\kappa}$ is a projector: $\hat{\boldsymbol{P}}_{\kappa}=I-\frac{\kappa \kappa^T}{\kappa^T \kappa}$ and $\boldsymbol{e}_{\kappa}$ is a random unit vector. In 2D $\boldsymbol{e}_{\kappa}=(\cos{(\theta_{\kappa})}, \sin{(\theta_{\kappa})})$ with $\theta_{\kappa} \sim \mathcal{U}[0, 2\pi]$.
    \item The velocity field is finally obtained by a transformation into physical space ($\boldsymbol{u}_h^1 = \text{FFT}^{-1}(\hat{\boldsymbol{u}})$), and a further projection step such that $\textbf{div}_h \boldsymbol{u}_h^1=\boldsymbol{0}$.
\end{enumerate}

\begin{figure}[htpb!]
    \centering
    \subfloat[$\boldsymbol{B}_h \bu_h^{1}(\boldsymbol{\mu}^{(1)})$]{\includegraphics[width=0.2\linewidth]{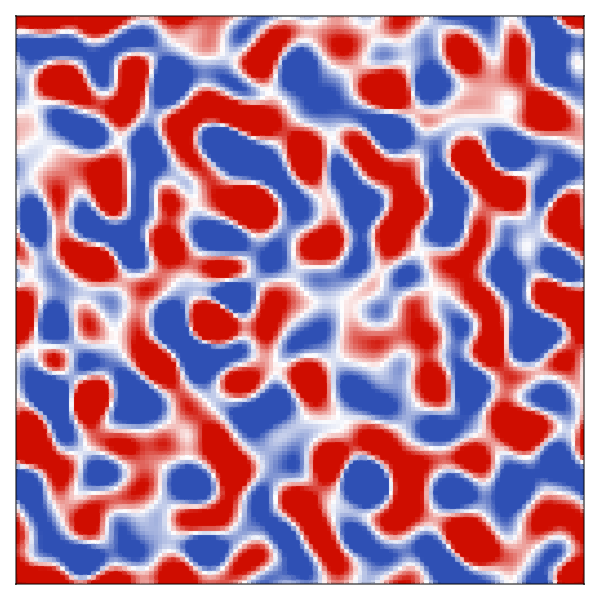}
    }
    \hspace{0.5cm}
    \subfloat[$\boldsymbol{B}_h \bu_h^{1}(\boldsymbol{\mu}^{(5)})$]{\includegraphics[width=0.2\linewidth]{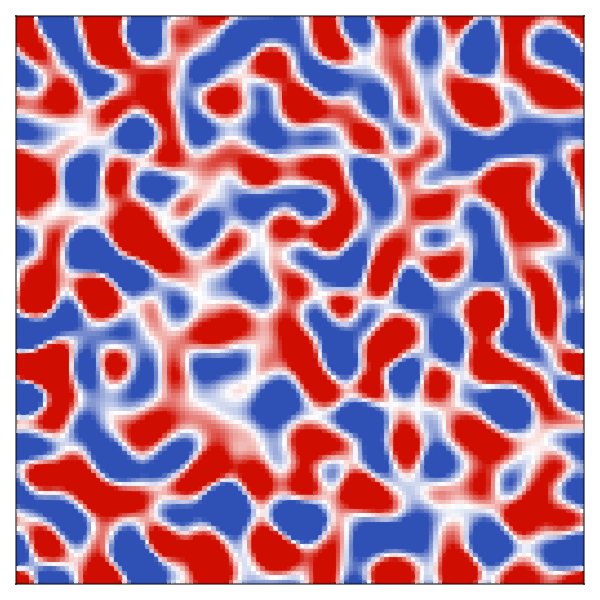}
    }
    \hspace{0.5cm}
    \subfloat[$\boldsymbol{B}_h \bu_h^{1}(\boldsymbol{\mu}^{(10)})$]{\includegraphics[width=0.2\linewidth]{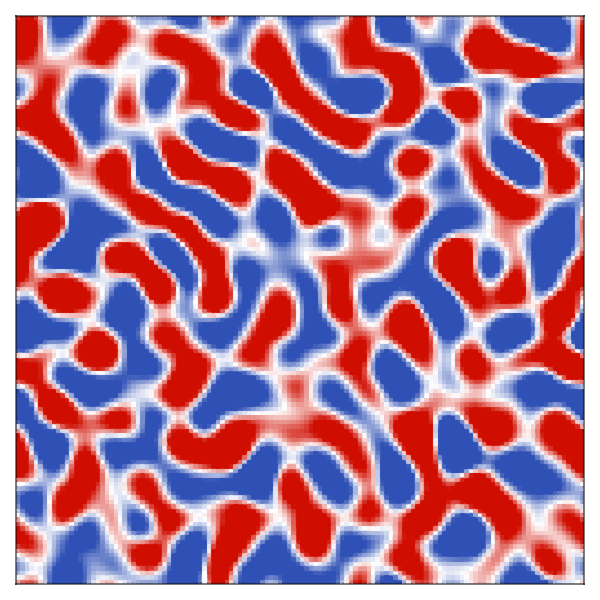}}
     \hspace{0.1cm}\subfloat{\includegraphics[width=0.05\linewidth]{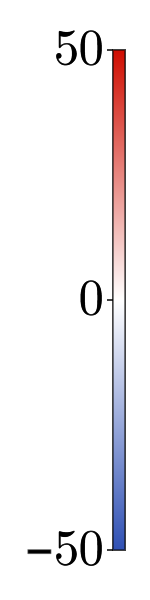}}
    \caption{Vorticity fields for initial time step, at different random initializations.}
    \label{fig:random-init}
\end{figure}

Following the notation in subsection \ref{dd-filter}, the initial velocity field can be written as $\bu_h^{1}(\boldsymbol{\mu}^{(j)})$, where the superscript $j$ represents the random seed used to generate $\theta_{\kappa}$ and $\tau_{\kappa}$ for each wavenumber $\kappa$.
Three examples of random initializations are represented in Figure \ref{fig:random-init}.

\subsubsection{Data-driven EFR in the case of full data}
\label{results-dd-efr-case1}

This part of the manuscript is dedicated to the results of \ddEF{} and \energyEFR{} approaches (presented in Sections \ref{dd-filter} and \ref{sec:tuning-chi-energy}, respectively).

The first \emph{offline} step involves constructing the filter matrix in the frequency domain (Equation \eqref{eq:f_star}) considering a full data regime ($\Itrain=10$ random initializations).
Additionally, we remark that the snapshot data is subsampled by retaining one snapshot every 10 time steps to reduce the computational burden.

\begin{figure}[htpb!]
    \centering
    \includegraphics[width=\linewidth]{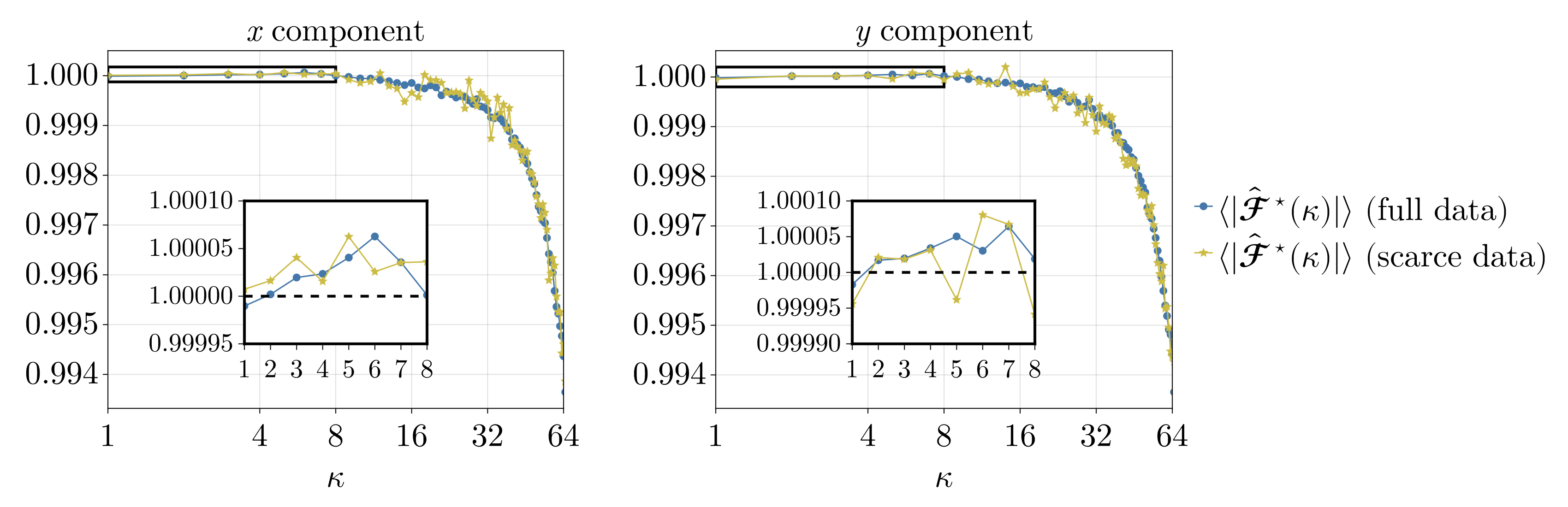}
    \caption{2D average of components $x$ and $y$ of matrix $\hat{\boldsymbol{\mathcal{F}}}^{\star}$, for the cases $\Ttrain=3$ s and $\Itrain=10$ (full data), and $\Ttrain=1$ s and $\Itrain=1$ (scarce data).}
    \label{fig:fstars}
\end{figure}

We analyze the effect of the filter on the velocity field in Figure \ref{fig:fstars} by considering the 2D average of the filter elements $\hat{\boldsymbol{\mathcal{F}}}^{\star}$ for both $x$ and $y$ components. We first associate an average filter value to all modes $\kappa$ with wavenumber magnitude $\kappa=\| \boldsymbol{k}\|= \sqrt{k_x^2 + k_y^2}$. Such average value is computed as the mean magnitude of the filter coefficients $\hat{f}^{\star}_i$ across all modes:
\begin{equation}
    \langle |\boldsymbol{\hat{{\mathcal{F}}}^{\star}}(\kappa)|\rangle=\frac{1}{|\mathcal{B}_{\kappa}|}\sum_{i=1}^{|\mathcal{B}_{\kappa}|}\hat{f}^{\star}_i \mathcal{M}_i(\kappa),
\end{equation}
where $\mathcal{B}_k$ collects the modes satisfying $\kappa-0.5<\|\boldsymbol{k}\|<\kappa+0.5$, while $\mathcal{M}_i$ is a mask defined as:
\[\mathcal{M}_i(\kappa)=
\begin{cases}
    1 \text{ if }\kappa-0.5<\|\boldsymbol{k}_i\|<\kappa+0.5,\\
    0 \text{ otherwise,}
\end{cases}
\]
and $|\mathcal{B}_{\kappa}|=\sum_{i}\mathcal{M}_i(\kappa)$ is the number of Fourier modes in $\mathcal{B}_{\kappa}$.

Values greater than one indicate that the filter amplifies the corresponding frequency components, while values smaller than one correspond to attenuation. A magnitude close to one implies that the filter leaves those frequencies mostly unchanged.
The plot shows that the filter slightly amplifies the low frequencies, which are the highest-energy ones, while it attenuates the medium-to-high frequencies.
The data-driven filter is highly influenced by the amount of data considered in the least squares problem. Indeed, the plot highlights that the filter built in the \emph{data scarcity} regime presents more oscillations and the optimization may converge to values bigger than one also at medium frequencies ($\kappa=9, 10, 14$ for the $y$ component). The energy injection at medium scales may increase significantly the velocity magnitude, and the simulation might blow up. This instability may be overcome by introducing the energy and/or enstrophy dissipation constraints in the relax step, as presented in Section \ref{sec:tuning-chi-energy} and \ref{sec:tuning-chi-enstrophy}. We will analyse the combined effect of such constraints in Section \ref{results-scarce-dd-efr-case1}.

\begin{figure}[htpb!]
    \centering
    \includegraphics[width=\linewidth]{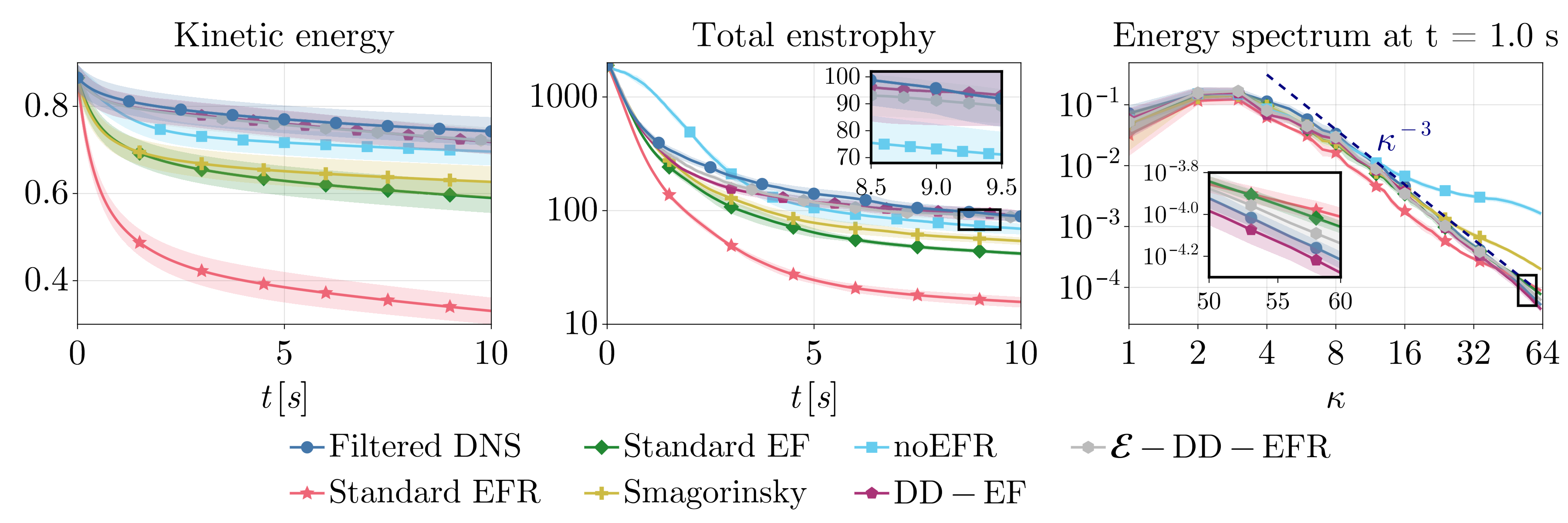}
    \caption{Time evolution of total kinetic energy, total enstrophy, and spectrum of the kinetic energy at time $t=1$ s, for the proposed methodologies (\ddEF{} and \energyEFR{}) and for state-of-the-art approaches.
    In particular, we show the case $\Ttrain=3$ s.
    For each method, the Figure shows the average values (solid lines) and the 95$\%$ confidence interval among 5 test configurations, in decaying turbulence test case.}
    \label{fig:test-general}
\end{figure}

\begin{figure}[htpb!]
    \centering
    \subfloat[Filtered DNS]{\includegraphics[width=0.18\linewidth]{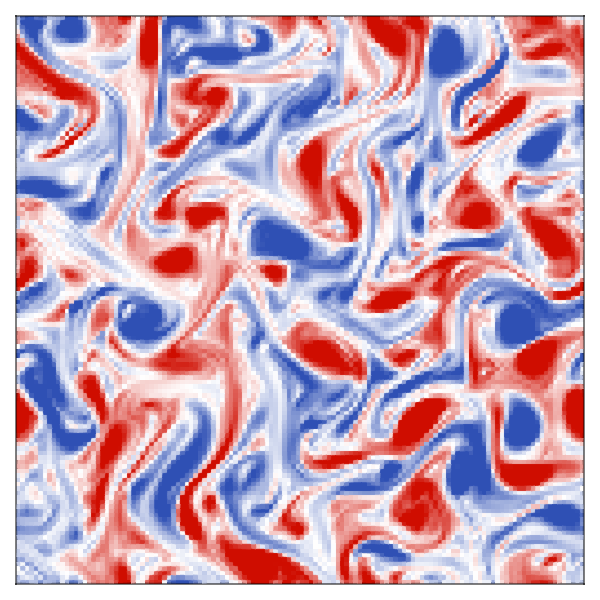}\label{filtered-dns-t03}}
    \hspace{0.5cm}
    \subfloat[\centering \ddEF \\ $\Ttrain=3$]{\includegraphics[width=0.18\linewidth]{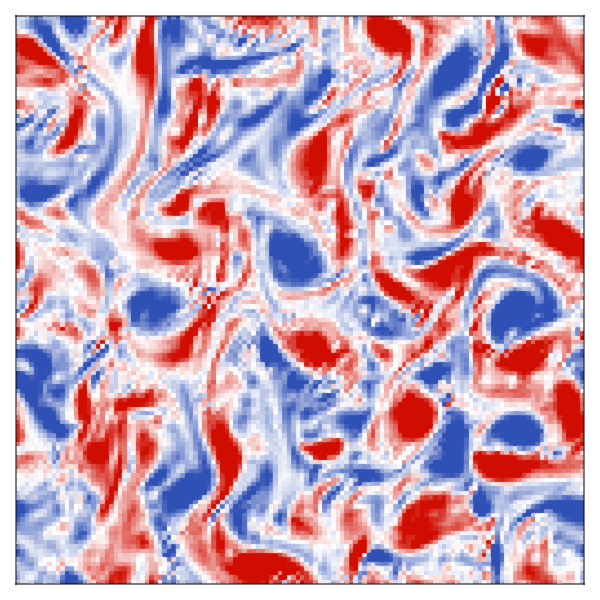} \label{dd-efr-t03}}
     \hspace{0.5cm}
     \subfloat[\centering \energyEFR \\ $\Ttrain=3$]{\includegraphics[width=0.18\linewidth]{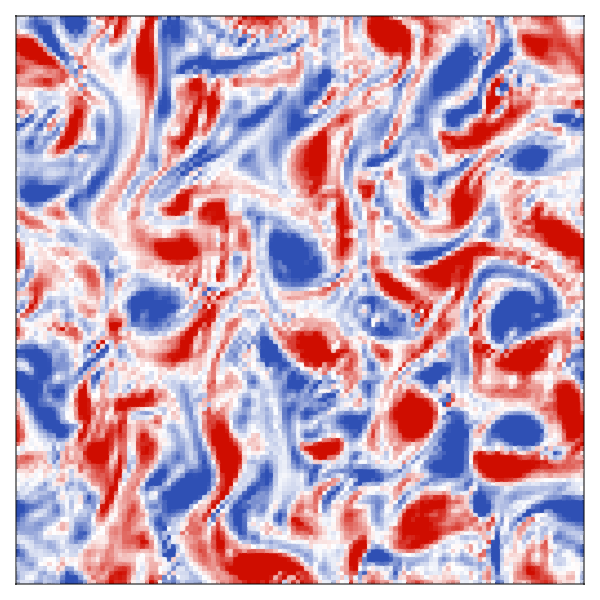}\label{energy-dd-efr-t03}}
     \hspace{0.5cm}\subfloat{\includegraphics[width=0.04\linewidth]{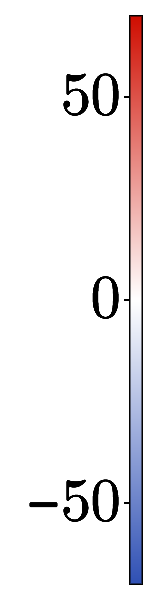}}\\
     \subfloat[noEFR]{\includegraphics[width=0.18\linewidth]{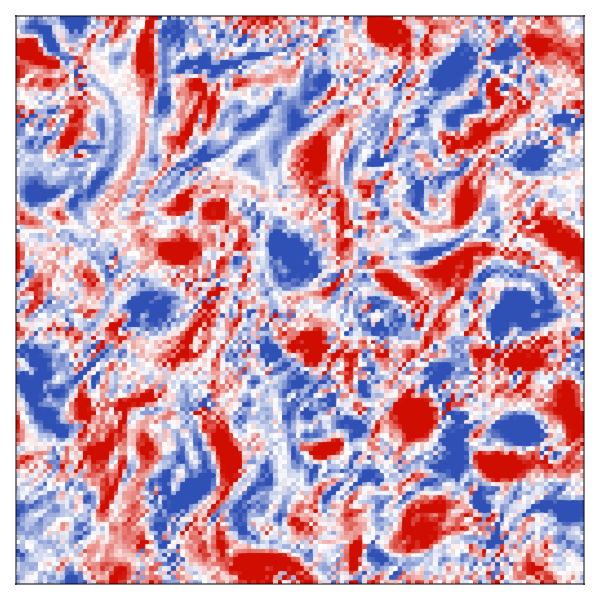}
    \label{noefr-t03}}
    \hspace{0.5cm}
    \subfloat[Smagorinsky]{\includegraphics[width=0.18\linewidth]{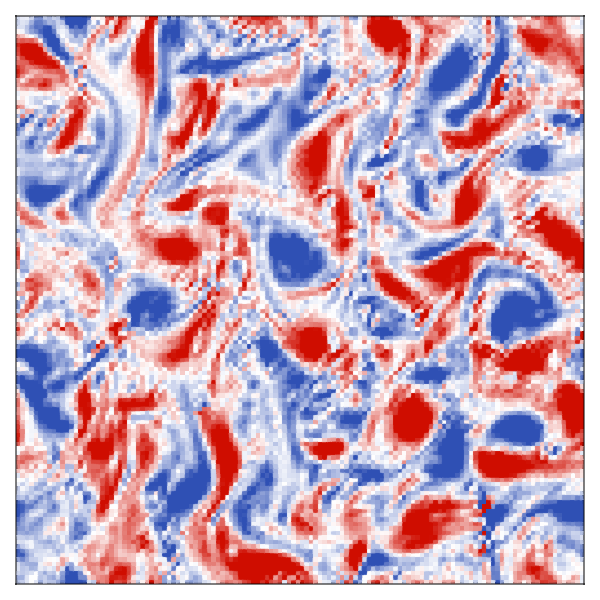}\label{smag-t03}}
    \hspace{0.5cm}
     \subfloat[Standard EFR]{\includegraphics[width=0.18\linewidth]{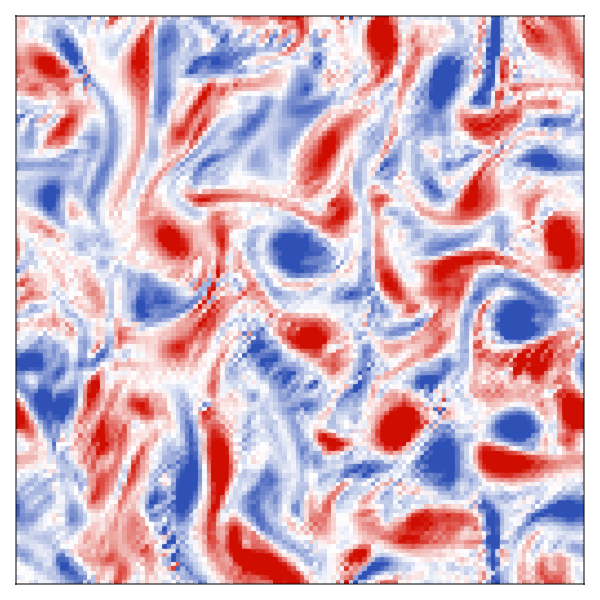}\label{std-efr-t03}}
     \hspace{0.5cm}
     \subfloat[Standard EF]{\includegraphics[width=0.18\linewidth]{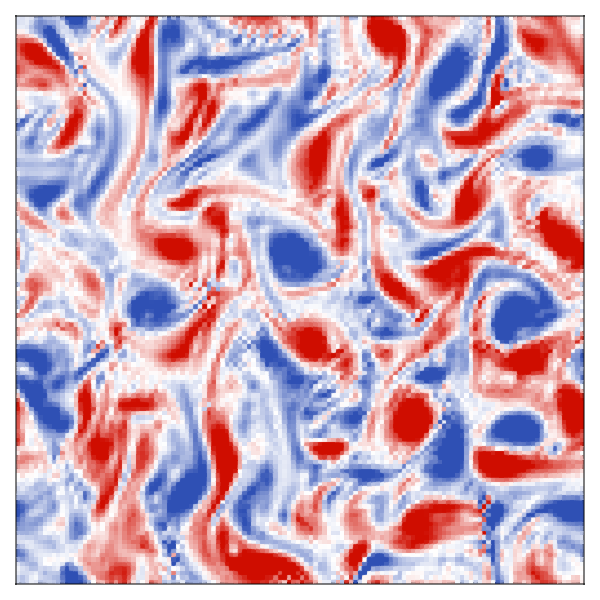}\label{std-ef-t03}}
     \hspace{0.5cm}
\subfloat{\includegraphics[width=0.04\linewidth]{images/decayingturbulence/colorbar_70.png}}
    \caption{Vorticity fields at $t=0.3$ for different methodologies in a test configuration, in decaying turbulence test case.}
    \label{fig:fields-t03}
\end{figure}

\begin{figure}[htpb!]
    \centering
    \subfloat[Filtered DNS]{\includegraphics[width=0.18\linewidth]{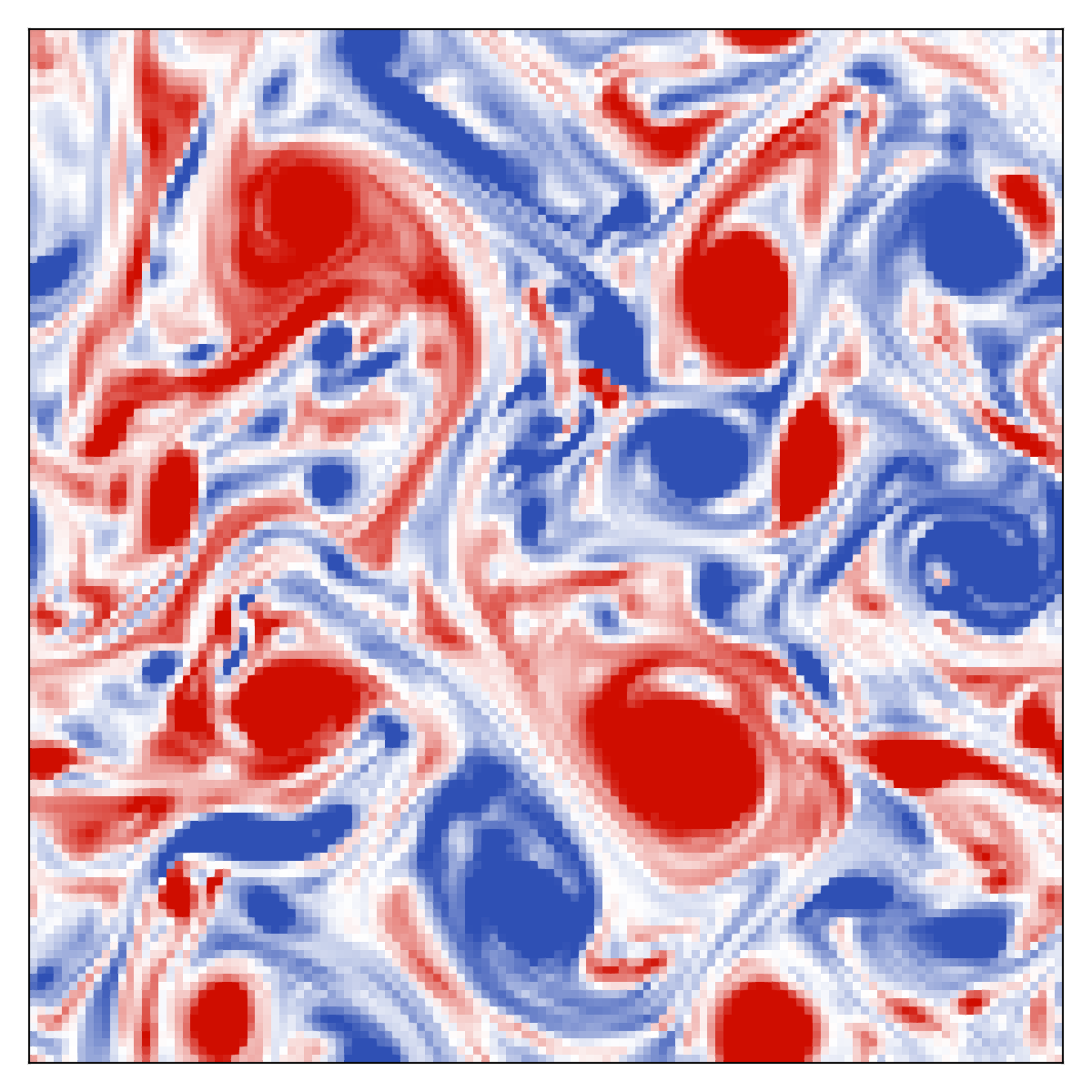}\label{filtered-dns-t1}}
    \hspace{0.5cm}
    \subfloat[\centering \ddEF \\ $\Ttrain=3$]{\includegraphics[width=0.18\linewidth]{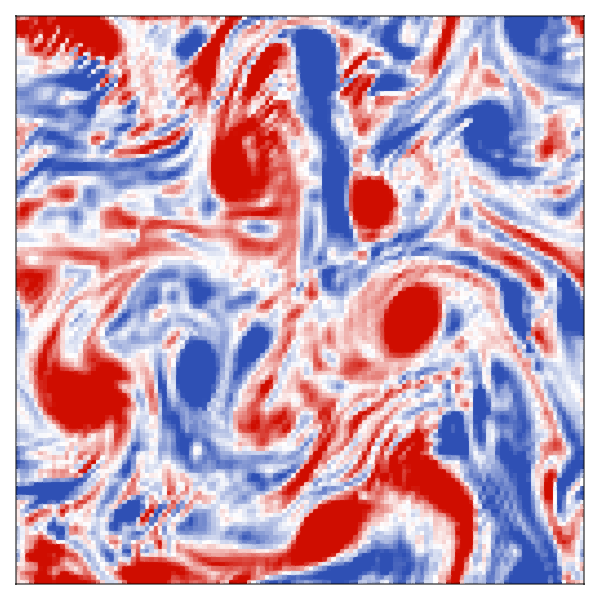}\label{dd-efr-t1}}
     \hspace{0.5cm}
     \subfloat[\centering \energyEFR \\ $\Ttrain=3$]{\includegraphics[width=0.18\linewidth]{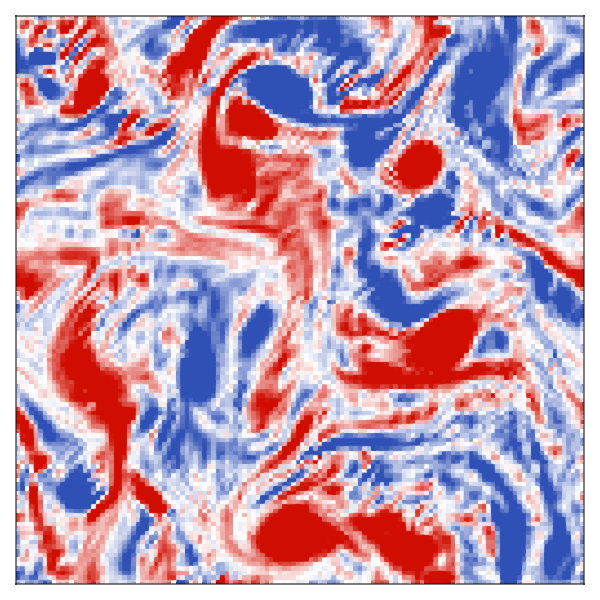}\label{energy-dd-efr-t1}}
     \hspace{0.5cm}\subfloat{\includegraphics[width=0.04\linewidth]{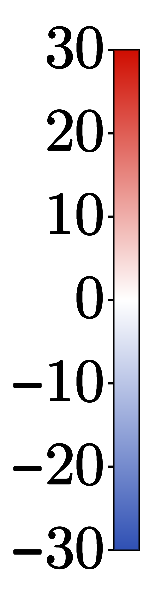}}\\
     \subfloat[noEFR]{\includegraphics[width=0.18\linewidth]{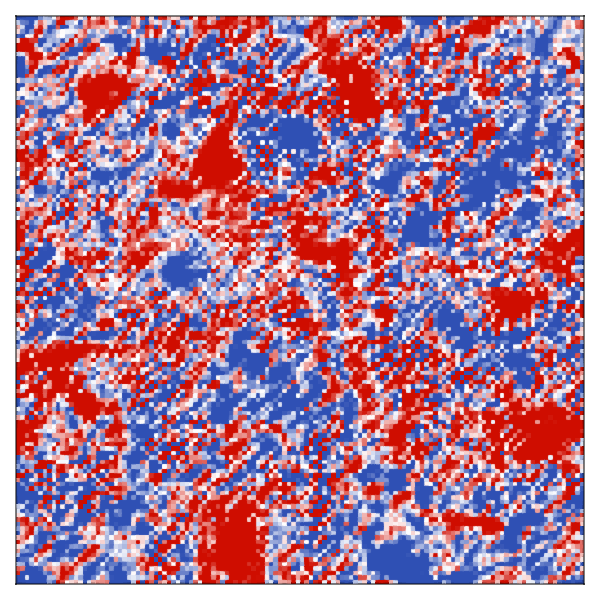}
     \label{noefr-t1}
    }
    \hspace{0.5cm}
    \subfloat[Smagorinsky]{\includegraphics[width=0.18\linewidth]{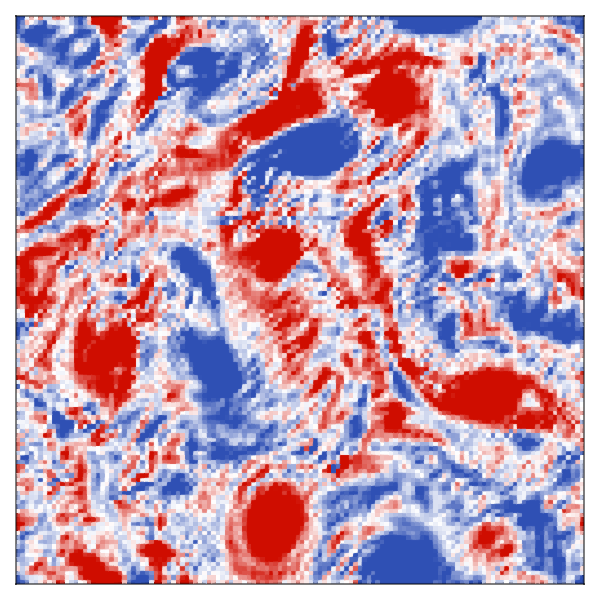}
    \label{smag-t1}}
    \hspace{0.5cm}
     \subfloat[Standard EFR]{\includegraphics[width=0.18\linewidth]{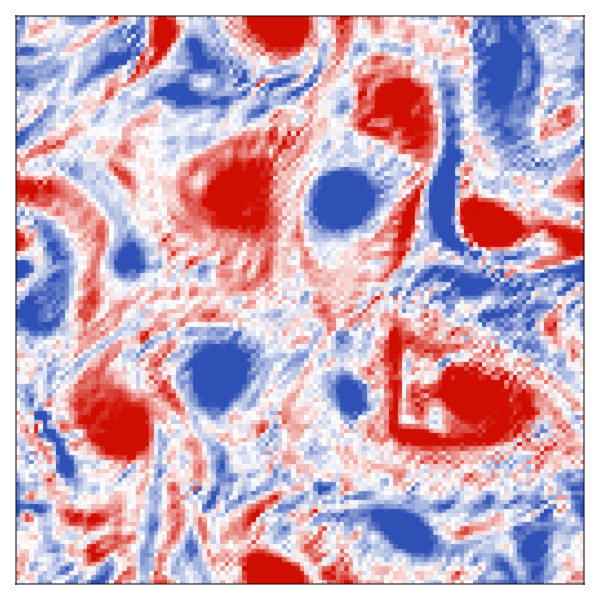}
     \label{std-efr-t1}}
     \hspace{0.5cm}
     \subfloat[Standard EF]{\includegraphics[width=0.18\linewidth]{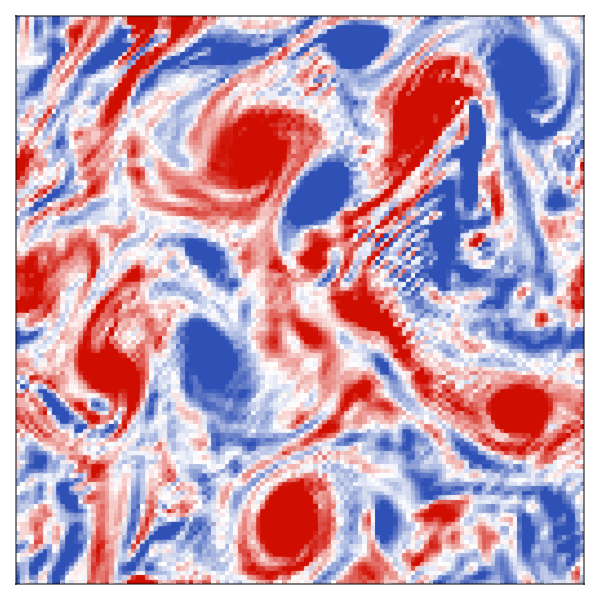}
     \label{std-ef-t1}}
     \hspace{0.5cm}
\subfloat{\includegraphics[width=0.04\linewidth]{images/decayingturbulence/colorbar.png}}
    \caption{Vorticity fields at $t=1$ for different methodologies in a test configuration, in decaying turbulence test case.}
    \label{fig:fields-t1}
\end{figure}

After building the filter, we perform the \emph{online} EFR simulations.
In particular, Figure \ref{fig:test-general} shows the performance of the method in terms of the total kinetic energy and enstrophy trend in time, and of the energy spectrum at fixed time $t=1$ s. In Figures \ref{fig:fields-t03} and \ref{fig:fields-t1} we can see the corresponding vorticity fields at fixed times $t=0.3$ and $t=1$ s.

From the plots we can draw the following conclusions:
\begin{itemize}
    \item \ddEF{} and \energyEFR{} outperform all existing methods in terms of kinetic energy, global enstrophy, and spectrum. The \ddEF{} and \energyEFR{} vorticity fields at initial times are close to the filtered DNS in local features (Figure \protect\ref{dd-efr-t03} and \protect\ref{energy-dd-efr-t03}), while being \emph{qualitatively similar} to the reference at later times (Figure \protect\ref{dd-efr-t1} and \protect\ref{energy-dd-efr-t1}).
    \item The \ddEF{} and \energyEFR{} methods exhibit comparable performance when trained in the full data regime. Using multiple random initializations and a long training time window helps stabilize the filter’s effect. In fact, this approach inherently yields decaying energy and enstrophy, without the need for explicit constraints.
    \item The most accurate methods in the spectrum are the DD-EF(R) approaches and standard EF (at fixed $\delta=\delta_{\text{opt}}$).
    \item All state-of-the-art methods appear overdiffusive (as seen in the analysis in Figure \ref{fig:plots-literature}), while the data-driven filter successfully retrieves the dissipated energy. From a qualitative point of view, the vorticity fields look all close to the reference at initial times (Figure \protect\ref{smag-t03}, \protect\ref{std-efr-t03}, and \protect\ref{std-ef-t03}), but overdiffusive at later time $t=1$ s (Figure \protect\ref{smag-t1}, \protect\ref{std-efr-t1}, and \protect\ref{std-ef-t1}).
    \item The noEFR simulation exhibits large enstrophy values in the time window $[0, 3]$ s.
    This is also reflected in the presence of noise and wiggles in the vorticity field (Figures \protect\ref{noefr-t03} and \protect\ref{noefr-t1}).
    As time evolves and turbulence is dissipated, the noEFR enstrophy decreases and closely follows the filtered DNS behaviour. However, the noEFR energy is the closest one to the filtered DNS in the global kinetic energy with respect to traditional filters. Therefore, relying solely on global kinetic energy as a metric can be misleading, as it may conceal deficiencies in accurately capturing the flow’s fine-scale structures and local dynamics, which are better reflected in enstrophy and vorticity-related quantities.
\end{itemize}

\begin{figure}[htpb!]
    \centering
    \subfloat{\includegraphics[height=80pt]{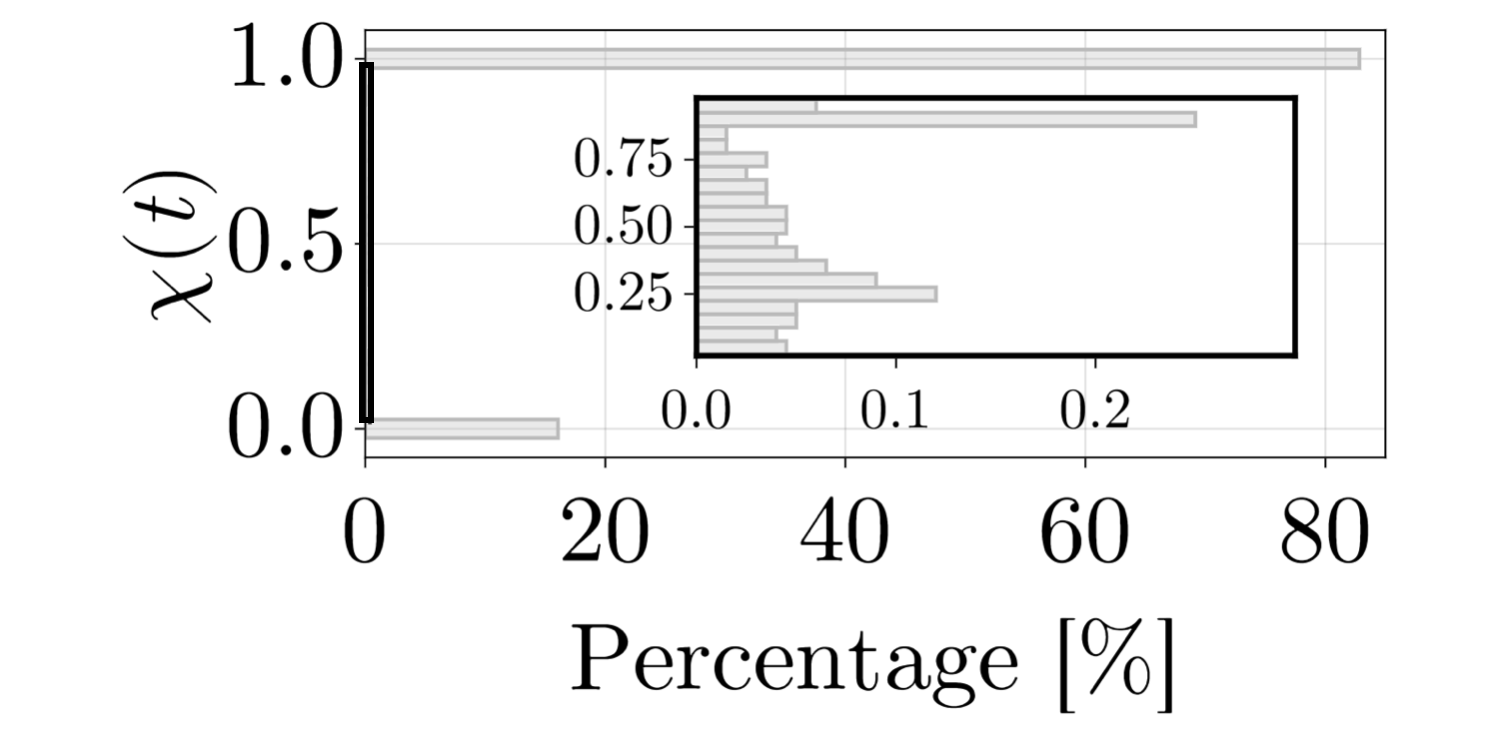}}
    \subfloat{\includegraphics[height=80pt]{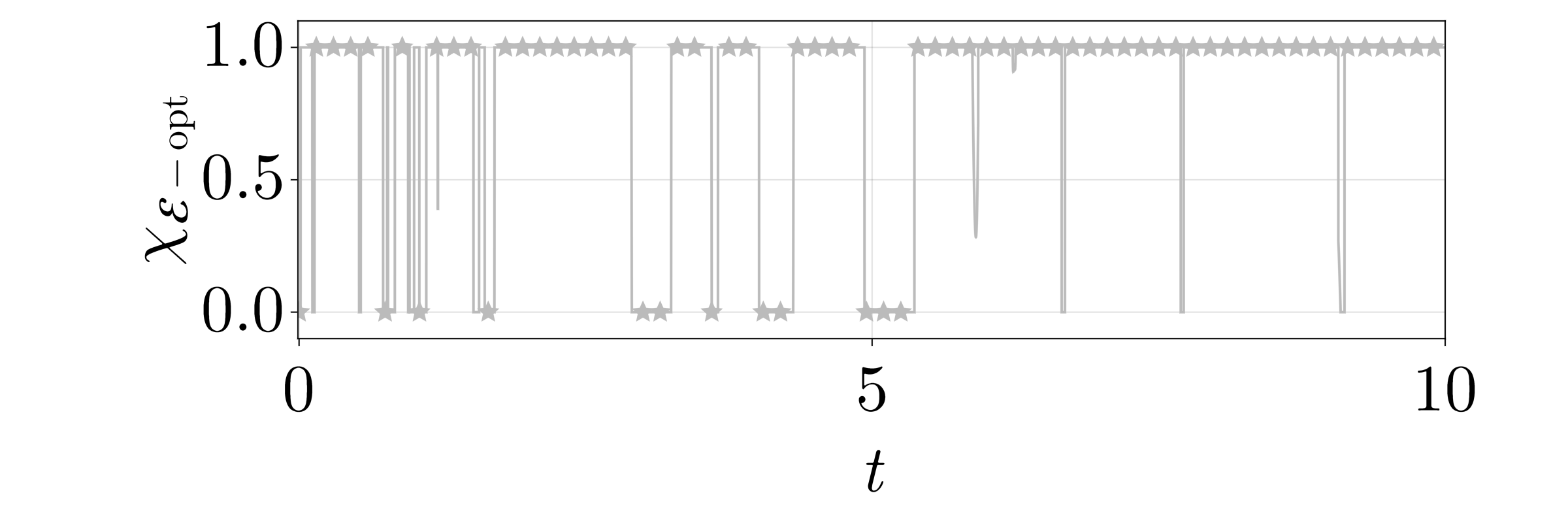}}
    \caption{Distribution and time history of the optimal relax parameter \chienergy{} for a test initialization. The plots refer to the case $\Ttrain=3$ s, $\Itrain=10$ in decaying turbulence test case.}
    \label{fig:chi-seed-11}
\end{figure}

It is also useful to analyse the online time evolution of the optimal relaxation parameter $\chienergy$. Figure \ref{fig:chi-seed-11} shows the distribution and time history of $\chienergy$ in the case of $\Ttrain=3$ s, for a specific test initialization. $\chienergy$ is mostly oscillating between its admissible lower and upper bounds, namely $0$ and $1$.
As pointed out in Section \ref{sec:tuning-chi-energy}, the enforcement of the energy constraint implies that, whenever $\tilde{b}_n>0$, we impose $\chienergy=0$, which is the only admissible value satisfying the constraint and corresponds to noEFR.
Additionally, we can see that the most frequent value ($\sim 85\%$) is $\chienergy=1$, which coincides to \ddEF{}. This confirms the similar performances of the methods observed in Figure \ref{fig:test-general}.

When analyzing the qualitative features at later times (Figure \ref{fig:fields-t1}), we observe more coherent vortex structures in \ddEF{} and \energyEFR{} with respect to state-of-the-art approaches. This is confirmed by the spectral analysis at $t=1\, s$ in Figure \ref{fig:test-general}. \ddEF{} and \energyEFR{} models show a better match with the filtered DNS in terms of energy at intermediate-to-high wavenumbers, indicating reduced numerical dissipation and improved preservation of vortex structures. Smagorinsky retains excessive energy at the large wavenumbers, while standard EFR exhibits a steeper decay in the spectrum at intermediate wavenumbers.

Additional results about the data-driven approaches in the full data regime in case of smaller training time window ($\Ttrain=1$ s), are available in the supplementary material (\ref{supplementary-case1}).

\subsubsection{Data-driven EFR in the case of scarce data}
\label{results-scarce-dd-efr-case1}

We now analyze the performance of the methods in the case of scarce data, namely when the filter is computed from a single random initialization ($\Itrain=1$).

\begin{figure}[htpb!]
    \centering
    \includegraphics[width=\linewidth]{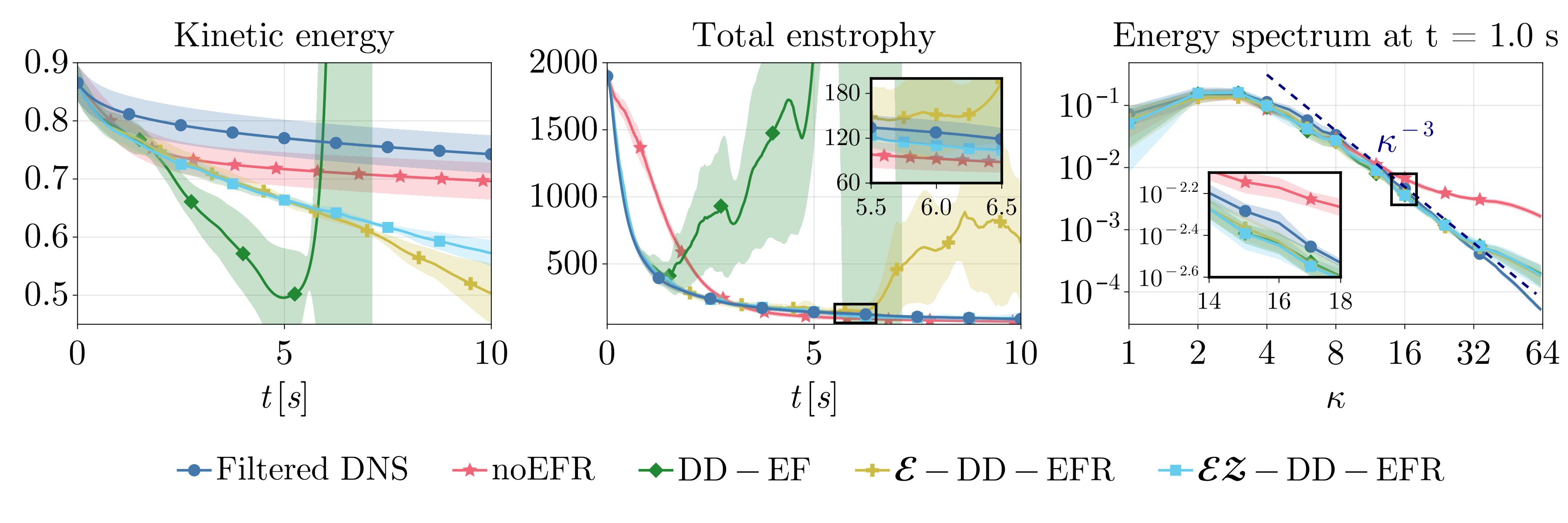}
    \caption{Time evolution of total kinetic energy, total enstrophy, and spectrum of the kinetic energy at time $t=1$ s, for the proposed DD-EF(R) methodologies and for state-of-the-art approaches.
    In particular, we show the case $\Ttrain=1$ s and $\Itrain=1$.
    For each method, the Figure shows the average values (solid lines) and the 95$\%$ confidence interval among 5 test configurations, in decaying turbulence test case.}
    \label{fig:scarce-test-general}
\end{figure}

Figure \ref{fig:scarce-test-general} presents the statistics for kinetic energy, enstrophy, and energy spectrum at fixed time, and we focus on the case $\Ttrain=1$ s, which is the most challenging one.

Figure \ref{fig:scarce-test-general} highlights the following aspects:
\begin{itemize}
    \item \ddEF{} is unstable both in terms of energy and enstrophy, which increase and blow up after $t=5$ s. It confirms the fact that the scarce-data-filter can be unstable, as it can also inject energy at medium wave numbers, as previously showed in Figure \ref{fig:fstars}.
    \item \energyEFR{}, which satisfies the energy dissipation constraint, presents a stable decreasing energy behaviour, but it is inaccurate in the global enstrophy after $t=6$ s.
    \item \enstrophyEFR{} is stable both in the energy and in the enstrophy. It provides an accurate approximation of the filtered DNS in the enstrophy. However, both \energyEFR{} and \enstrophyEFR{} lead to overdissipative results in the energy at large times. This is a consequence of using too little training data for the filter construction.
\end{itemize}

\begin{figure}[htpb!]
    \centering
    \subfloat{\includegraphics[height=80pt, trim={0 0 1cm 0}, clip]{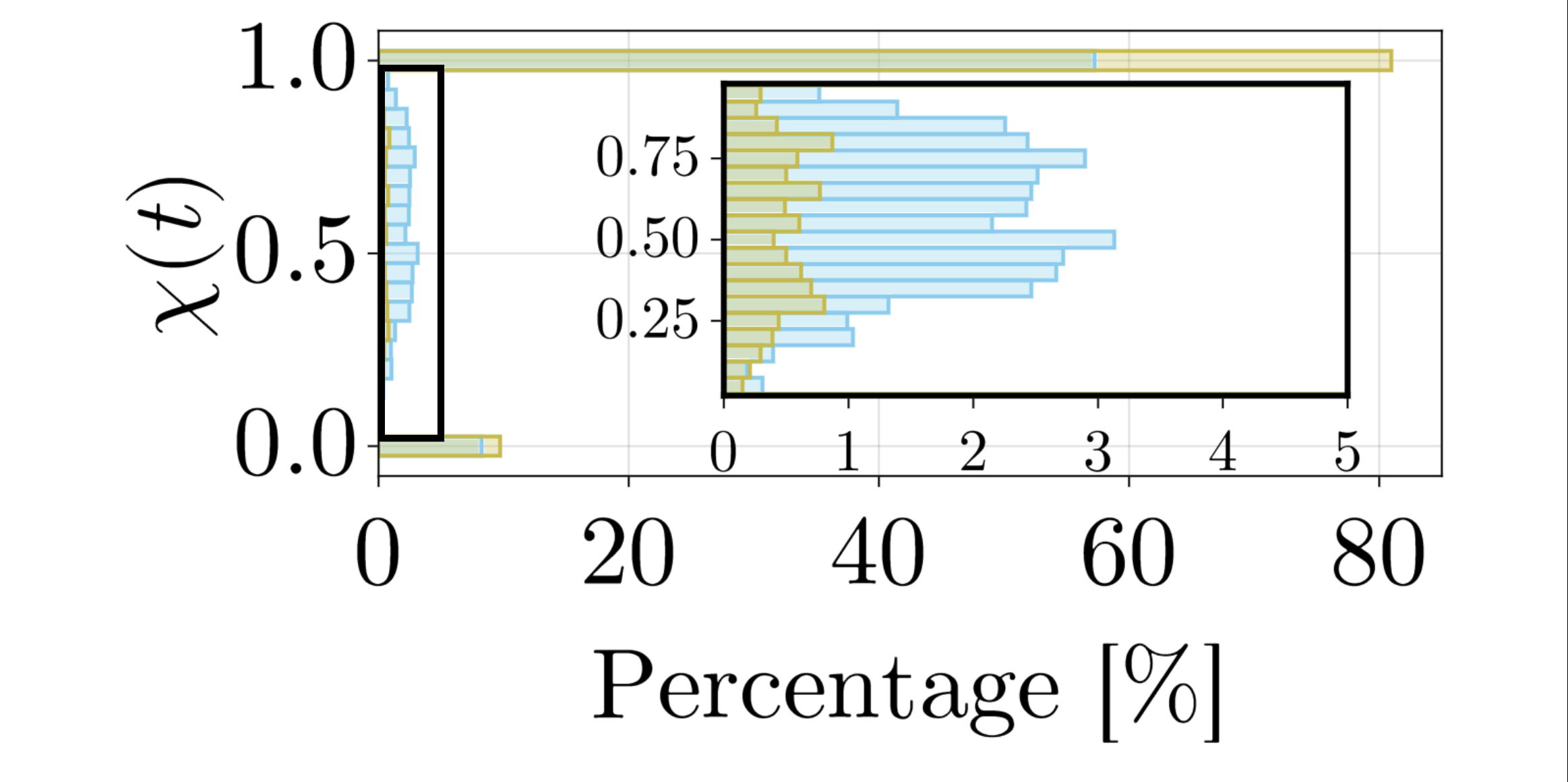}}
    \subfloat{\includegraphics[height=80pt]{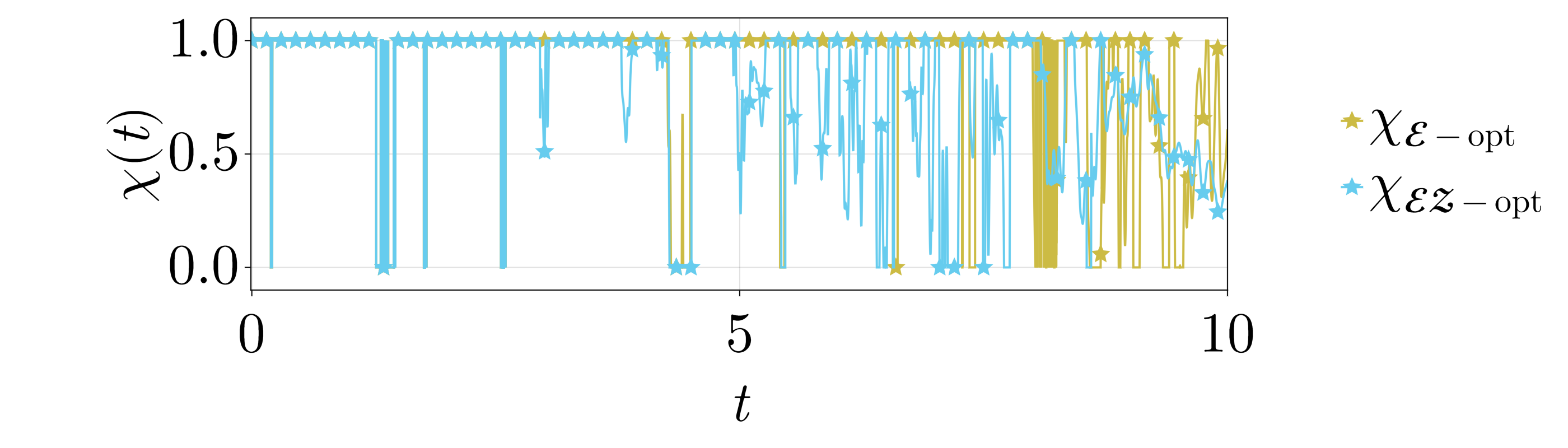}}
    \caption{Distribution and time history of the optimal relax parameter \chienergy{} and \chienstrophy{} for a test initialization. The plots refer to the case $\Ttrain=1$ s, $\Itrain=1$ in decaying turbulence test case.}
    \label{fig:scarce-chi-seed-11}
\end{figure}

Figure \ref{fig:scarce-chi-seed-11} illustrates the time evolution of $\chienergy$ and $\chienstrophy$ for \reviewerB{a test} initialization. $\chienergy$ exhibits stronger temporal oscillations with respect to the previous full data-case. Indeed, the use of a reduced amount of data leads to increasing instability and error of the filtering approach in terms of kinetic energy and enstrophy, as previously noticed from the analysis of the filter in Figure \ref{fig:fstars}.

Moreover, since \chienstrophy{} fulfills both energy and enstrophy dissipation constraints, it has more temporal variability. $\chienstrophy$ also assumes more values within the interval $(0, 1)$ especially at $t>5$ s, while $\chienergy$ mostly converges to either $0$ (noEFR) or $1$ (\ddEF).

We can conclude that the scarce-data filter may lead to unstable results due to the small amount of data, but the imposition of the energy and/or enstrophy constraints \reviewerB{(Equations \eqref{eq:condition-energy} and \eqref{eq:condition-enstrophy}, respectively)} stabilizes the simulation. It is worth remarking that such constraints are only imposed \emph{online} during the simulation by modifying the relax step, and without affecting the filter itself. This highlights the \emph{flexibility} of the approach, allowing stability to be enforced without altering the underlying filter design.

Additional results on data-driven approaches in the scarce data regime in the case of larger training time window ($\Ttrain=3$ s), are available in the supplementary material (\ref{supplementary-case1}).

\medskip

Given the observed differences in performance between filters trained with full and scarce data, the questions we will try to address are:
\emph{(i)} ``\emph{How do the scarce-data and full-data filters compare in terms of overall accuracy across different physical metrics?}", and \emph{(ii)} ``\emph{To what extent can constraints bridge the gap in terms of overall error?}".

To answer these questions, we analyse the following error metrices.

\begin{itemize}
    \item Relative errors in global energy and enstrophy averaged in time, computed as:
    \begin{equation}
    \begin{split}
        &\text{err}_{\mathcal{E}}(\boldsymbol{\mu}^{(i)})=\frac{1}{N_{t=3}} \sum_{n=1}^{N_{t=3}}\left|\dfrac{\mathcal{E}^n_h(\boldsymbol{\mu}^{(i)})-\mathcal{E}^n_{\text{ref}}(\boldsymbol{\mu}^{(i)})}{\mathcal{E}^n_{\text{ref}}(\boldsymbol{\mu}^{(i)})}\right|, \\ &\text{err}_{\mathcal{Z}}(\boldsymbol{\mu}^{(i)})=\frac{1}{N_{t=3}} \sum_{n=1}^{N_{t=3}}\left|\dfrac{\mathcal{Z}^n_h(\boldsymbol{\mu}^{(i)})-\mathcal{Z}^n_{\text{ref}}(\boldsymbol{\mu}^{(i)})}{\mathcal{Z}^n_{\text{ref}}(\boldsymbol{\mu}^{(i)})}\right|,
        \end{split}
        \label{eq:errs-energy-enstrophy}
    \end{equation}
    where $\boldsymbol{\mu}^{(i)}$ represents the $i$-th test random initialization, $\mathcal{E}_{\text{ref}}^n$ and $\mathcal{Z}_{\text{ref}}^n$ are the energy and enstrophy of the filtered DNS solution at time step $t_n$.
    \item Absolute logarithmic errors in the energy spectrum averaged in the wave numbers and in time, computed as:
    \begin{equation}
        \text{err}_{\text{spectrum}}=
        \frac{1}{N_{t=3}} \sum_{n=1}^{N_{t=3}}
        \frac{1}{K} \sum_{\kappa=1}^{K}
\log_{10}{\left( \frac{\mathcal{E}^n_h(\kappa)(\boldsymbol{\mu}^{(i)})}{\mathcal{E}^n_{\text{ref}}(\kappa)(\boldsymbol{\mu}^{(i)})} \right)}.
        \label{eq:errs-spectrum}
    \end{equation}
\end{itemize}

The value $N_{t=3}$ in expressions \eqref{eq:errs-energy-enstrophy} and \eqref{eq:errs-spectrum} is the number of time steps corresponding to $t=3$ s. We choose the time window
$[0,3]$ for error computation because in the first test case the energy dissipates significantly at later times, and the flow loses its turbulent characteristics.

\begin{figure}[htpb!]
    \centering
    \subfloat{\includegraphics[height=125px]{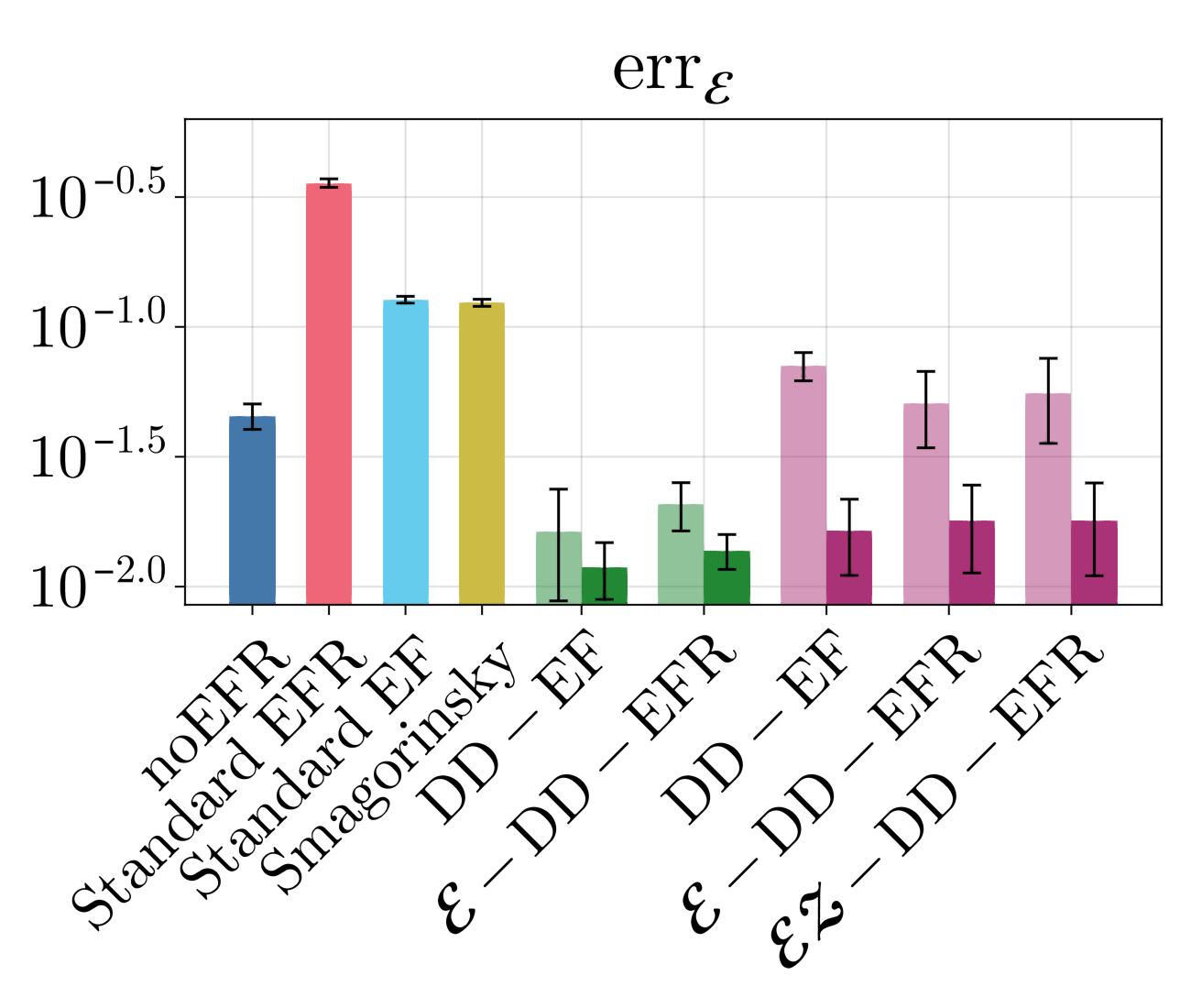}}
    \subfloat{\includegraphics[height=125px]{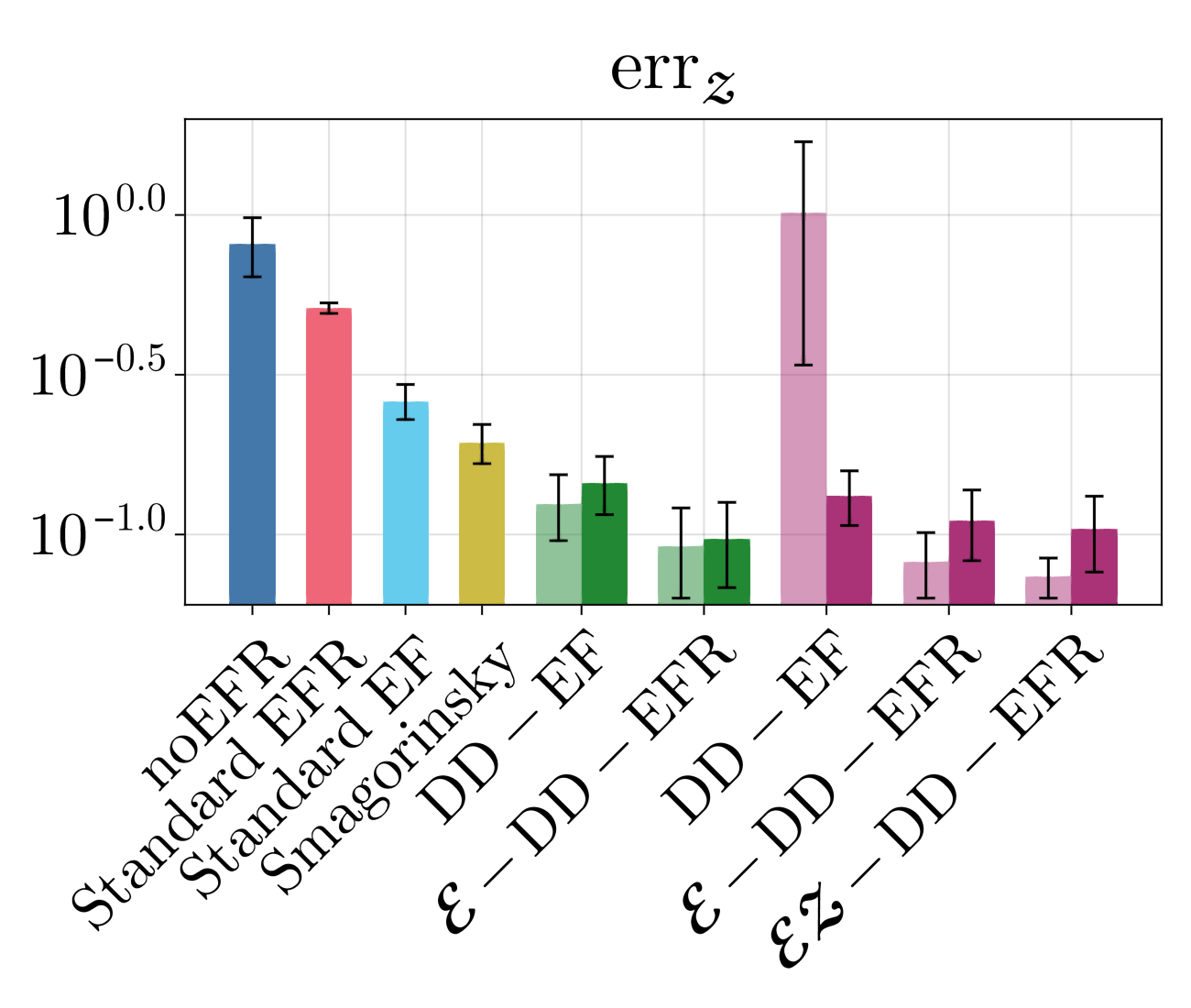}}
   \subfloat{\includegraphics[height=125px]{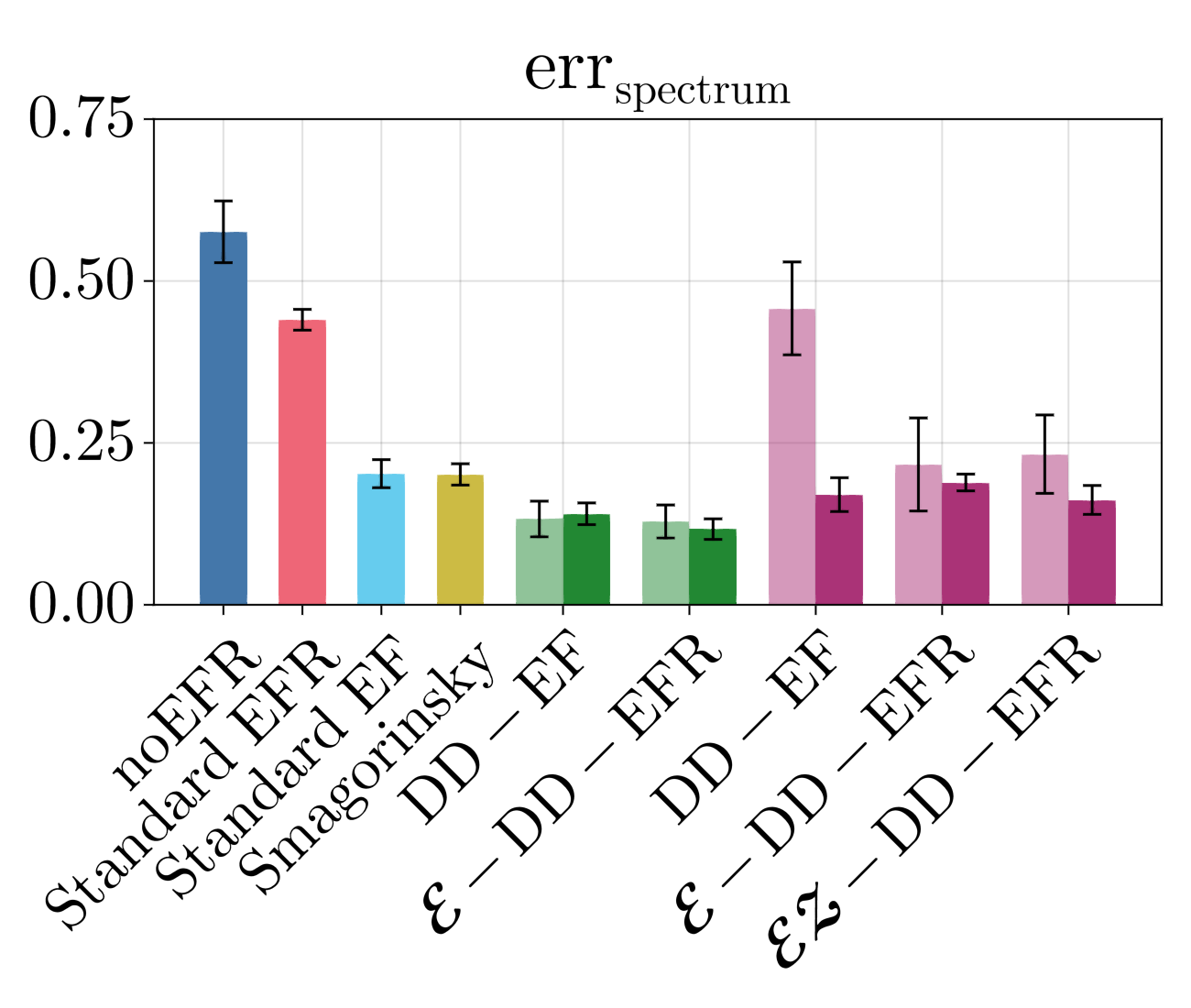}}\\
   \subfloat{\includegraphics[height=27px]{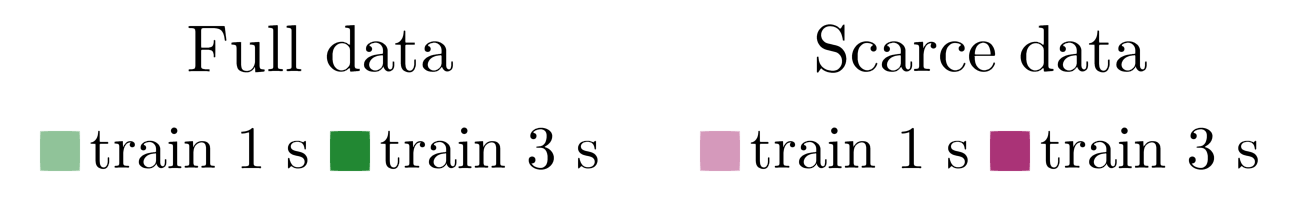}}
    \caption{Relative error of the total kinetic energy, total enstrophy, and absolute error in logarithmic scale of the spectrum, of the different methodologies with respect to the filtered DNS. The errors are averaged in the time interval $[0, 3]$, and then evaluated as average values and confidence intervals (error bars) among $\Itest=5$ test configurations. This Figure refers to decaying turbulence test case.}
    \label{fig:errs-all-3s}
\end{figure}

The plots in Figure \ref{fig:errs-all-3s} display the above-mentioned errors in their average value and confidence interval among $\Itest=5$ initializations.

From Figure \ref{fig:errs-all-3s}, we can draw the following conclusions:
\begin{itemize}
    \item DD-EF(R) approaches with $\Ttrain=3$ s are the best performing methods in all the metrics considered.
    \item Among the data-driven methods, \ddEF{} and \energyEFR{} with $\Ttrain=3, \Itrain=10$ (full data) reach the highest precision with respect to the filtered DNS in terms of global energy and spectrum.
    \item As previously noticed from Figure \ref{fig:scarce-test-general}, \energyEFR{} and \enstrophyEFR{} when $\Ttrain=1, \Itrain=1$ (scarce data) significantly improve the \ddEF{} results, especially in the enstrophy. On the other hand, only a slight improvement is observed when $\Ttrain=3, \Itrain=1$, likely because the filter, built from a larger dataset, exhibits greater stability.
    \item In both full and scarce data, we obtain the largest accuracy in energy and spectrum when $\Ttrain=3$ s, as it aligns with the increased amount of data. Hovever, we have more accurate results in enstrophy when $\Ttrain=1$ s. A possible explanation is that the system exhibits more chaotic behavior in the initial phase, while dissipation dominates as time progresses. Therefore, a filter built using data from a shorter time window may be more accurate in capturing gradients compared to those constructed over longer time intervals.
\end{itemize}

Hence, we can conclude that the answers to the previous questions are:\\
\emph{(i)} In general, full-data filters are more accurate in terms of kinetic energy and energy spectrum.\\
\emph{(ii)} In the scarce-data regime, the constraints significantly improve the global enstrophy accuracy, reaching even better results than full-data regimes.

\subsection{Test case 2: two-dimensional Kolmogorov flow}
\label{test-case-2}

The second test case has a similar setting to the first one, but in the presence of a forcing term.
To prevent energy decay during long simulations, we add a Kolmogorov-type body force to inject energy into the system at specific wave numbers, as in \cite{Agdestein2025MLLES, chandler2013invariant, list2025differentiability, kochkov2021MLLES}. This force is only characterized by an horizontal component
$$f_x(y) = c_f\sin{(2\pi k_f y)},$$ where $k_f=4$ is the wave number at which the energy is injected. The forcing amplitude $c_f$ is selected to ensure that the system attains a statistically stationary energy level sufficiently high to sustain turbulence without leading to unbounded energy growth. In our case $c_f=0.65$.

To test extrapolation of our approach across test cases, in the DD-EF(R) we keep the filter matrix $\hat{\boldsymbol{\mathcal{F}}}^{\star}$ computed in the first test case, effectively performing \emph{extrapolation in this test scenario} as well.

For the state-of-the-art methods we also consider the optimal parameters computed in the first test case.

\subsubsection{Data-driven EFR in the case of full data}
\label{results-dd-efr-case2}

As done for the first test case in Section \ref{results-dd-efr-case1}, we dedicate this Section to discuss the \ddEF{} and \energyEFR{} approaches in full data regime ($\Itrain=10$). For the sake of brevity we only consider the case $\Ttrain=3$ s.

\begin{figure}[htpb!]
    \centering
    \includegraphics[width=\linewidth]{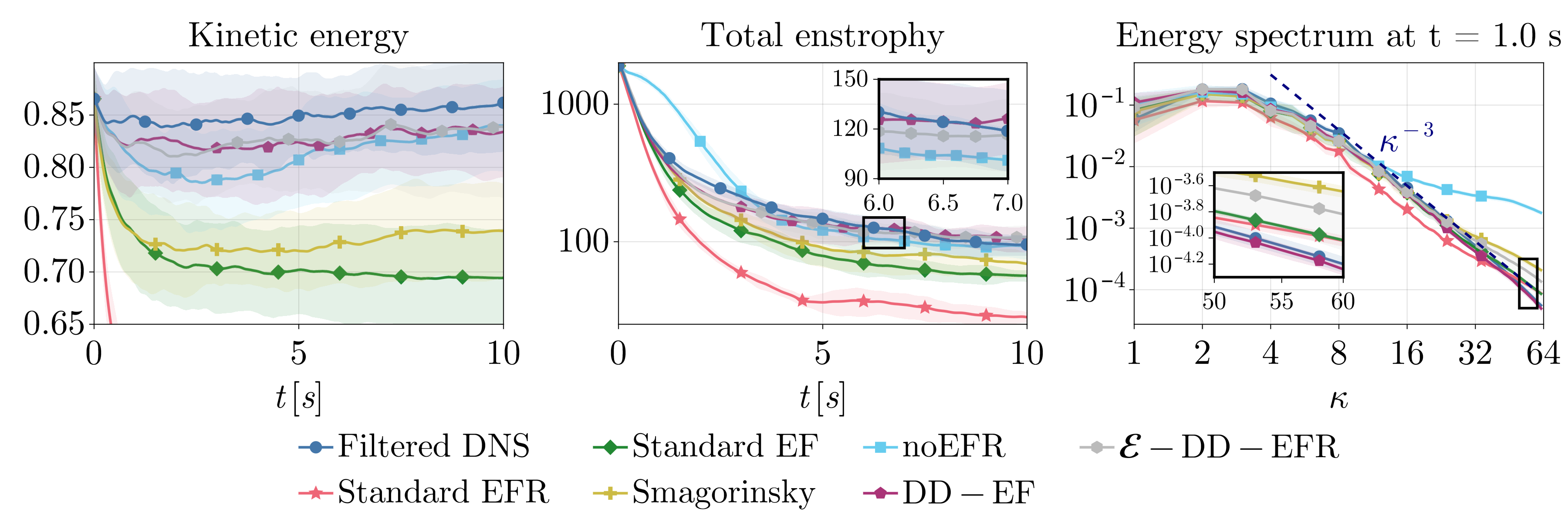}
    \caption{Time evolution of total kinetic energy, total enstrophy, and spectrum of the kinetic energy at time $t=1$ s, for the proposed DD-EF(R) methodologies and for state-of-the-art approaches.
    In particular, we show the case $\Ttrain=3$ s, $\Itrain=10$.
    For each method, the Figure shows the average values (solid lines) and the 95$\%$ confidence interval among 5 test configurations, in the Kolmogorov test case.}
    \label{fig:case2-test-general}
\end{figure}

In particular, we analyse the time evolution of the global energy and enstrophy, and the spectrum at fixed time $t=1$ s, comparing such metrics with the filtered DNS reference (Figure \ref{fig:case2-test-general}).

The Figure shows that:
\begin{itemize}
    \item The \ddEF{} and \energyEFR{} provide the best-performing solutions in all metrics considered. While both have similar accuracy in energy and enstrophy, \energyEFR{} is closer to the reference in the spectrum, especially at large wave numbers.
    \item As expected, the reference kinetic energy in the Kolmogorov case remains stable in time, while state-of-the-art methods provide overdiffusive solutions, as in the first test case.
    \item The reference global enstrophy has a similar trend as before, mostly decreasing in time. As for the previous test case, the noEFR solution is characterized by larger gradients and enstrophy at initial time, which is dissipated at $t>3$ s.
\end{itemize}

The time evolution of $\chienergy$ in the \energyEFR{} approach is similar to the one obtained in the previous test case (Figure \ref{fig:chi-seed-11}). The same holds for the qualitative behavior of the solutions, that the reader can find in the supplementary material (\ref{supplementary-case2}).

\subsubsection{Data-driven EFR in the case of scarce data}
\label{results-scarce-dd-efr-case2}

In this part, we show the results when training the filter with scarce data (one random initialization), and we focus on $\Ttrain=3$ s.

\begin{figure}[htpb!]
    \centering\includegraphics[width=\linewidth]{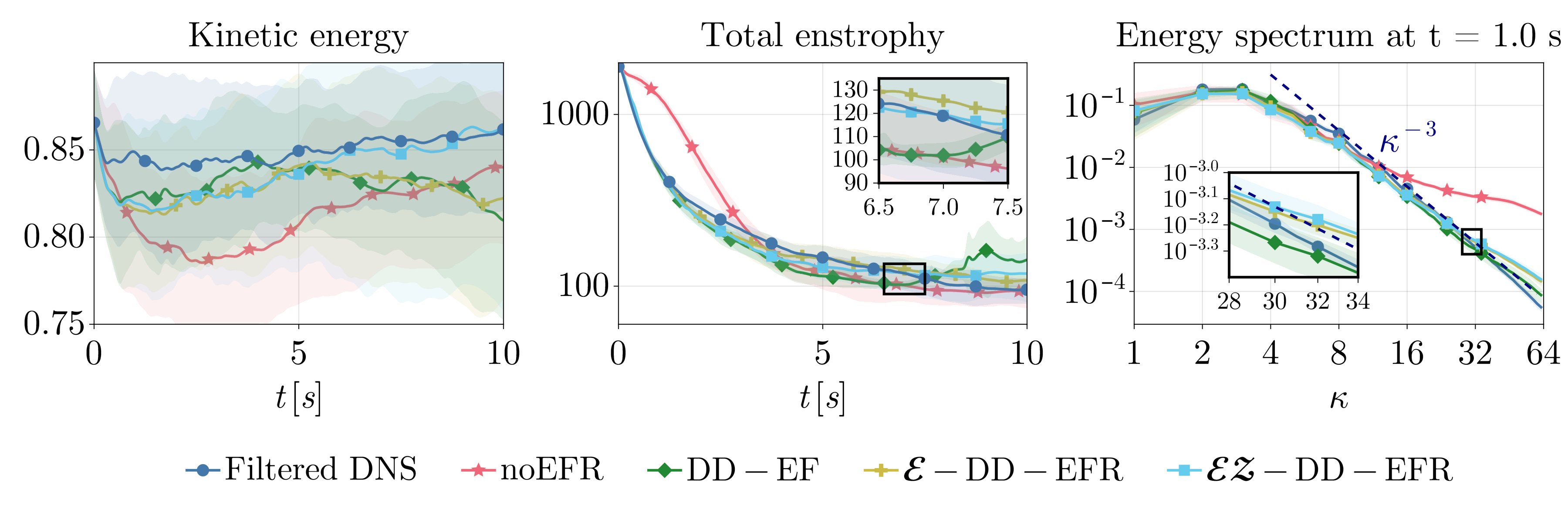}
    \caption{Time evolution of total kinetic energy, total enstrophy, and spectrum of the kinetic energy at time $t=1$ s, for the proposed DD-EF(R) methodologies and for state-of-the-art approaches.
    In particular, we show the case $\Ttrain=3$ s, $\Itrain=1$.
    For each method, the Figure shows the average values (solid lines) and the 95$\%$ confidence interval among 5 test configurations, in the Kolmogorov test case.}
    \label{fig:case2-scarce-test-general}
\end{figure}

Differently from what noticed in Section \ref{results-scarce-dd-efr-case1}, where we considered the case $\Ttrain=1$ s, Figure \ref{fig:case2-scarce-test-general} shows more accurate and less dissipative DD-EF(R) solutions. This is consistent with the increased amount of time instances used to compute the filter.
More in detail, \ddEF{} and \energyEFR{} have similar energy and enstrophy time evolution. However, \ddEF{} exhibits instabilities and inaccuracies in the global enstrophy, especially at $t>6$ s. The addition of the energy and enstrophy constraints stabilizes the enstrophy at large time and better matches the energy spectrum at medium wave numbers.
The further addition of the enstrophy constraint improves the accuracy in the kinetic energy, especially after $t=5$ s.

\begin{figure}[htpb!]
    \centering
    \subfloat{\includegraphics[height=80pt]{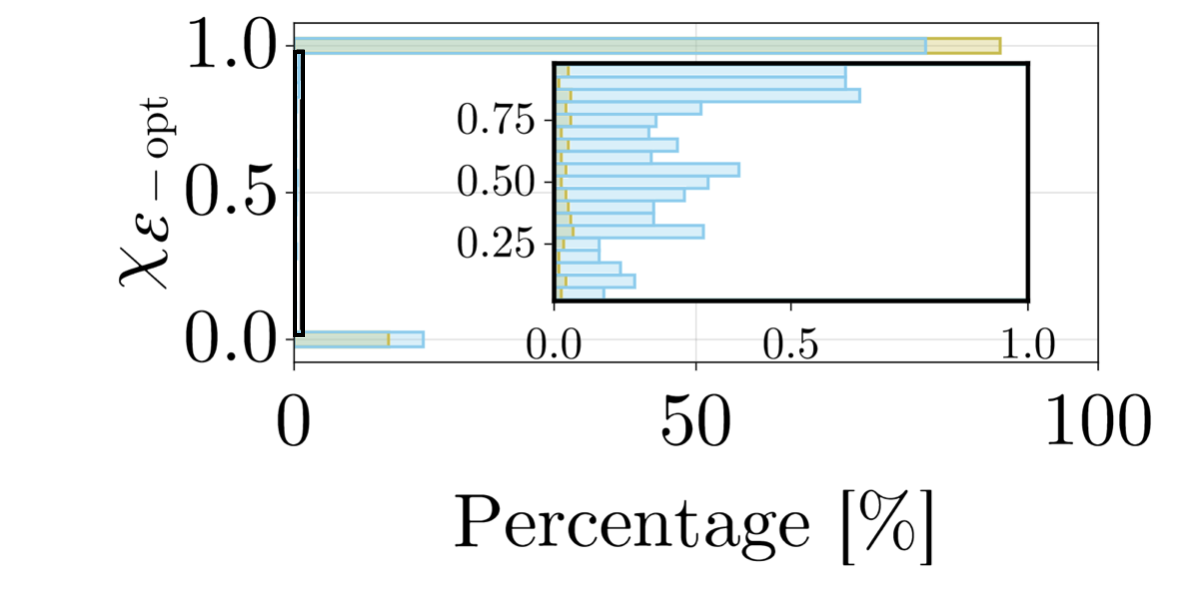}}
    \subfloat{\includegraphics[height=80pt]{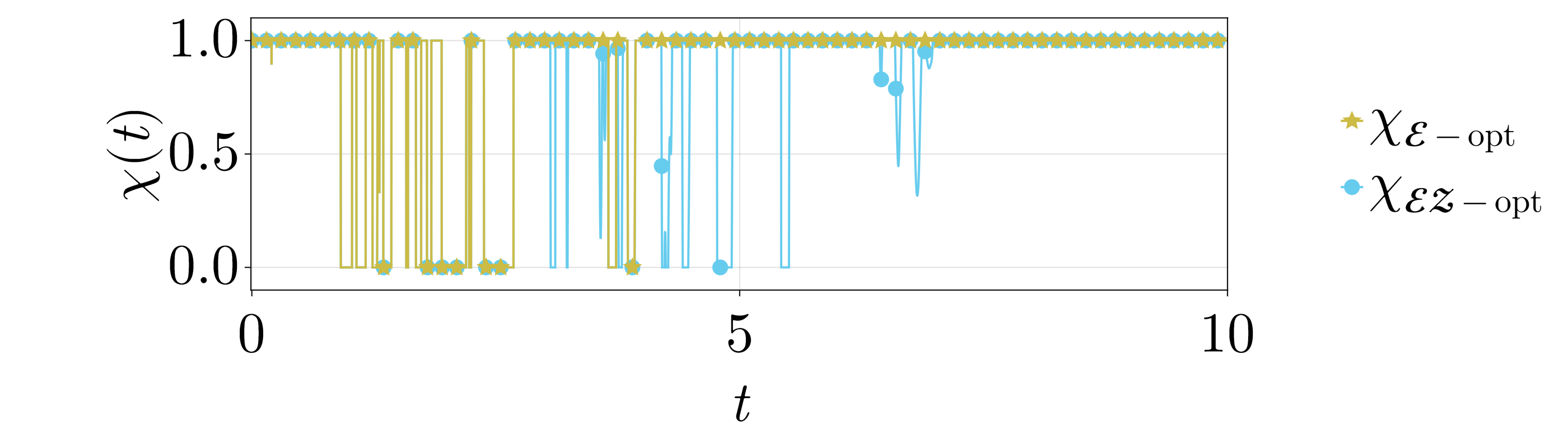}}
    \caption{Distribution and time history of the optimal relax parameter \chienergy{} and \chienstrophy{} for a test initialization. The plots refer to the case $\Ttrain=3$ s, $\Itrain=1$ in the Kolmogorov test case.}
    \label{fig:case2-scarce-chi-seed-11}
\end{figure}

The time evolution of the optimal $\chienergy$ and $\chienstrophy$ is displayed in Figure \ref{fig:case2-scarce-chi-seed-11}. As in Figure \ref{fig:scarce-chi-seed-11}, in both cases the optimal value is mostly $\chi=1$, but $\chienstrophy$ has more oscillations in time, especially at large simulation time.

From the comparison between Figure \ref{fig:scarce-chi-seed-11} ($\Ttrain=1$ s) and \ref{fig:case2-scarce-chi-seed-11} ($\Ttrain=3$ s), we can notice that a larger training window leads to increased stability, reduced oscillations in the optimal $\chi$ and fewer intermediate values.

Based on the observed results, we can conclude that the novel filtering approach, originally computed on a different configuration, retains its effectiveness in the current test case, both in full and scarce data regimes. This validates the method’s capability to \emph{extrapolate} and \emph{adapt} to qualitatively similar flow conditions.

\medskip

\begin{figure}[htpb!]
    \centering
    \subfloat{\includegraphics[height=125px]{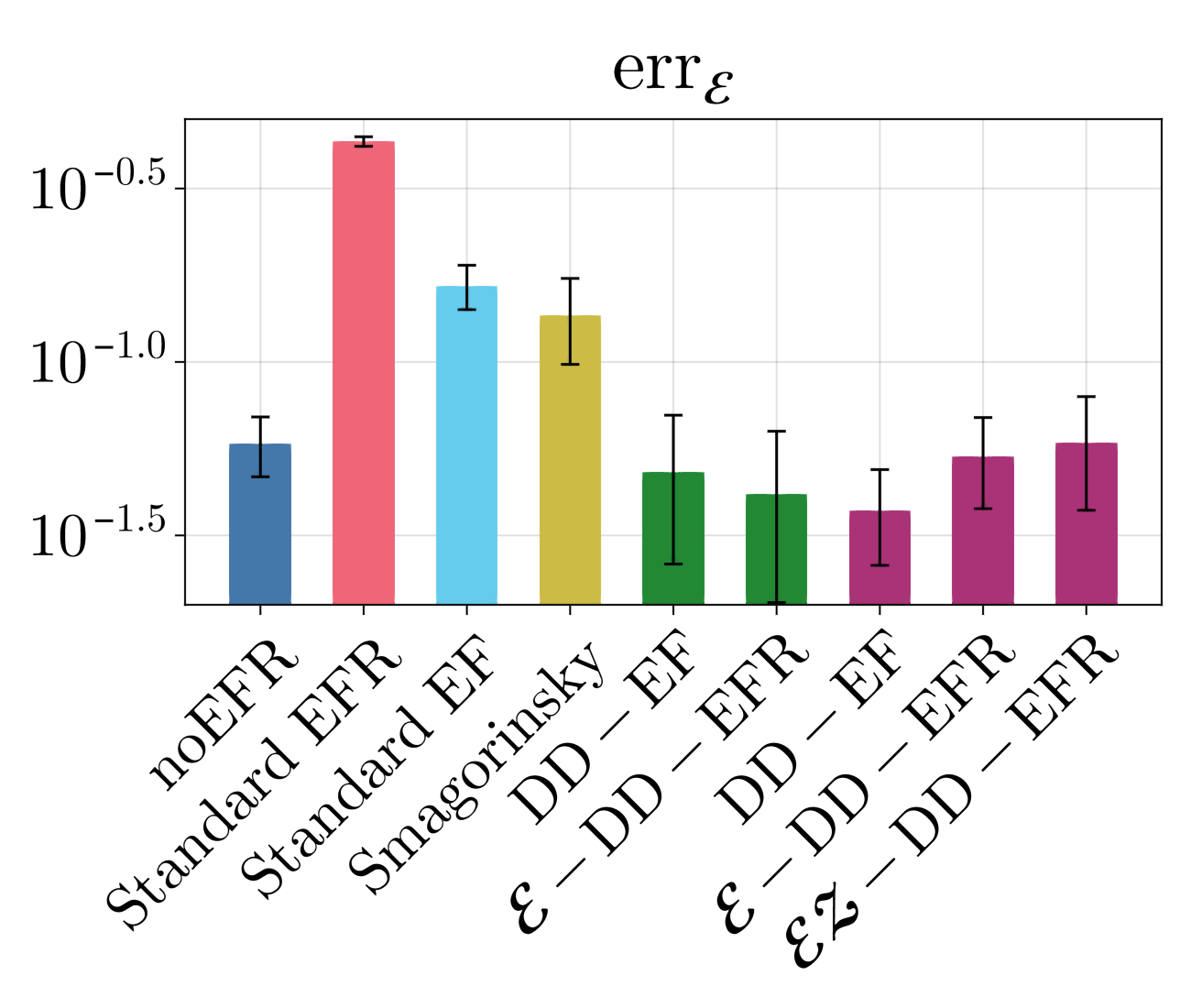}}
    \subfloat{\includegraphics[height=125px]{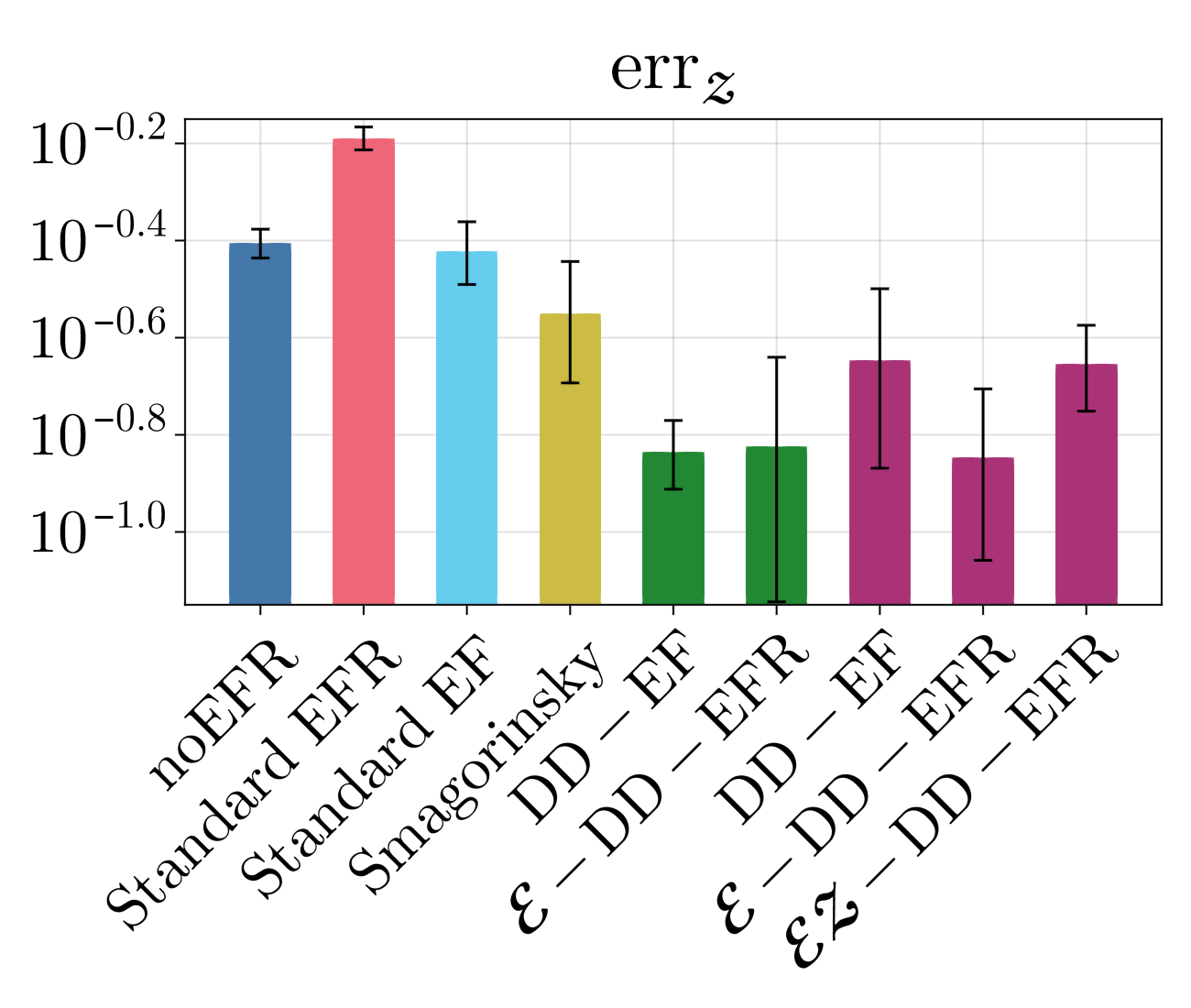}}
   \subfloat{\includegraphics[height=125px]{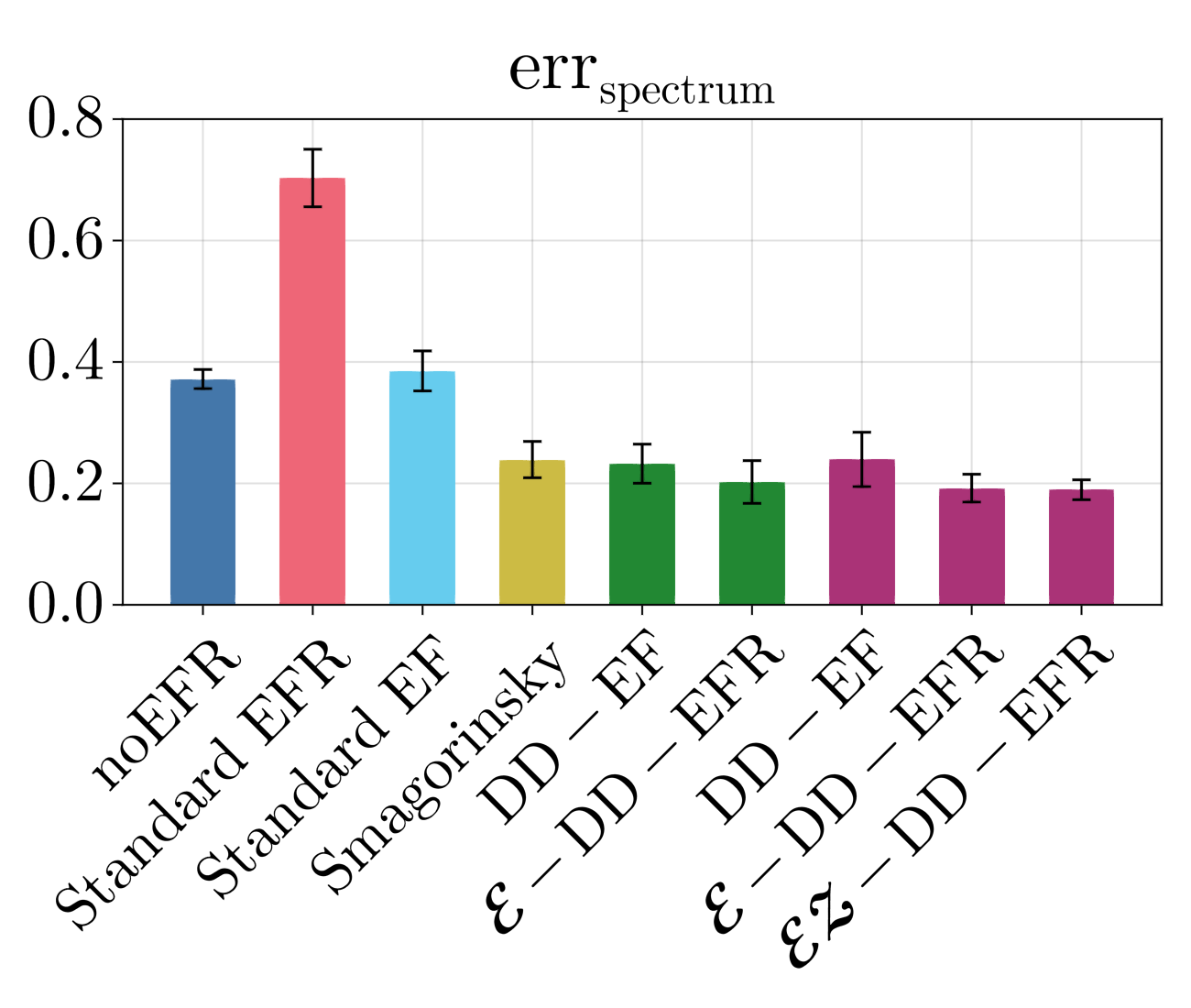}}\\
   \subfloat{\includegraphics[height=20px, trim={0 1cm 0 1cm}, clip]{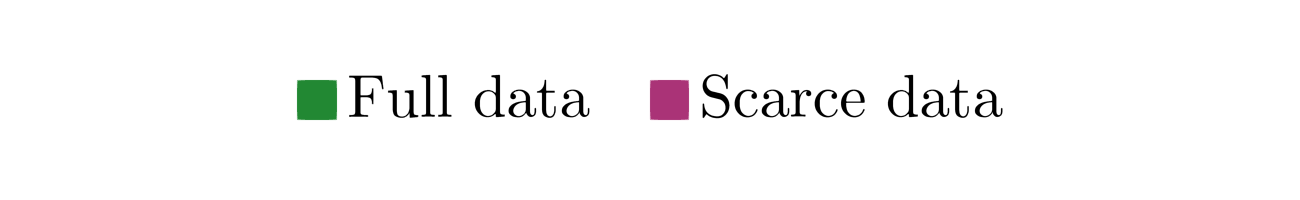}}
    \caption{Relative error of the total kinetic energy, total enstrophy, and absolute error in logarithmic scale of the spectrum, of the different methodologies with respect to the filtered DNS. The errors are averaged in the time interval $[0, 10]$, and then evaluated as average values and confidence intervals (error bars) among $5$ test configurations. This Figure refers to the Kolmogorov test case.}
    \label{fig:case2-errs-all-10s}
\end{figure}

As done for first test case, we now compare the full and scarce data regimes in terms of global errors (expressions \eqref{eq:errs-energy-enstrophy} and \eqref{eq:errs-spectrum}).
Figure \ref{fig:case2-errs-all-10s} depicts the global error analysis for the Kolmogorov test case, which is similar to the one showed in \ref{fig:errs-all-3s}.
We specify that the total number of time steps considered in this case for error computation is $N$ instead of $N_{t=3}$, and we only show the errors related to the best-performing DD-EF(R) strategies (namely, those with $\Ttrain=3$ s).

The main difference with respect to Figure \ref{fig:errs-all-3s} lies in the noEFR accuracy, which improves when considering $N$ instead of $N_{t=3}$. As it no longer exhibits spurious oscillations and noisy behavior at later times, it is qualitatively closer to the reference solution and has smaller error.

\subsection{Computational time considerations}

In this Section, we discuss the computational cost required to construct the proposed methods and compare their efficiency with that of state-of-the-art approaches \footnote{The simulations have been performed on an Apple M2 chip (8-core), 16GB unified memory.}.

The offline and online wall times are reported in Table \ref{tab:cpu-time}. We call \emph{offline} wall time the time needed to perform the parameters' optimization in the case of state-of-the-art approaches (standard EFR/EF and Smagorinsky), and the time used to compute the filter matrix $\hat{\boldsymbol{\mathcal{F}}}^{\star}$ in DD-EF(R) methods.
We refer to \ref{literature-tuning} for more details on the parameters' tuning.

The offline wall time in DD-EF(R) includes:
\begin{itemize}
    \item the time needed to evolve the filtered DNS solutions of one time step, namely to obtain the snapshot matrix $W^{\text{true}}$;
    \item the time needed to solve the least squares minimization problem and find the filter components as in Equation \eqref{eq:f_star}.
\end{itemize}

We refer to \emph{online} wall time as the time needed to perform the simulation at $T=3$ s, with fixed $\Delta t=\num{5e-4}$ s, to have a fair comparison among the methods.

As from Table \ref{tab:cpu-time}, the offline time needed to perform parameter tuning in state-of-the-art methods may be large, depending on the optimization strategy and the amount of data considered in the optimization. Here, we used an approximate gradient descent, although one may also rely on literature-based parameter choices to avoid this step.

On the other hand, the cost to compute the data-driven filter matrix is relatively small, as it only involves the resolution of a least-squares problem. We remark that the filter is the same in \ddEF{}, \energyEFR{} and \enstrophyEFR{} approaches and the offline time only depends on the given amount of data (\Ttrain{} and \Itrain{}) used to compute it, as specified in Table \ref{tab:cpu-time}.

\begin{table}[htpb!]
    \centering
    \begin{tabular}{ccc}
    \toprule
       {\textbf{Method}}&{\textbf{Offline wall time} [s]}&{\textbf{Online wall time} [s]} \\
       \midrule
       {Filtered DNS}&{\textbf{-}}&224.9\\
       \midrule
       {noEFR}&{\textbf{-}}&13.5\\
       \midrule
       {Standard EFR/EF}&$\mathcal{O}(10^5)$&213.4\\
       \midrule
       {Smagorinsky}&$\mathcal{O}(10^2)$&20.6\\
       \midrule
       {\ddEF{}}&\multirow{3}{*}{\makecell[l]{3.8 ($\Ttrain=1$, $\Itrain=1$)\\
       42.7 ($\Ttrain=1$, $\Itrain=10$)\\
       12.1 ($\Ttrain=3$, $\Itrain=1$)\\
       145.4 ($\Ttrain=3$, $\Itrain=10$)}}&18.3\\
       \cmidrule{1-1}\cmidrule{3-3}
       {\energyEFR{}}&&25.7\\
       \cmidrule{1-1}\cmidrule{3-3}
       {\enstrophyEFR{}}&&33.9\\
       \bottomrule
    \end{tabular}
    \caption{Offline and online wall times for the strategies considered. The wall time is measured for online simulations with $T=3$ s, at fixed time step $\Delta t=\num{5e-4}$ s.}
    \label{tab:cpu-time}
\end{table}

The online cost is significantly large for standard EFR/EF, since it involves the solution of the differential filter equation. As one would expect, noEFR is the fastest approach, since it does not include any filtering operation or closure terms, followed by \ddEF{} and Smagorinsky.
The addition of the energy and enstrophy constraints in data-driven approaches slightly increases the computational time, but still remains in the same order of magnitude of \ddEF{}. It is important to note that the energy and enstrophy dissipation constraints depend solely on the data from the online simulation, and not on the reference filtered DNS. As a result, the optimization of $\chi(t)$ in \energyEFR{} and \enstrophyEFR{} introduces no significant computational overhead, since it only involves computing the system's energy and enstrophy.

\begin{figure}[htpb!]
    \centering
    \includegraphics[width=\linewidth]{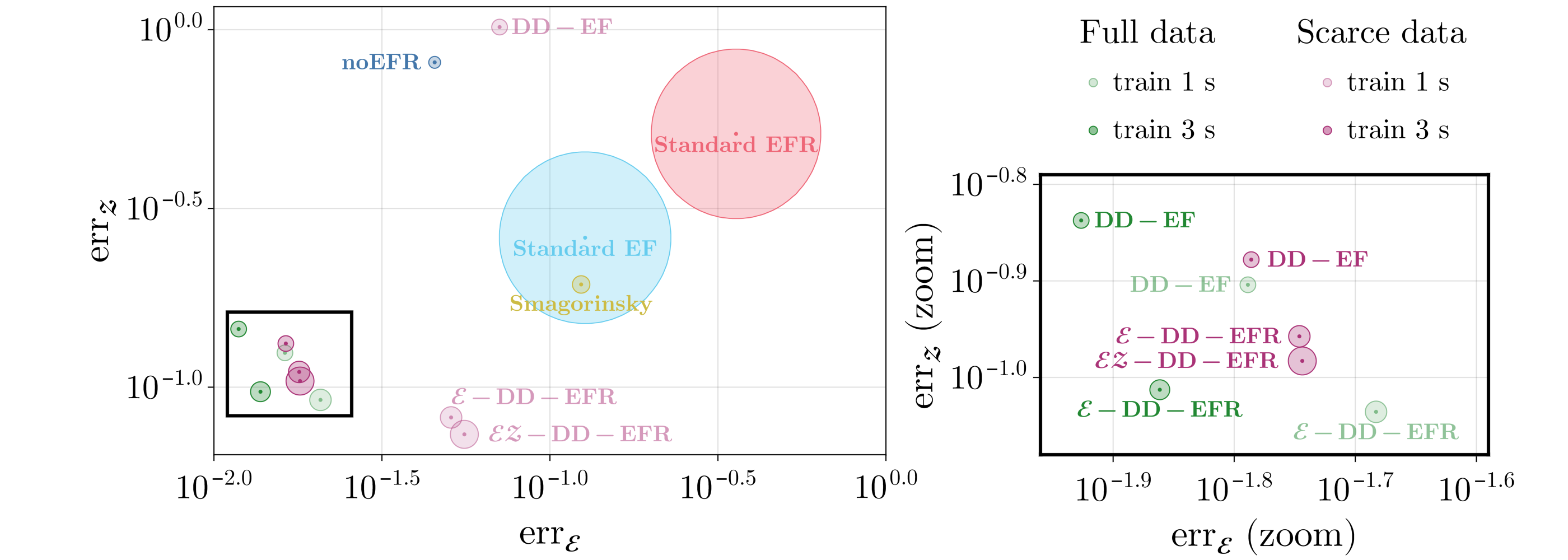}
    \caption{Average energy-enstrophy accuracy of the simulations, where the radius of circles is proportional to the average wall time needed to perform online simulations. The errors are measured in decaying turbulence test case.}
    \label{fig:time-accuracy}
\end{figure}

We include Figure \ref{fig:time-accuracy} as a graphic comparison between the methods. The plot shows the energy and enstrophy errors on the $x$ and $y$ axis, respectively, as presented in \eqref{eq:errs-energy-enstrophy}.
The average error points are surrounded by circles whose size is proportional to the online simulations wall time.
The most accurate methods are located in the bottom-left part of the plots, while the cheapest simulations have smallest circles. Hence, an ideal approach is located in the bottom left part (high accuracy) and is characterized by a small circle (low wall time).

We can see that all DD-EF(R) methods, except for the case $\Ttrain=1$ s and $\Itrain=1$ s (pink circles in the plot), exhibit good accuracy in both energy and enstrophy and are $\sim 10$ times faster with respect to standard EFR and EF, as already formalized in Table \ref{tab:cpu-time}.

The DD-EF(R) simulations in the case of scarce data and $\Ttrain=1$ (pink circles) are less accurate with respect to the other data-driven approaches. However, the results highlight a large improvement in the enstrophy accuracy in the presence of energy and/or enstrophy constraints (\energyEFR{} and \enstrophyEFR{}), while keeping a similar CPU time.

Based on the observed results, we can conclude that the novel DD-EF(R) methods outperform all the existing methods, not only in the accuracy but also in the computational effort.


\section{Conclusion and Outlook}
\label{sec:conclusions}
In this work, we proposed a novel data-driven extension of the Evolve–Filter–Relax (EFR) framework for coarse-grid simulation of turbulent flows.

First, we propose a new learned filter operator, computed through efficient offline least-squares optimization from filtered DNS data, and integrated into a structure-preserving staggered grid finite-volume discretization.

Second, we propose a novel way to choose the relaxation parameter $\chi(t)$, by enforcing global dissipation of energy and/or enstrophy during online simulations. The resulting \energyEFR{} and \enstrophyEFR{} approaches ensure numerical stability and preserve key physical properties of the flow. Specifically, \energyEFR{} enforces global energy dissipation at each time step, while \enstrophyEFR{} extends this approach by additionally enforcing enstrophy dissipation.
These constraints act directly on the online simulation data and help stabilize the learned filtering strategy, making the methods robust, especially when the training data is scarce.

The proposed methodology was tested on two benchmark cases: decaying homogeneous isotropic turbulence and the Kolmogorov flow. In both scenarios, we employed the same filter operator, which was computed solely from the first test case. Despite this, the approach yielded excellent results also in the second case, clearly demonstrating the capability of the learned filter to successfully extrapolate across different flow configurations.

The proposed data-driven method lead to significant improvements over traditional approaches, such as Smagorinsky and standard EFR strategies, particularly in terms of accuracy in energy spectra, global energy, and enstrophy.

Moreover, the methodology maintains low computational cost: the learned filter avoids the need for solving linear systems (as in the traditional differential filter framework), and the adaptive computation of the relaxation parameter relies solely on online simulation data. As such, the method proves both accurate and efficient, especially when compared with classical differential filtering techniques.

Future work will focus on extending the approach to non-periodic boundary conditions, unstructured grids, and three-dimensional flows, as well as integrating alternative learning strategies and uncertainty quantification tools to further enhance robustness and generalization. 

\reviewerALL{Regarding three-dimensional flows, the presented methodology can be directly extended, as the formulation itself is dimension-independent. The main challenge lies in the increased computational cost associated with adding another spatial dimension.

For unstructured or non-uniform grids, the use of FFTs must be revised, as it relies heavily on the Cartesian nature of the grid. In the current implementation, the filter is learned and applied in Fourier space. For unstructured discretizations, a natural alternative is to express the filter in the eigenbasis of the discrete diffusion operator, whose eigenfunctions generalize the Fourier modes. However, projecting onto the full eigenbasis at each time step can be computationally expensive. To alleviate this cost, a reduced spectral basis can be employed, for example obtained via singular value decomposition or by selecting a subset of dominant eigenmodes. This approach significantly reduces the number of operations required to apply the filter, thereby mitigating the added computational burden.

Finally, to handle non-periodic boundary conditions, one could take inspiration from the approach proposed in \cite{rosenberger2023no}, where lifting functions are used to satisfy the boundary constraints. 
}

\section{CRediT authorship contribution}

\textbf{Anna Ivagnes:} Conceptualization, Methodology, Software, Writing - original draft.
\textbf{Toby van Gastelen:} Conceptualization, Methodology, Software, Writing - original draft.
\textbf{Syver D{\o}ving Agdestein:} Software.
\textbf{Benjamin Sanderse:} Conceptualization, Methodology, Writing - review \& editing, Funding acquisition.
\textbf{Giovanni Stabile:} Writing - review \& editing, Funding acquisition.
\textbf{Gianluigi Rozza:} Writing - review \& editing, Funding acquisition.

\section*{Acknowledgements}
This project has been conducted as part of the CWI (Centrum Wiskunde $\&$ Informatica, Amsterdam) Internship Program 2025 and is partially funded by CWI.

AI, GS and GR acknowledges INdAM-GNCS: Istituto Nazionale di Alta Matematica –– Gruppo Nazionale di Calcolo Scientifico.

TvG and BS acknowledge the financial support of Dutch Research Council (NWO) under the project ``Unraveling Neural Networks with Structure-Preserving Computing” (with project number OCENW.GROOT.2019.044
of the research programme NWO XL). In addition, we acknowledge the Eindhoven University of Technology which funded part of this project. 

GS acknowledges the financial support under the National Recovery and Resilience Plan (NRRP), Mission 4, Component 2, Investment 1.1, Call for tender No. 1409 published on 14.9.2022 by the Italian Ministry of University and Research (MUR), funded by the European Union – NextGenerationEU– Project Title ROMEU – CUP P2022FEZS3 - Grant Assignment Decree No. 1379 adopted on 01/09/2023 by the Italian Ministry of Ministry of University and Research (MUR) and acknowledges the financial support by the European Union (ERC, DANTE, GA-101115741). Views and opinions expressed are however those of the author(s) only and do not necessarily reflect those of the European Union or the European Research Council Executive Agency. Neither the European Union nor the granting authority can be held responsible for them.

\newpage

\bibliography{main}
\bibliographystyle{abbrv}

\newpage

\appendix
\section{Appendix}
\label{sec:appendix}
\subsection{Parameter optimization in state-of-the-art approaches}
\label{literature-tuning}

This Section is dedicated to describe additional details about the selection of optimal parameters in well-known state-of-the-art approaches to deal with underresolved turbulence simulations. In particular, we consider:

\begin{itemize}
    \item \textbf{Smagorinsky} model \cite{smagorinsky1963general, lilly1962numerical, lilly1966representation}, which is one of the earliest and most widely used subgrid-scale models in LES of turbulent flows. It belongs to the class of eddy viscosity models and introduces an eddy viscosity given by:
    
    $$\nu_t = \theta^2 \overline{\Delta}^2\sqrt{2\, \text{tr}(\overline{S}\, \overline{S})},$$
    
    where $\overline{\Delta}$ is the filter width, typically related to the grid size, $\overline{S}=\frac{1}{2}(\nabla \overline{\bu} + \nabla \overline{\bu}^T)$ is the strain rate tensor, and $\theta \in [0, 1]$ is the \emph{Smagorinsky coefficient}. Such parameter is typically fitted to the reference data \cite{shankar2023differentiable, guan2023learning}, in our case the filtered DNS.
    \item \textbf{Standard EFR} model, introduced in Section \ref{subsec:efr}, which relies on two additional steps after the \emph{evolve} one (corresponding to discretized Navier--Stokes):
    
    \begin{itemize}
        \item A differential filter step, depending on the filter radius $\delta$. Common choices for $\delta$ are either the minimum mesh size $h_{\text{min}}$ or the Kolmogorov length scale $\eta$ \cite{mou2023energy}.
        \item A relax step depending on a parameter $\chi$. A standard choice is $\chi \sim \Delta t$ \cite{ervin2012numerical}, where $\Delta t$ is the time step of the simulation, but other choices have been explored in \cite{bertagna2016deconvolution, strazzullo2023new, girfoglio2019finite}.
    \end{itemize}
    In this case, we select $\delta=h=\num{7.81e-3}$ and only optimize $\chi$.
    \item \textbf{Standard EF} model, which corresponds to EFR with $\chi=1$. In this case, we optimize the filter radius $\delta$.
\end{itemize}

Each of the above-mentioned models depends on one or more parameters that require careful calibration. Here, we use a gradient-based optimization to determine the optimal values. In particular, we employ a \emph{finite-difference stochastic gradient descent}, where the loss function is derived from $\Itrain=10$ ensemble-averaged simulations, and the optimization seeks to minimize the enstrophy mismatch by adjusting the parameter of interest.

The gradient descent algorithm used for the gradient descent is proposed in \ref{alg:opt}.
In the proposed algorithm, we consider a generic parameter $\alpha$ that needs to be optimized, and we call $\mathcal{M}(\alpha)$ the numerical simulation depending on it. For example, in the EF model, $\mathcal{M}(\alpha)$ is the EF simulation and $\alpha$ coincides with $\delta$.

The loss function is the relative mean squared error (RMSE) in the global enstrophy. More in detail, the enstrophy is evaluated for $I_{\text{train}}=10$ different train simulations which differ only in the random initialization $\boldsymbol{\mu}^{(j)}$ ($j=1, \dots, \Itrain$), and then averaged.

We selected the enstrophy error as objective function because other metrics like the energy error do not provide significant insights on the performance of the method.
The global enstrophy instead is a meaningful metric in fluid dynamics because it quantifies the intensity of vorticity in a flow field, providing direct insight into small-scale structures and dissipative mechanisms.

\begin{figure}
    \centering
    \includegraphics[width=\linewidth]{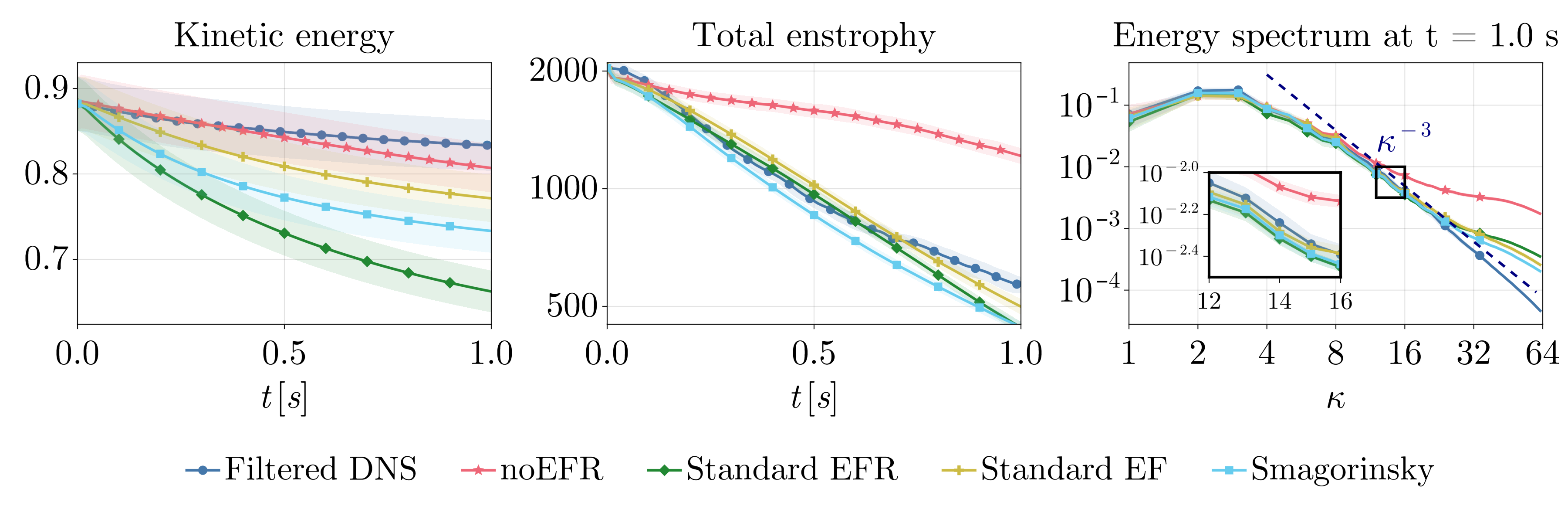}
    \caption{Time evolution for kinetic energy and global enstrophy, spectrum at $t=1$. The plot represents the average and $95 \%$ confidence intervals over $I_{\text{train}}=10$ initializations, in decaying turbulence test case.}
    \label{fig:plots-literature}
\end{figure}

Figure \ref{fig:plots-literature} confirms that noEFR is indeed the least overdiffusive method. In fact, if the kinetic energy error is taken as the objective function, the state-of-the-art approaches effectively reduce to noEFR — that is, setting $\delta=0, \, \chi=0$ in the standard EF/EFR formulations, and $\theta=0$ in the Smagorinsky model.
This observation supports the fairness of the comparison between the enstrophy-based state-of-the-art approaches and the proposed data-driven methods.

\begin{algorithm}[htpb!]
\caption{Gradient descent for optimizing a parameter $\alpha$ in a method of type $\mathcal{M}(\alpha)$.}
\begin{algorithmic}[1]
\Require Initial parameter $\alpha_0$; step size $\beta$; perturbation $\varepsilon$; max iterations $N_{\max}$; tolerance $\eta$; bounds $[\alpha_{\min}, \alpha_{\max}]$; $\boldsymbol{\mu}^{(i)}$, $i=1, \dots, I_{\text{train}}$ random initializations; $n=1, \dots, N_{\text{train}}$ training time steps.
\State Set $\alpha \gets \alpha_0$
\For{$k = 1$ to $N_{\max}$}
    \For{$i = 1$ to $I_{\text{train}}$}
        \State $\mathcal{Z}_h^n(\boldsymbol{\mu}^{(i)}) \gets \mathcal{M}(\alpha, \boldsymbol{\mu}^{(i)})$ \Comment{Run simulation for each initialization}
        \State
        $
        \text{RMSE}^{(i)}(\alpha) = \frac{1}{N_{\text{train}}} \sum_{n=1}^{N_{\text{train}}} \frac{\left( \mathcal{Z}^n_h(\boldsymbol{\mu}^{(i)}) - \mathcal{Z}^n_{\text{ref}}(\boldsymbol{\mu}^{(i)}) \right)^2}{(\mathcal{Z}^n_{\text{ref}}(\boldsymbol{\mu}^{(i)}))^2}$
    \Comment{Compute RMSE in global enstrophy}
    \EndFor
    \State $\mathcal{L}(\alpha) \gets \dfrac{1}{I_{\text{train}}} \sum_{i=1}^{I_{\text{train}}} \text{RMSE}^{(i)}(\alpha)$
    \State $\alpha_{\text{perturbed}} \gets \text{clamp}(\alpha + \varepsilon, \alpha_{\min}, \alpha_{\max})$ \Comment{Restrict parameter to the bounds}
    \For{$i = 1$ to $I_{\text{train}}$}
        \State $\mathcal{Z}_{h,\text{ perturbed}}^n(\boldsymbol{\mu}^{(i)}) \gets \mathcal{M}(\alpha_{\text{perturbed}}, \boldsymbol{\mu}^{(i)})$
        \State $
        \text{RMSE}^{(i)}(\alpha_{\text{perturbed}}) = \frac{1}{N_{\text{train}}} \sum_{n=1}^{N_{\text{train}}} \frac{\left( \mathcal{Z}^n_{h,\text{ perturbed}}(\boldsymbol{\mu}^{(i)}) - \mathcal{Z}^n_{\text{ref}}(\boldsymbol{\mu}^{(i)}) \right)^2}{(\mathcal{Z}^n_{\text{ref}}(\boldsymbol{\mu}^{(i)}))^2}$
    \EndFor
    \State $\mathcal{L}(\alpha_{\text{perturbed}}) \gets \dfrac{1}{I_{\text{train}}} \sum_{i=1}^{I_{\text{train}}} \text{RMSE}^{(i)}(\alpha_{\text{perturbed}})$

    \State $g \gets \dfrac{\mathcal{L}(\alpha_{\text{perturbed}}) - \mathcal{L}(\alpha)}{\varepsilon}$ \Comment{Gradient estimation}
    \State $\alpha_{\text{new}} \gets \text{clamp}(\alpha - \beta g, \alpha_{\min}, \alpha_{\max})$ \Comment{Gradient descent update}
    \If{$|\alpha_{\text{new}} - \alpha| < \eta$}
        \State \textbf{Converged}; \textbf{break}
    \EndIf
    \State $\alpha \gets \alpha_{\text{new}}$
\EndFor
\State \Return $\alpha$
\end{algorithmic}
\label{alg:opt}
\end{algorithm}

\begin{figure}[htpb!]
    \centering
    \subfloat[Smagorinsky model]{\includegraphics[width=0.31\linewidth]{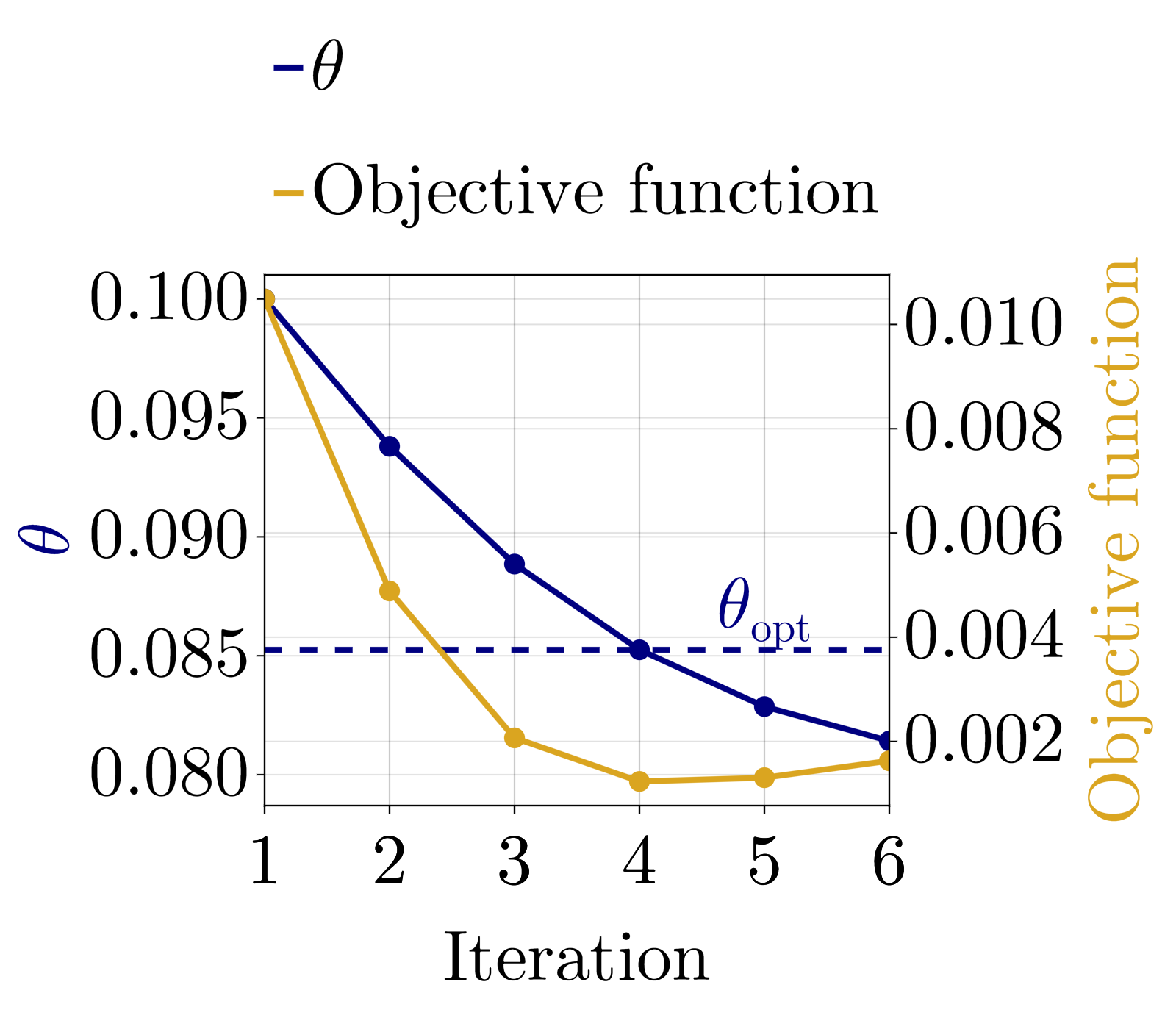}}
    \hspace{0.3cm}
    \subfloat[Standard EF model]{\includegraphics[width=0.31\linewidth]{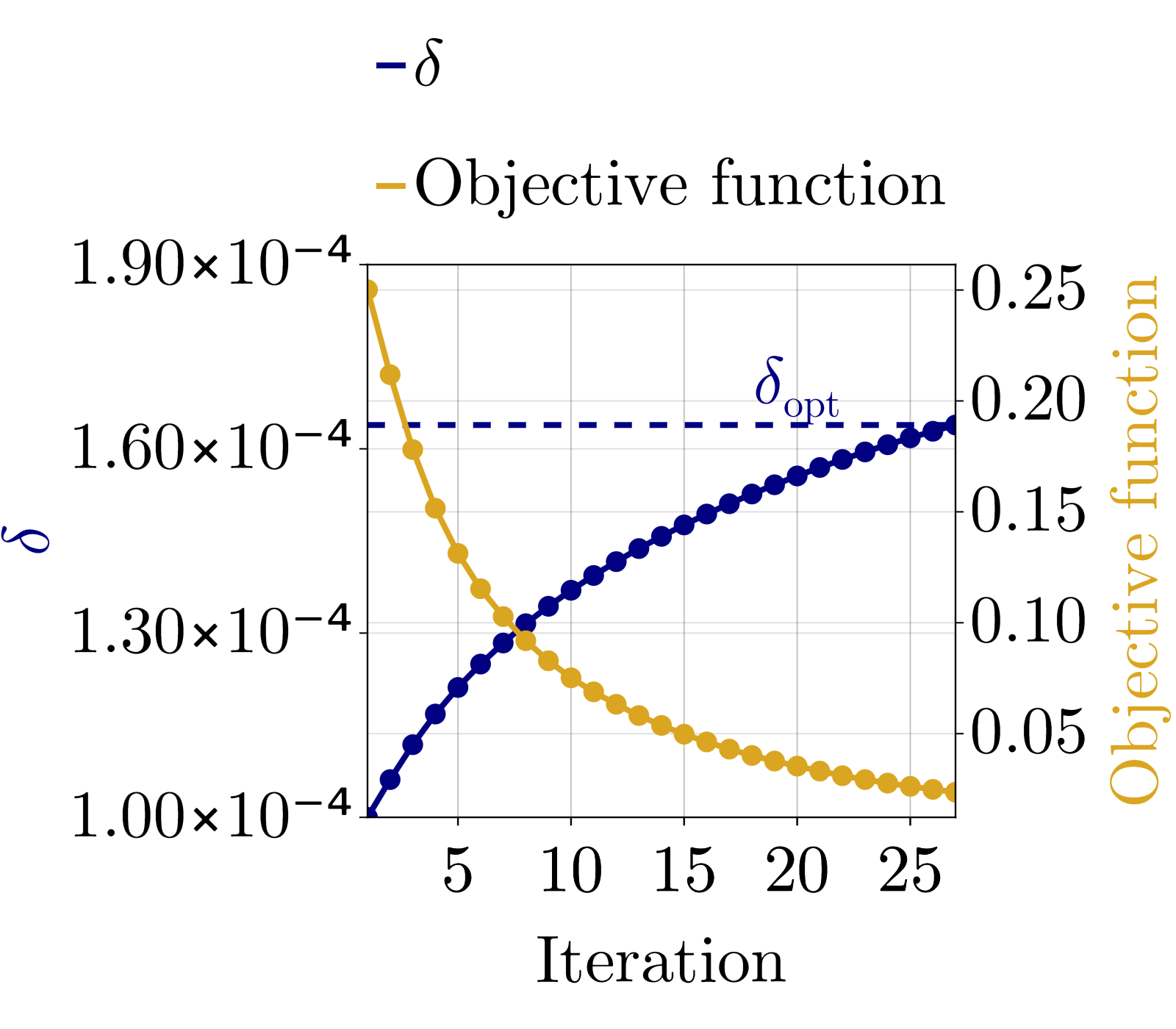} \label{plot:std-ef}}
    \hspace{0.3cm}
    \subfloat[Standard EFR model]{\includegraphics[width=0.31\linewidth]{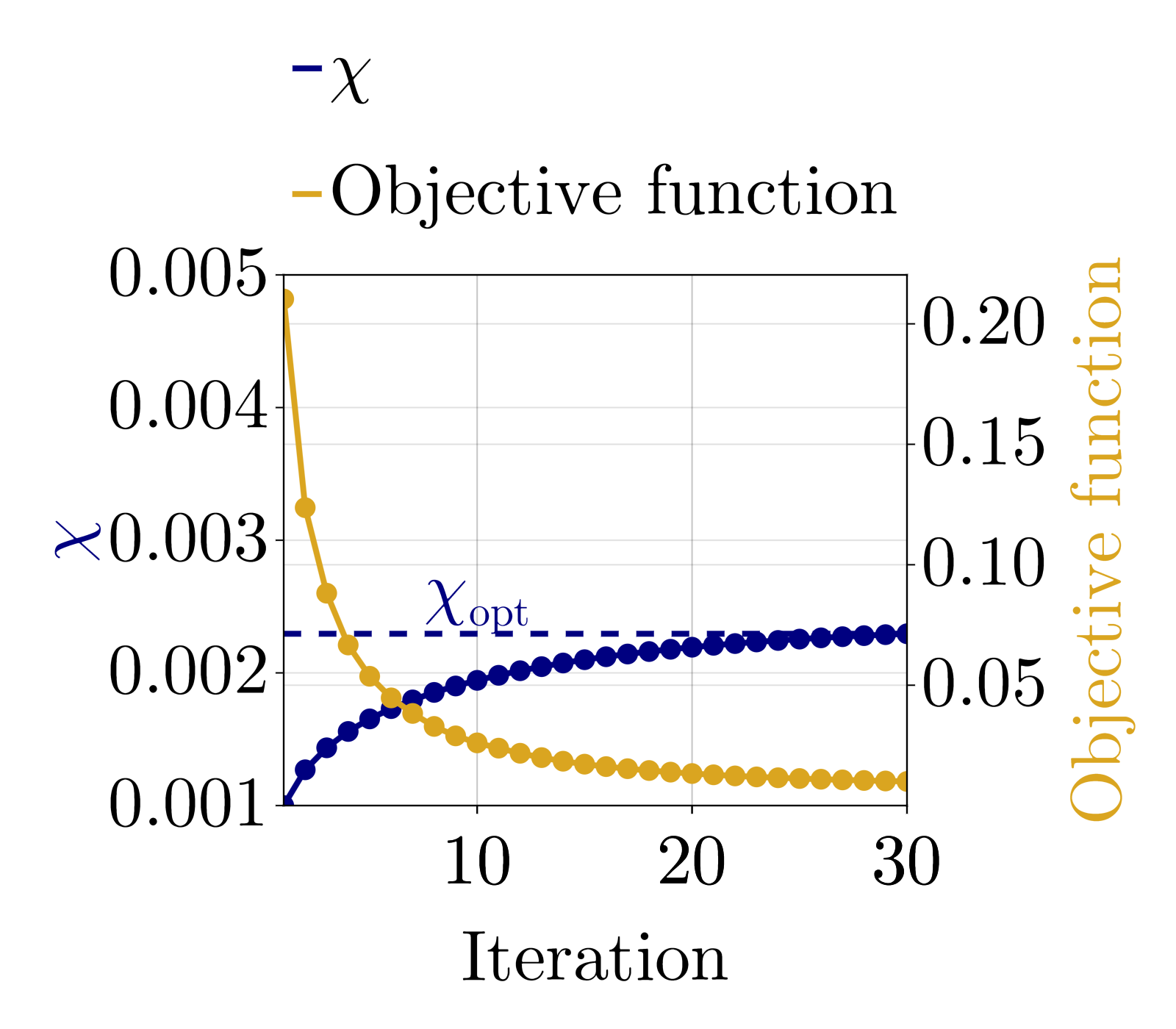}}
    \caption{Parameter tuning of state-of-the-art approaches using approximate gradient descent method.}
    \label{fig:literature-comparison}
\end{figure}

Figure \ref{fig:literature-comparison} shows how the parameters and the objective function are updated during the iterations of the gradient descent.
We emphasize that each state-of-the-art method involves parameters with different orders of magnitude. For a fair comparison, we manually tune the tolerance $\eta$, the step size $\beta$, and the admissible range $[\alpha_{\text{min}}, \alpha_{\text{max}}]$ based on the typical scale of the corresponding quantity. The seemingly slower convergence of standard EF in Figure~\ref{plot:std-ef} is due to the fact that $\delta \sim \mathcal{O}(\num{1e-4})$, whereas in standard EFR we have $\chi \sim \mathcal{O}(\num{1e-3})$, and in the Smagorinsky model $\theta \sim \mathcal{O}(\num{1e-1})$. As a result, the tolerance used in standard EF is slightly higher in relative terms, which can influence the perceived convergence rate.

Table \ref{tab:optimal-literature-parameters} shows an overview of the optimal values obtained for the parameters, and the corresponding admissible range considered in the optimization.

\begin{table}[htpb!]
    \centering
    \begin{tabular}{cccc}
    \toprule
         \textbf{Model} & \textbf{Parameter} & \textbf{Parameter range} & \textbf{Optimal value} \\
         \midrule
         Smagorinsky & $\theta$ &  $[0, 1]$ & $0.085$\\
         \midrule
         Standard EF & $\delta$ & $[0.1 \, \eta, 10\, \eta]$ & $\num{1.64e-4}\simeq 0.5\eta$\\
         \midrule
         Standard EFR & $\chi$ ($\delta=h$)& $[0.1\, \Delta t, 10\, \Delta t]$& $\num{2.3e-3} \simeq 4.5 \, \Delta t$\\
         \bottomrule
    \end{tabular}
    \caption{Optimized parameters in traditional approaches methodologies.}
    \label{tab:optimal-literature-parameters}
\end{table}

Figure \ref{fig:plots-literature} represents the kinetic energy, global enstrophy trends in time, and the energy spectrum at $t=1.$, corresponding to the final time used for the optimization.
The results indicate good agreement in the global enstrophy, as this quantity is directly targeted during the optimization process. However, the energy time history reveals that all methods exhibit excessive dissipation, with noEFR showing the closest match to the filtered DNS.

This motivates the development of our data-driven filtering approach in the following part. 

\subsection{Additional results on data-driven EFR strategies}
\label{supplementary-ddefr}

This part of the supplementary material is dedicated to additional results related to the data-driven approaches presented in the paper.

\subsubsection{Decaying homogeneous turbulence test case}
\label{supplementary-case1}

\begin{figure}[htpb!]
    \centering\includegraphics[width=\linewidth]{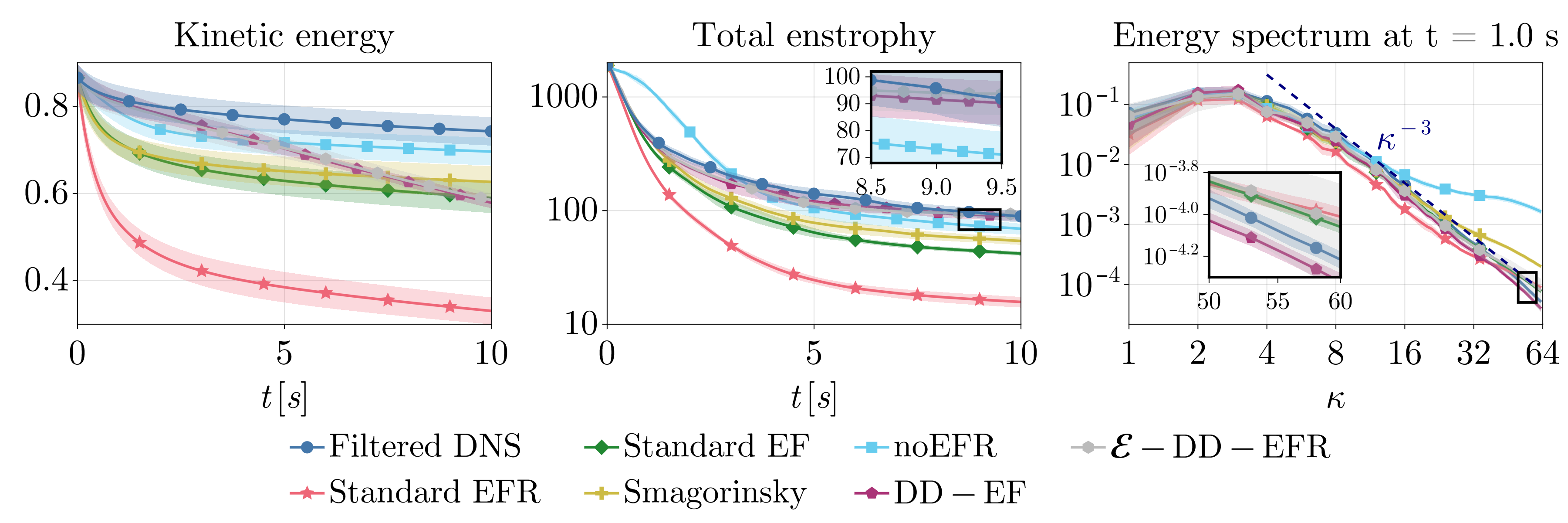}
    \caption{Time evolution of total kinetic energy, total enstrophy, and spectrum of the kinetic energy at time $t=1$ s, for the proposed DD-EF(R) methodologies and traditional approaches.
    In particular, we show the case $\Ttrain=1$ s.
    For each method, the Figure shows the average values (solid lines) and the 95$\%$ confidence interval among 5 test configurations, in decaying turbulence test case.}
    \label{fig:test-general-1s}
\end{figure}

\begin{figure}[htpb!]
    \centering\includegraphics[width=\linewidth]{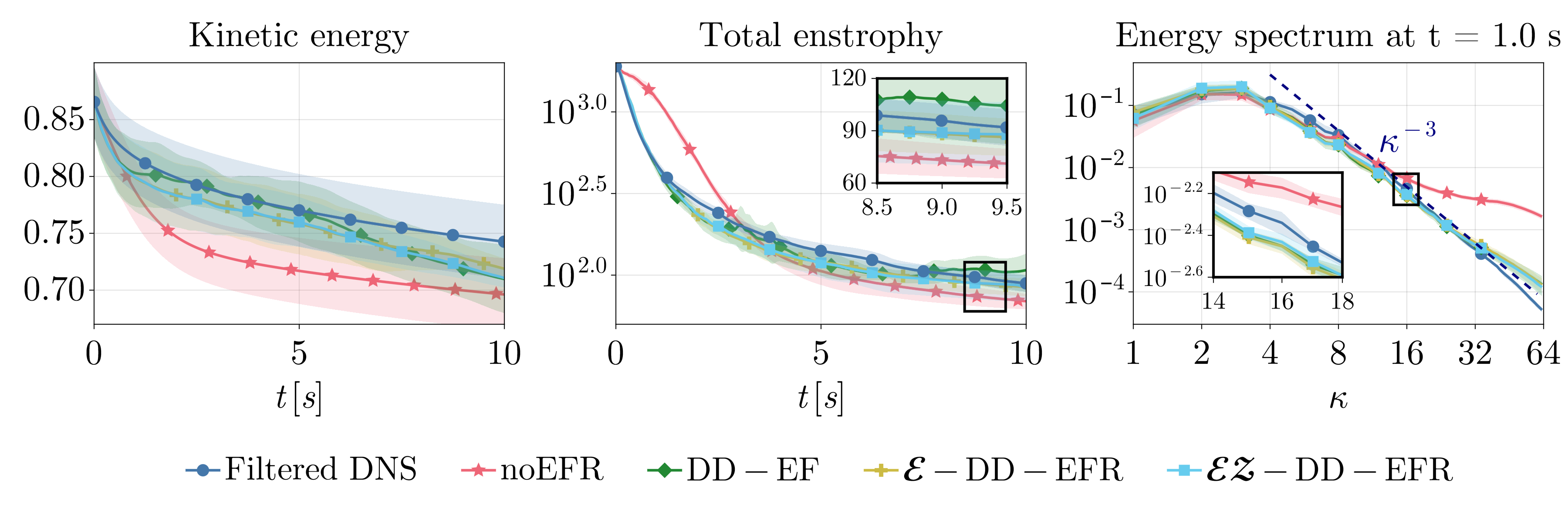}
    \caption{Time evolution of total kinetic energy, total enstrophy, and spectrum of the kinetic energy at time $t=1$ s, for DD-EF(R) methodologies and for state-of-the-art approaches.
    In particular, we show the case  $\Ttrain=3$ s.
    For each method, the Figure shows the average values (solid lines) and the 95$\%$ confidence interval among 5 test configurations, in decaying turbulence test case.}
    \label{fig:scarce-test-general-3s}
\end{figure}

This Section outlines extra results on the data-driven techniques using $\Ttrain=1$ s in the full data regime (Figure \ref{fig:test-general-1s}) and $\Ttrain=3$ s in the scarce data regime (Figure \ref{fig:scarce-test-general-3s}).

In the first case (Figure \ref{fig:test-general-1s}), we can notice good accuracies of \ddEF{} before $t=5$ s, and an overdissipative behaviour for larger times. This is due to the fact the model takes as input less data and has less extrapolation capability.

In the second case (Figure \ref{fig:scarce-test-general-3s}), even if in scarce regime, the DD-EF(R) methods are approximating well the solution also at large times, due to the larger training time. Moreover, \ddEF{} is characterized by some instabilities in the kinetic energy, which are alleviated by the addition of the  energy constraint.

\begin{figure}[htpb!]
    \centering
     \subfloat[\centering \ddEF \\ $\Ttrain=3$ s]{\includegraphics[width=0.14\linewidth]{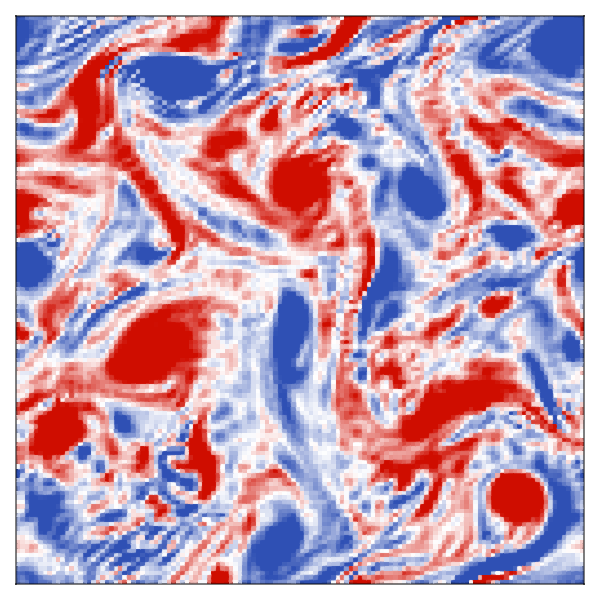}}
     \hspace{0.2cm}
     \subfloat[\centering \energyEFR \\ $\Ttrain=3$ s]{\includegraphics[width=0.14\linewidth]{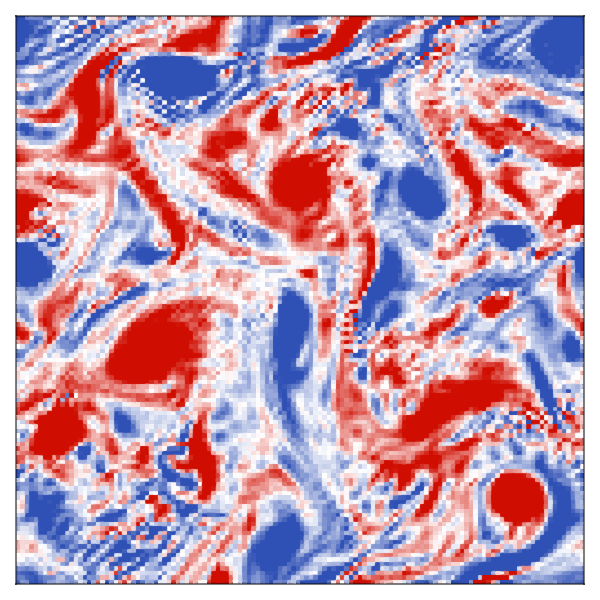}}
     \hspace{0.2cm}
     \subfloat[\centering \enstrophyEFR \\ $\Ttrain=3$ s]{\includegraphics[width=0.14\linewidth]{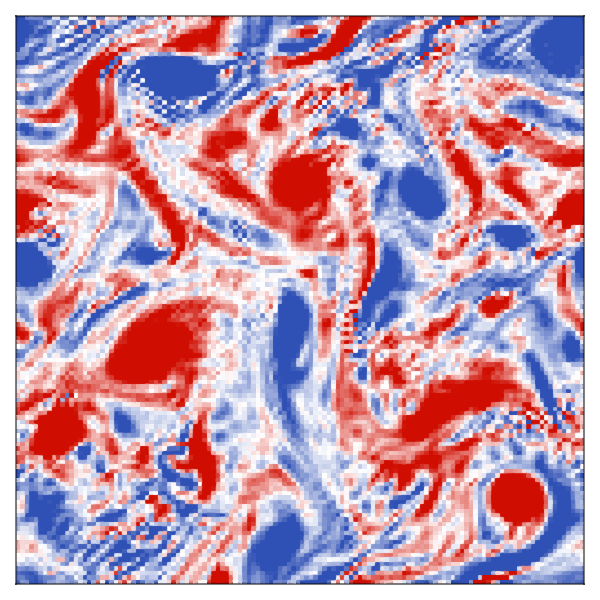}}
     \hspace{0.2cm}
     \subfloat[\centering \ddEF \\ $\Ttrain=1$ s]{\includegraphics[width=0.14\linewidth]{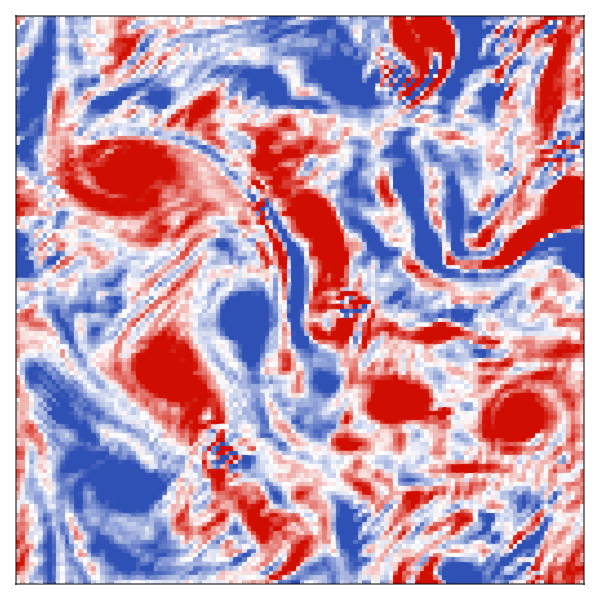}}
     \hspace{0.2cm}
     \subfloat[\centering \energyEFR \\ $\Ttrain=1$ s]{\includegraphics[width=0.14\linewidth]{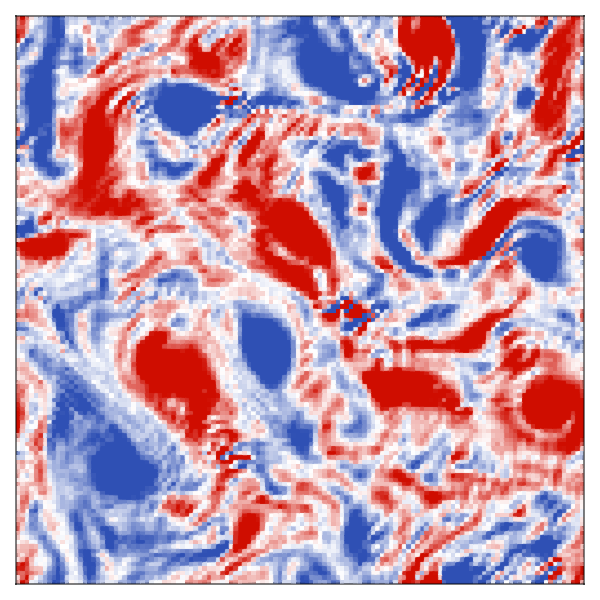}}
     \hspace{0.2cm}
     \subfloat[\centering \enstrophyEFR \\ $\Ttrain=1$ s]{\includegraphics[width=0.14\linewidth]{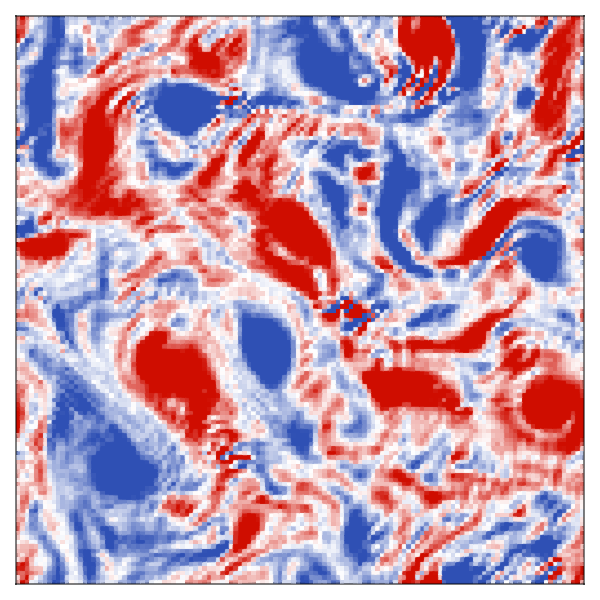}}
     \hspace{0.1cm} \subfloat{\includegraphics[width=0.035\linewidth]{images/decayingturbulence/colorbar.png}}
    \caption{Vorticity fields at $t=1$ s, for different methodologies in a test configuration, in decaying turbulence regime.}
    \label{fig:fields-t1-scarce}
\end{figure}

The solution fields at $t=1$ s for \ddEF{}, \energyEFR{}, and \enstrophyEFR{} in scarce data regimes are represented in Figure \ref{fig:fields-t1-scarce}. The vorticity solutions are qualitatively very similar, and exhibit some spurious oscillations in the case $\Ttrain=1$.

\subsubsection{Kolmogorov test case}
\label{supplementary-case2}

\begin{figure}[htpb!]
    \centering\includegraphics[width=\linewidth]{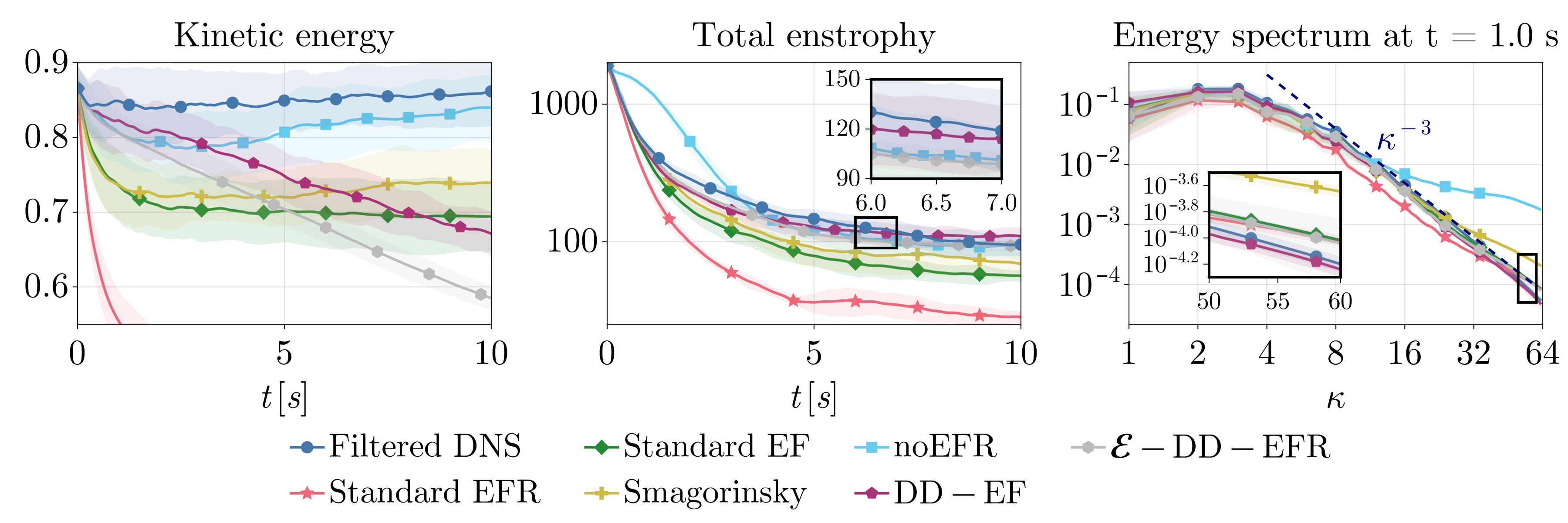}
    \caption{Time evolution of total kinetic energy, total enstrophy, and spectrum of the kinetic energy at time $t=1$ s, for the proposed methodologies and for state-of-the-art approaches.
    In particular, we show the case  $\Ttrain=1$ s, $\Itrain=10$. 
    For each method, the Figure shows the average values (solid lines) and the 95$\%$ confidence interval among 5 test configurations, in the Kolmogorov test case.}
    \label{fig:case2-test-general-1s}
\end{figure}

\begin{figure}[htpb!]
    \centering\includegraphics[width=\linewidth]{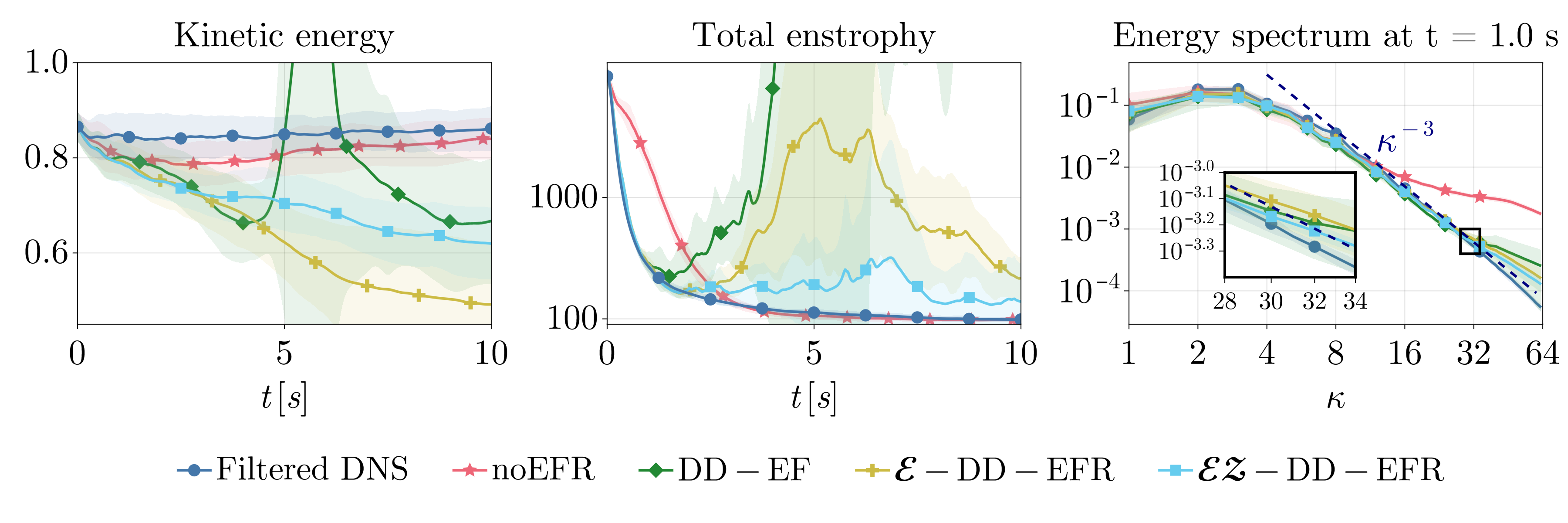}
    \caption{Time evolution of total kinetic energy, total enstrophy, and spectrum of the kinetic energy at time $t=1$ s, for the proposed methodologies  and for state-of-the-art approaches.
    In particular, we show the case $\Ttrain=1$ s, $\Itrain=1$ (scarce data). 
    For each method, the Figure shows the average values (solid lines) and the 95$\%$ confidence interval among 5 test configurations, in the Kolmogorov test case.}
    \label{fig:case2-scarce-test-general-1s}
\end{figure}

This Section is dedicated to extra results in the Kolmogorov test case, specifically in the case $\Ttrain=1$, both for full and scarce data (Figures \ref{fig:case2-test-general-1s} and \ref{fig:case2-scarce-test-general-1s}, respectively).

As noticed in decaying turbulence case, when $\Ttrain$ is smaller, the \ddEF{} result either in stable but overdiffusive solutions (full data) or in unstable solutions (scarce data). The energy and enstrophy constraints alleviate such oscillations and instabilities.

\begin{figure}[htpb!]
    \centering
    \subfloat[Filtered DNS]{\includegraphics[width=0.14\linewidth]{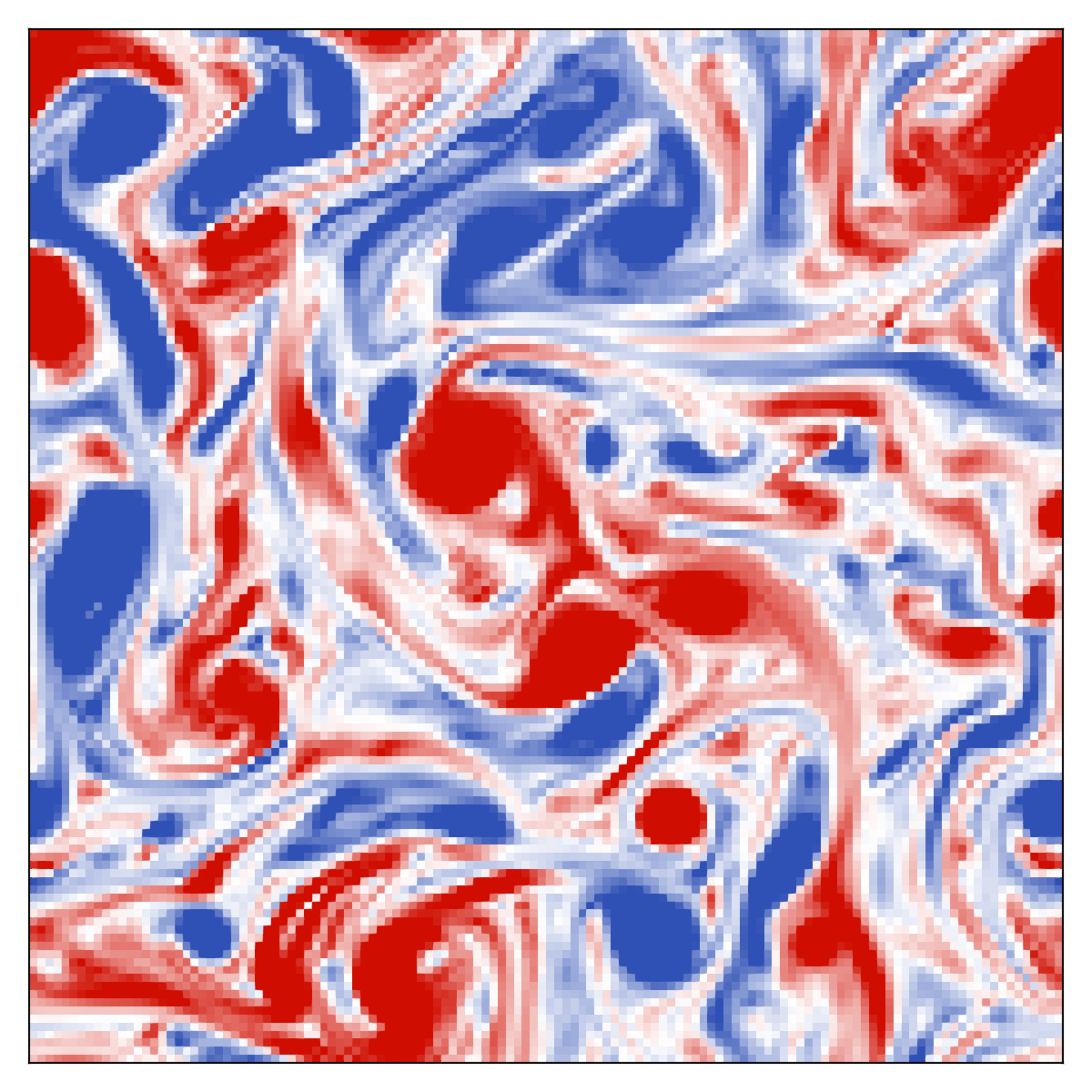}}
    \hspace{0.5cm}
    \subfloat[noEFR]{\includegraphics[width=0.14\linewidth]{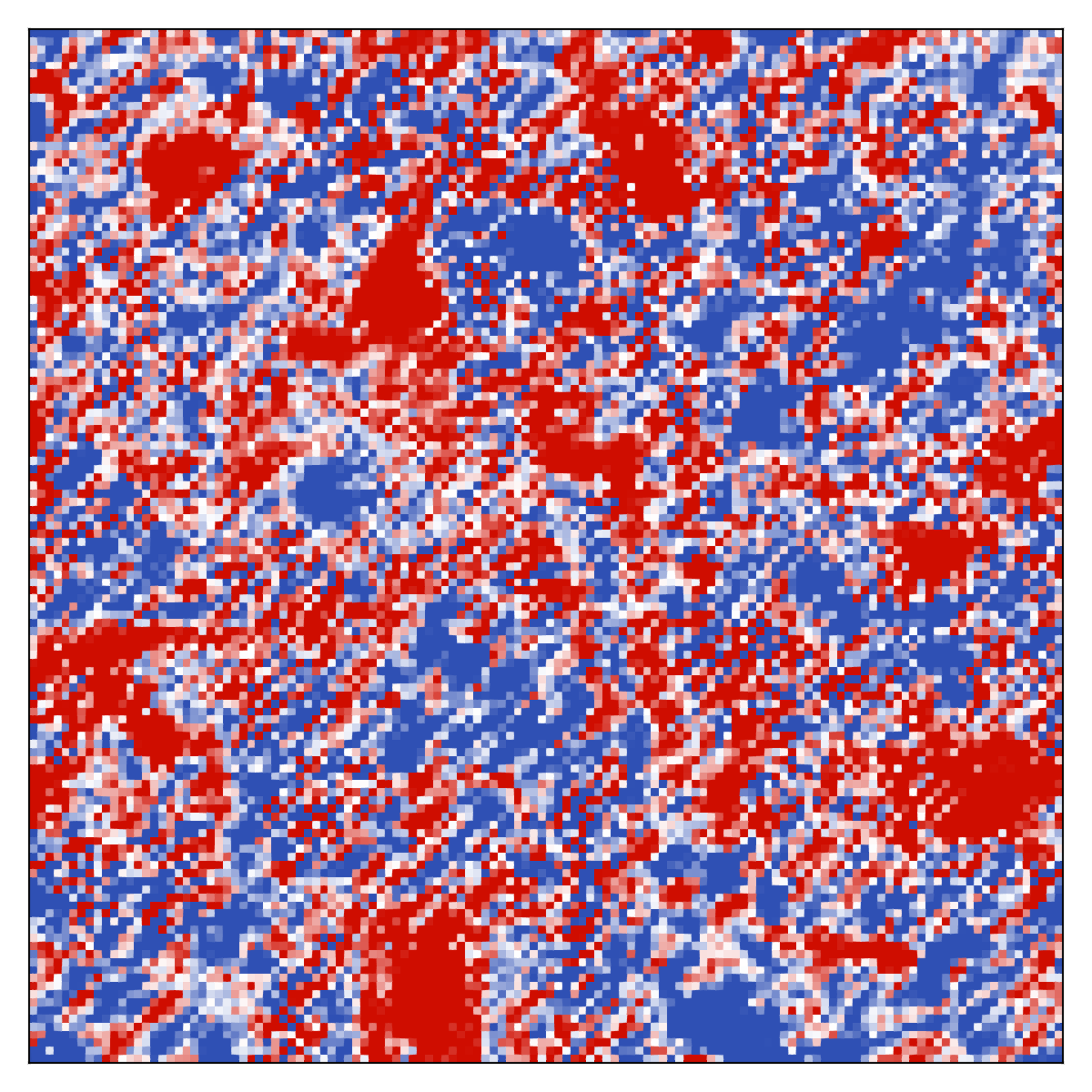}
    }
   \hspace{0.5cm}
    \subfloat[Smagorinsky]{\includegraphics[width=0.14\linewidth]{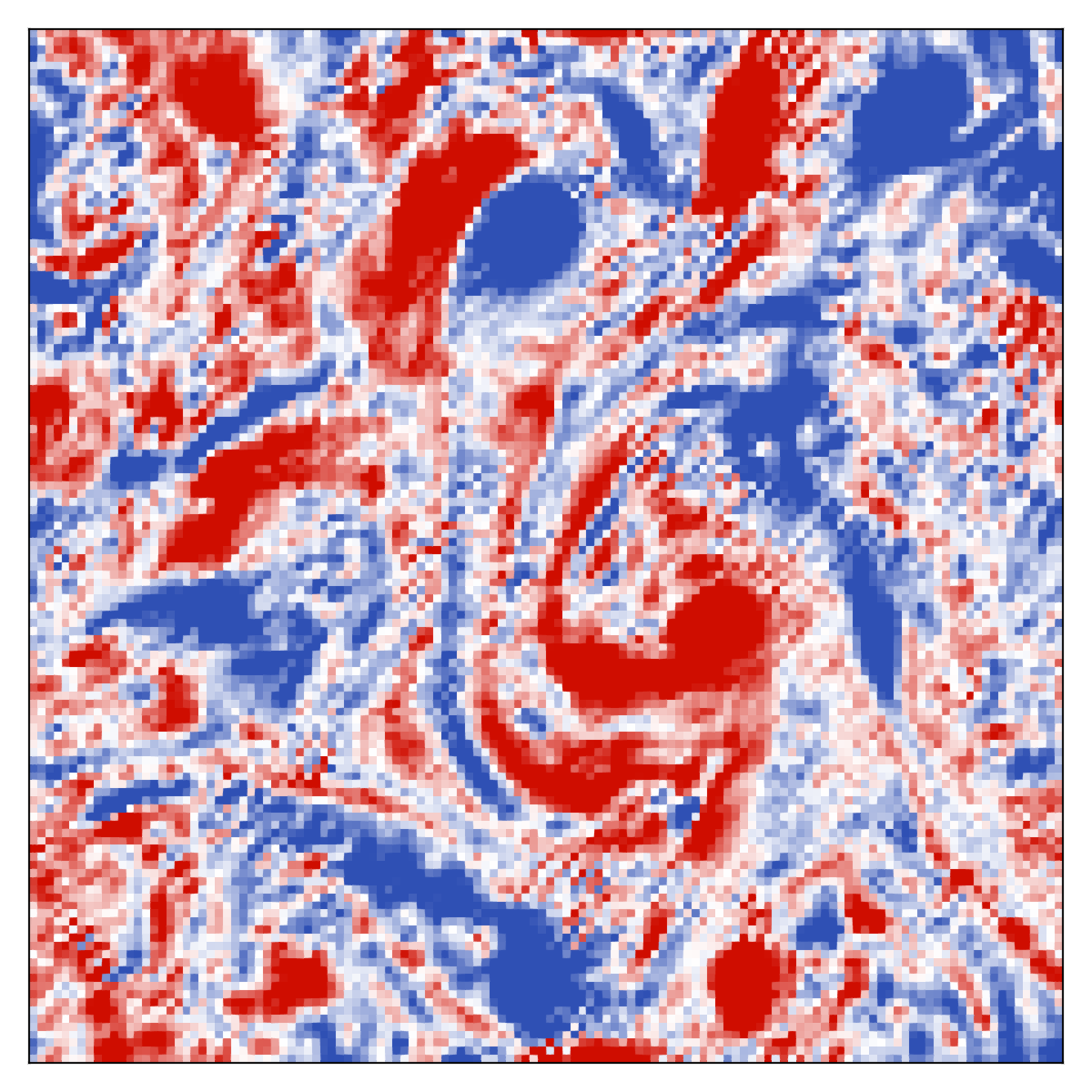}}
    \hspace{0.5cm}
     \subfloat[Standard EFR]{\includegraphics[width=0.14\linewidth]{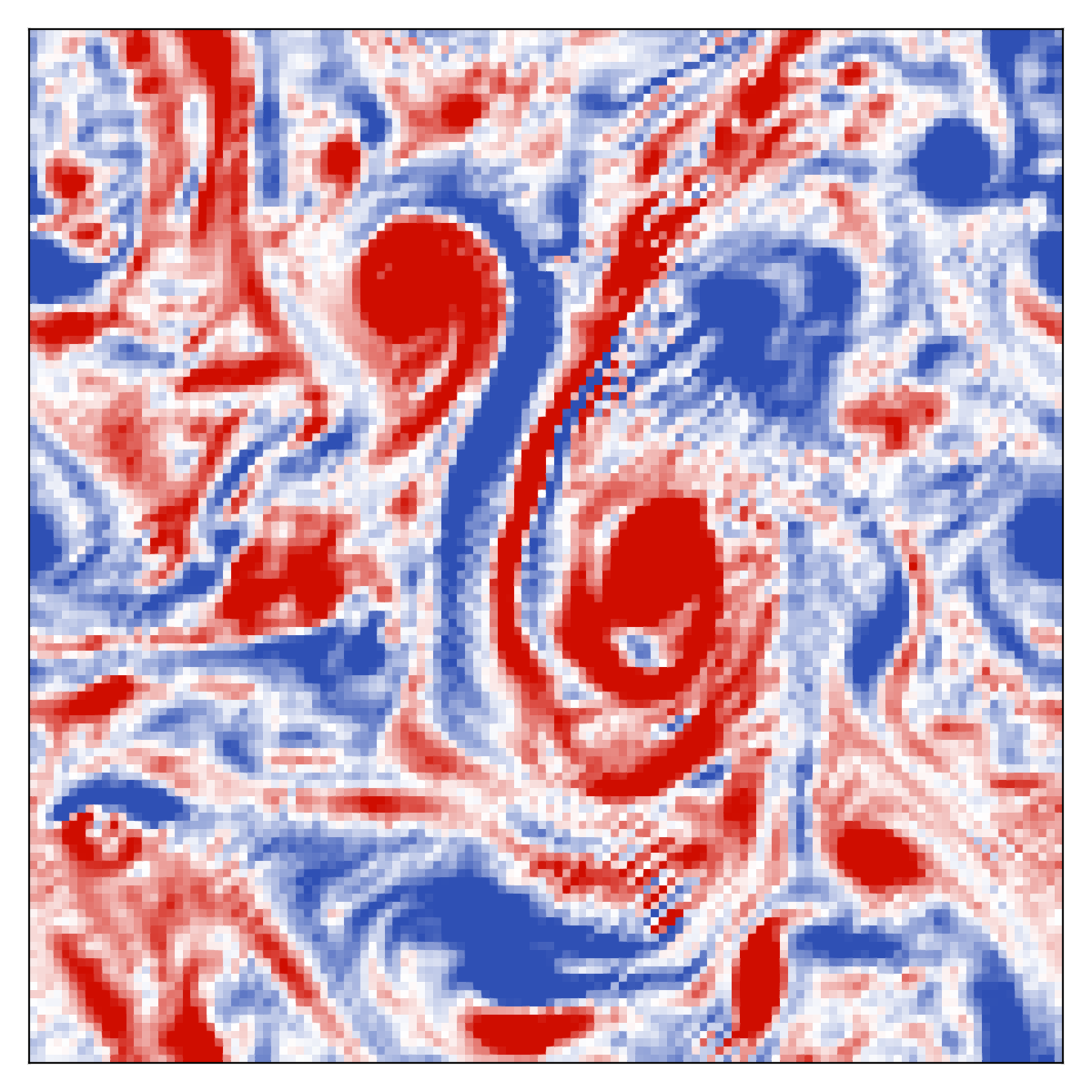}}
     \hspace{0.5cm}
     \subfloat[Standard EF]{\includegraphics[width=0.14\linewidth]{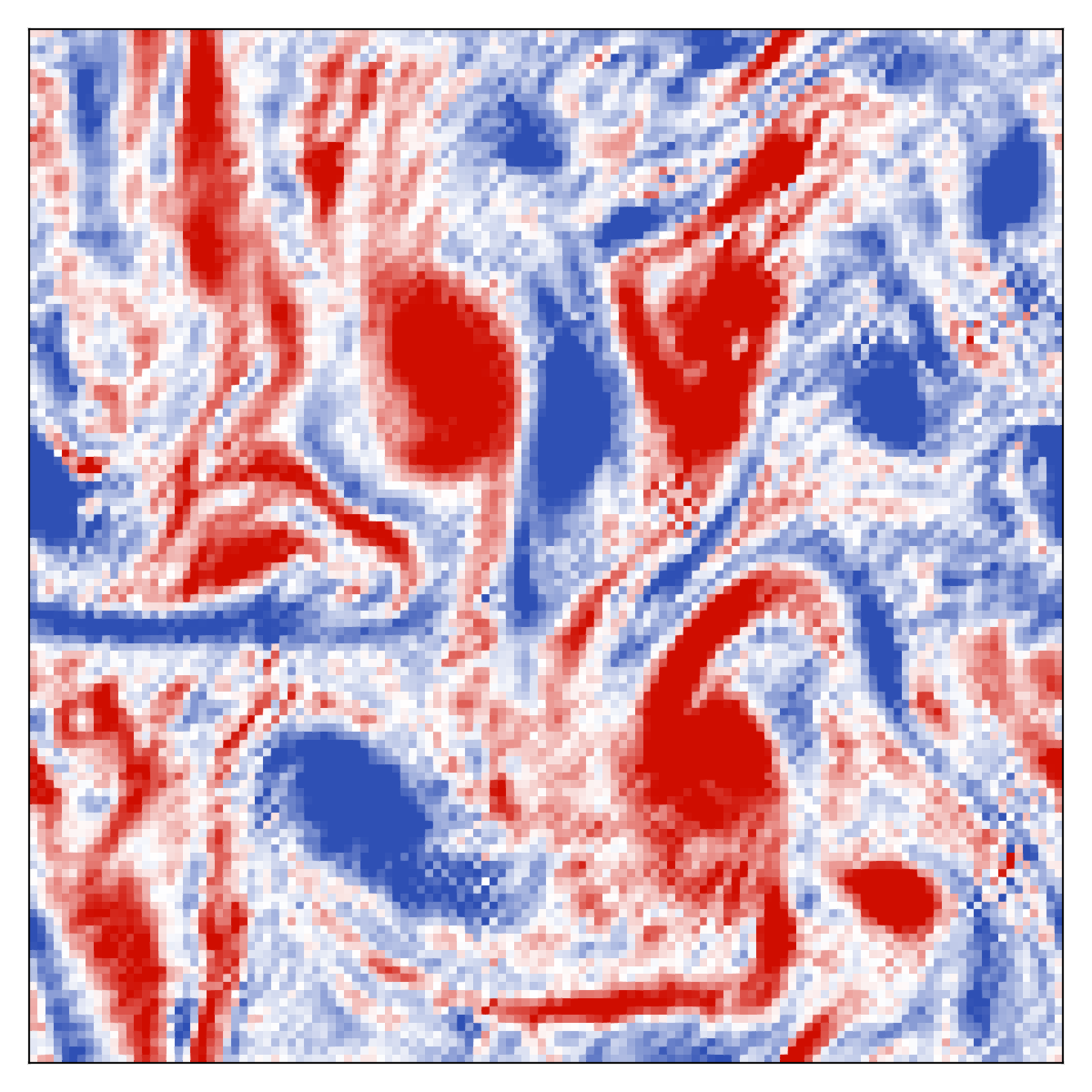}}
     \hspace{0.5cm}
     \subfloat{\includegraphics[width=0.035\linewidth]{images/decayingturbulence/colorbar.png}}\\
     \hspace{0.5cm}
     \subfloat[\centering \ddEF \\ $\Ttrain=3$]{\includegraphics[width=0.14\linewidth]{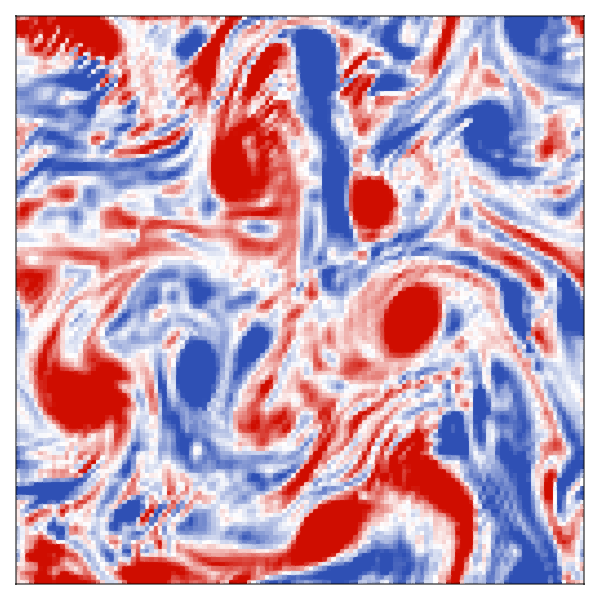}}
     \hspace{0.5cm}
     \subfloat[\centering \energyEFR \\ $\Ttrain=3$]{\includegraphics[width=0.14\linewidth]{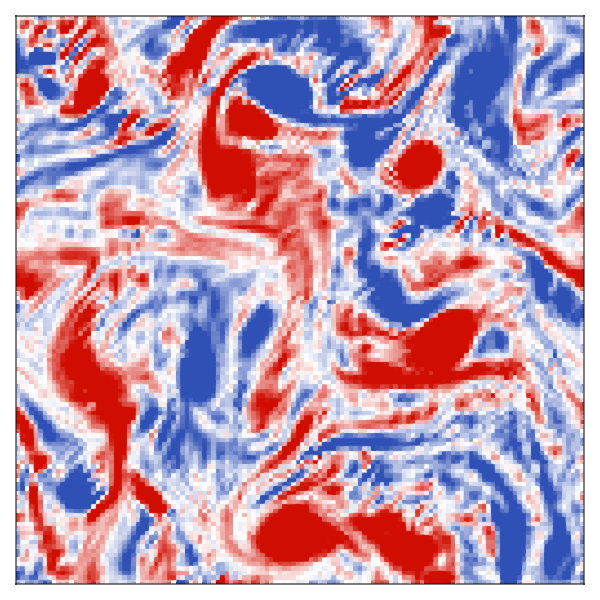}}
     \hspace{0.5cm}
     \subfloat[\centering \ddEF \\ $\Ttrain=1$]{\includegraphics[width=0.14\linewidth]{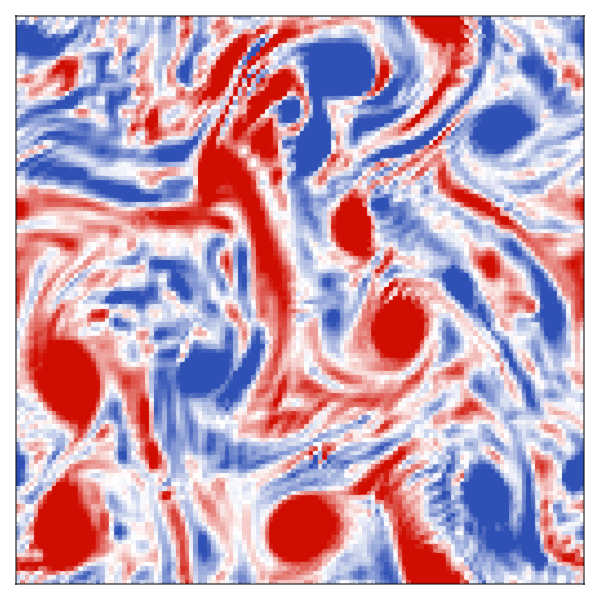}}
     \hspace{0.5cm}
     \subfloat[\centering \energyEFR \\ $\Ttrain=1$]{\includegraphics[width=0.14\linewidth]{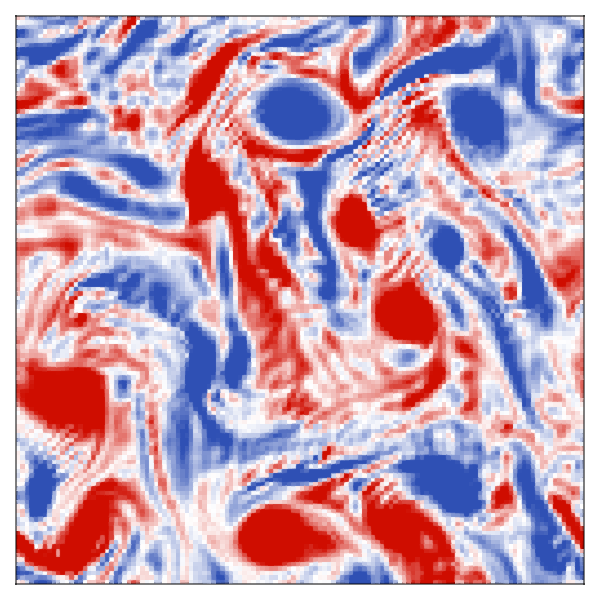}}
     \hspace{0.5cm}
     \subfloat{\includegraphics[width=0.035\linewidth]{images/decayingturbulence/colorbar.png}}
    \caption{Vorticity fields at $t=1$ s, for different methodologies in a test configuration.}
    \label{fig:case2-fields-t1}
\end{figure}

The solutions at fixed time instance $t=1$ s are represented in Figure \ref{fig:case2-fields-t1}.
\end{document}